\documentclass[a4paper,10pt]{article}

\title{\large \bf $(a,b)$-Koszul algebras}
\author{Andrea Rey - Andrea Solotar}
\date{}

\usepackage[all]{xy}
\usepackage{amsmath,amsfonts,amssymb}
\usepackage[fleqn,tbtags]{mathtools}
\usepackage{color}

\usepackage{listings}
\lstnewenvironment{listing}[1][]
   {\lstset{#1}\pagebreak[0]}{\pagebreak[0]}
   
\newtheorem{thm}{\bf Theorem}[section]
\newtheorem{prop}[thm]{\bf Proposition}
\newtheorem{lema}[thm]{\bf Lemma}
\newtheorem{rmk}[thm]{\bf Remark}
\newtheorem{defi}[thm]{\bf Definition}
\newtheorem{coro}[thm]{\bf Corollary}
\newtheorem{example}[thm]{\bf Example}
\newenvironment{prof}{\noindent \bf Proof. \rm}

\numberwithin{equation}{section}

\def\qed{\hfill \mbox{$\square$}}
\def\place{{-}}

\def\ker{\mathop{\rm Ker}\nolimits}

\def\im{\mathop{\rm Im}\nolimits}

\def\Mod{\mathop{\rm Mod}\nolimits}
\def\gr{\mathop{\rm deg}\nolimits}
\def\Hom{\mathop{\rm Hom}\nolimits}
\def\hom{\mathop{\rm hom}\nolimits}
\def\tor{\mathop{\rm Tor}\nolimits}
\def\ext{\mathop{\rm ext}\nolimits}
\def\Ext{\mathop{\rm Ext}\nolimits}

\usepackage{hyperref}
\hypersetup{
colorlinks,
citecolor=blue,
linkcolor=black}

\begin{document}

\maketitle

\begin{abstract}
Let $a$ and $b$ be two integers such that $2\le a<b$. In this article we define the notion of $(a,b)$-Koszul algebra as a generalization of $N$-Koszul algebras. We also exhibit examples and we provide a minimal graded projective resolution of the algebra $A$ considered as $A$-bimodule, which allows us to compute the Hochschild homology groups for some examples of $(a,b)$-Koszul algebras.
\end{abstract}

\medskip

\textbf{Mathematical subject codes:} 16E05, 16E30, 16E40, 16S37, 16W50.

\tableofcontents


\section{\texorpdfstring{Introduction}{sec:introduction}}

The definition of Koszul algebras was given by S. Priddy \cite{P} in 1970, motivated by an article published by Koszul in the 50s \cite{K}. The study of these algebras and their generalizations has been developed in the last years because of their applications in algebraic geometry, Lie theory, quantum groups, algebraic topology and combinatorics (see for instance \cite{BGS2,F,HL}). Berger defined the notion of $N$-Koszul algebra in his article \cite{B2}, see also \cite{BG}. The condition of being $N$-Koszul is also described in \cite{GMMZ}, where the authors study quotients of tensor algebras over copies of the base field by an ideal generated by homogeneous elements of degree $N$. Moreover, similarities among the notions of Koszul algebras and $N$-Koszul algebras have been proved (see for example \cite{B2, GMMZ, BG, FV}). However, the $N$-homogeneity condition of the relations is the cause of various problems. In particular, the class of $N$-Koszul algebras, for $N>2$, is not closed under graded \"Ore extensions, normal regular extensions or tensor products. It is then natural to generalize the definition of $N$-Koszul algebras to algebras $A=T(V)/I$ where $V$ is a finite dimensional vector space and $I$ is an ideal generated by homogeneous elements of two different degrees. 

The contents of the article are as follows. We start by recalling in Section \S 2 well-known definitions and results about graded modules and we prove some technical lemmas. 

In Section \S 3 we construct a minimal graded projective resolution of the trivial $A$-module $k$ for $A$ an $(a,b)$-homogeneous algebra and we study in detail the kernels of the morphisms of this resolution. 

Section \S 4 is devoted to the definition of $(a,b)$-Koszul algebras and to find equivalent conditions to this definition. Some of these conditions are related to distributive lattices. We generalize the notion of distributive triples to multidistributive tuples. We define the opposite algebra $A^{\circ}$ and the Koszul dual $A^{!}$ of an $(a,b)$-homogeneous algebra $A$, next we prove that $A^{\circ}$ is $(a,b)$-Koszul if and only if $A$ is. However, this is not the case for $A^{!}$. 

In \cite{CS}, the authors define the algebras $\mathcal{K}_2$ as a possible generalization of $N$-Koszul algebras. The $(a,b)$-Koszul algebras are related to the $\mathcal{K}_2$ algebras (see \cite{R}).

Notice also that the notion of $(a,b)$-Koszul algebra is different from the $(p,q)$-Koszul rings in \cite{BBK}.

The aim of Section \S 5 is to give several examples of $(a,b)$-Koszul algebras.

In Section \S 6 we construct a minimal $A$-bimodule resolution of an $(a,b)$-Koszul algebra. This resolution allows the computation of the Hochschild homology groups of the algebra; in fact, we compute the dimensions of the $k$-vector spaces $HH_i(A)$ for some $(a,b)$-Koszul algebras $A$. 

Finally, in Section \S 7 we exhibit an algorithm using \cite{WM} needed for the computations of the previous section.

Throughout this article, $k$ will denote a field, $V$ a finite dimensional $k$-vector space, $A$ an augmented associative $\mathbb{Z}$-graded $k$-algebra with unit and $I$ a two-sided ideal of $A$ generated by homogeneous elements of degrees $a$ and $b$ such that $2\le a<b$. We will denote the vector space $V^{\otimes n}$ by $V^{(n)}$ and an elementary tensor $v_1\otimes \cdots \otimes v_n\in V^{(n)}$ by $v_1\cdots v_n$. By ``module'' we will mean a $\mathbb{Z}$-graded left $A$-module, unless the contrary is stated.


\section{\texorpdfstring{Preliminaries and basic properties}{sec:preliminaries}}\label{sec:preliminaries}

Given an algebra $A=\bigoplus_{m\in \mathbb{Z}} A_m$, a $\mathbb{Z}$-graded left $A$-module $M=\bigoplus_{n\in \mathbb{Z}} M_n$ is a $\mathbb{Z}$-graded $k$-vector space such that $A_mM_n\subseteq M_{n+m}$. We will denote by $A$-grMod the abelian category of $\mathbb{Z}$-graded left $A$-modules, where the morphisms are the $A$-linear maps preserving the grading. The $A$-module $M$ is said to be \textbf{left bounded} if there exists an integer $m$ such that $M_n=0$ for all $n<m$. The graded left $A$-modules which are left bounded form a full subcategory of $A$-grMod.

Suppose that there exists a morphism of graded $k$-algebras $\varepsilon :A\rightarrow k$ such that $\varepsilon(1)=1$. Then, it induces a left $A$-module structure on $k$ such that $a\lambda :=\varepsilon (a)\lambda$ with $a\in A$ and $\lambda \in k$.

Given $l\in \mathbb{Z}$ and $M\in$ $A$-grMod, its shift $M[l]$ is defined by $(M[l])_n=M_{n+l}$ for all $n\in \mathbb{Z}$. The $A$-module $M$ is said to be \textbf{graded-free} if it has a basis of homogeneous elements. If $M$ is left bounded, then it is graded-free if and only if it is isomorphic to a direct sum of shifts $A[-l_i]$, where the subset of $\mathbb{Z}$ formed by the degrees $l_i$ is left bounded.

We next recall without proof some well-known results concerning the category $A$-grMod.

\begin{prop}\cite{B2, C}\label{prop:gradedfreeiifproj}
An object $M$ in $A$-grMod is projective if and only if it is graded-free.
\end{prop}

\begin{defi}\cite{B4, C}
A surjective morphism $f:M\rightarrow M'$ in $A$-grMod is called \textbf{essential} if for each morphism $g:N\rightarrow M$ in $A$-grMod such that $f\circ g$ is surjective, then $g$ is also surjective.
\end{defi}

\begin{defi}\cite{B4}
Let $M$ be an object in $A$-grMod. A \textbf{projective cover} of $M$ is a pair $(P,f)$ such that $P\in$ $A$-grMod is projective and $f:P\rightarrow M$ is an essential surjective morphism.
\end{defi}

\begin{prop}\cite{C}\label{prop:projcover}
Every $\mathbb{Z}$-graded $A$-module $M$ has a projective cover, which is unique up to isomorphism.
\end{prop}

\begin{prop} \cite[Graded version of the Nakayama Lemma]{B2,C}\label{prop:trivialmod}
Let $M$ be a $\mathbb{Z}$-graded left $A$-module such that $M_n=0$ for all $n<0$. If $k\otimes_A M=0$ then $M$ is also zero.
\end{prop}

We introduce the following definition that will be used throughout this article.

\begin{defi}\label{def:sconcentrated}
An object $M$ in $A$-grMod is called \textbf{$s$-concentrated} (respectively \textbf{$s$-pure}) in degrees $l_1, \cdots ,$ $l_s$ if there exist non-negative integers $l_1< \cdots <l_s$ with $M=M_{l_1}\oplus \cdots \oplus M_{l_s}$ (respectively $M=AM_{l_1}+ \cdots + AM_{l_s}$).
\end{defi}

\begin{rmk}
If $s=1$ we simply say that $M$ is a concentrated (respectively pure) module (cf. \cite{B2}).
\end{rmk}

\begin{rmk}
In \cite{GMMZ} the authors say that a graded module $M$ which is pure in degrees $l_1, \cdots ,l_s$ is \textbf{generated in degrees $l_1,\cdots ,l_s$}.
\end{rmk}

In both cases, the integers $l_1, \cdots ,l_s$ such that $M_{l_1},\cdots ,M_{l_s}$ are nonzero are uniquely determined whenever $M$ is a nontrivial module. It is evident that every module which is $s$-concentrated in degrees $l_1,\cdots ,l_s$ is $s$-pure in degrees $l_1,\cdots ,l_s$. Moreover, every module $s$-concentrated in degrees $l_1,\cdots ,l_s$ is isomorphic to a direct sum of shifts $k[-l_1]\oplus \cdots \oplus k[-l_s]$ that will be denoted $k[-l_1,\cdots , -l_s]$. On the other hand, given $M$ in $A$-grMod which is $s$-pure in degrees $l_1,\cdots ,l_s$, it is isomorphic to a direct sum of shifts $A[-l_1,\cdots ,-l_s]$ and so isomorphic to $(A\otimes_k M_{l_1})\oplus \cdots \oplus (A\otimes_k M_{l_s})$ where $M_{l_i}$ is concentrated in degree $l_i$.

We note that simple graded modules are isomorphic to $k[-l]$, so every $s$-concentrated module is semisimple.

\begin{prop}\label{prop:essentialineachdegree}
Let $f\in \Hom_{A\text{-grMod}}(M,M')$ be surjective. If $M$ is $s$-pure in degrees $l_1,\cdots ,l_s$, then $M'$ is also $s$-pure in degrees $l_1,\cdots ,l_s$. Moreover, $f$ is essential if and only if the induced morphisms $f_{l_i}:M_{l_i}\rightarrow M'_{l_i}$ are bijective for $1\le i \le s$.
\end{prop}

\begin{prof}
Let $m'\in M'$ and $m\in M$ be such that $f(m)=m'$. Since $M$ is $s$-pure in degrees $l_1,\cdots ,l_s$, there exist $a_i\in A$ and $m_i\in M_{l_i}$ for $1\le i \le s$ such that $m=a_1m_1+\cdots +a_sm_s$. Therefore, $m'=a_1f(m_1)+\cdots +a_sf(m_s)$ where $f(m_i)\in M'_{l_i}$. Thus, $M'$ is $s$-pure in degrees $l_1,\cdots ,l_s$.

The Nakayama Lemma in the graded category says that $f$ is essential if and only if the induced $k$-linear map
\[
\overline{f}:k\otimes_A M\longrightarrow k\otimes_A M'
\]

\noindent is bijective. However, since $M$ and $M'$ are $s$-pure then $k\otimes_A M$ and $k\otimes_A M'$ are canonically isomorphic to $k\otimes M_{l_1}+ \cdots + k\otimes M_{l_s}$ and $k\otimes M'_{l_1}+ \cdots +k\otimes M'_{l_s}$ respectively. Thus, the restrictions to each degree of $\overline{f}$ become $f_{l_i}:M_{l_i}\rightarrow M'_{l_i}$ for $1\le i\le s$. \qed 
\end{prof}

\bigskip

Next, we state three easy technical lemmas concerning $k$-vector spaces that will be very useful in the sequel.

\begin{lema}\label{lem:trivialspace}
Let $V$ be a nonzero finite dimensional $k$-vector space and let $W\subseteq V^{(n)}$ be a subspace. If there exists $m>0$ such that $V^{(m)}\otimes W=0$, then $W=0$.
\end{lema}

\begin{prof}
Let $\dim V=s\ne 0$. Then, $\dim V^{(m)}=ms$ and $\dim (V^{(m)}\otimes W)=ms\dim W$. If there exists $m>0$ such that $V^{(m)}\otimes W=0$ then $ms\dim W=0$. Hence, $\dim W=0$. \qed
\end{prof}

\bigskip

From now on we fix a basis $\mathcal{B}=\{ v_1,\cdots ,v_s\}$ of the $k$-vector space $V$.

\begin{lema}\label{lem:inthesubsp}
Let $W$ be a subspace of $V^{(n)}$ and $\sum_i \lambda_i z_i\otimes x_i$ a nonzero element of $W\otimes V^{(m)}$, where $\lambda_i\in k$, $z_i\in W$ and $x_i=\sum_{j=(j_1,\cdots ,j_m)} \mu_j^i v_{j_1}\cdots v_{j_m}$ with $\mu_j^i\in k$ and $v_{j_t}\in \mathcal{B}$. Then $\sum_{i, j=(j_1,\cdots ,j_m)} \lambda_i \mu_j^i z_i \in W$.
\end{lema}

\begin{prof}
Let $f:V^{(m)}\rightarrow k$ be the linear map such that $f(v_{i_1} \cdots v_{i_m})=1$ for all $v_{i_h}\in \mathcal{B}$. Then $\sum_{
i,j=(j_1,\cdots ,j_m)} \lambda_i \mu_j^i z_i$ is the image of $\sum_{i,j=(j_1,\cdots ,j_m)} \lambda_i z_i\otimes x_i$ by $1_W\otimes f$. \qed
\end{prof}

\bigskip

Take now $A=T(V)/I$.

\begin{lema}\label{lem:belongtoauxideal}
Given $x=\sum_{i=(i_1,\cdots ,i_n)} \lambda_i v_{i_1} \cdots v_{i_{n-1}}\otimes v_{i_n} \in I_{n-1}\otimes V$, with $\lambda_i \in k$ and $v_{i_j} \in \mathcal{B}$ ($1\le j\le n$), suppose that if $i\ne i'$, then $v_{i_n} \ne v_{i'_n}$. Therefore, $\lambda_i v_{i_1} \cdots v_{i_{n-1}} \in I_{n-1}$ for all $i=(i_1,\cdots ,i_n)$.
\end{lema}

\begin{prof}
By Lemma \ref{lem:inthesubsp}, we get that $\sum_{i=(i_1,\cdots ,i_n)} \lambda_i v_{i_1} \cdots v_{i_{n-1}} \in I_{n-1}$.

We fix a term $h=(h_1,\cdots ,h_n)$ and consider $y=(\sum_{i=(i_1,\cdots ,i_n)} \lambda_i v_{i_1} \cdots v_{i_{n-1}})\otimes v_{h_n} \in I_{n-1}\otimes V$. Then \[
x-y =\smashoperator{\sum \limits_{\substack{
i=(i_1,\cdots ,i_n)
\\
(i_1,\cdots ,i_n)\ne (h_1,\cdots ,h_n)}}} (\lambda_i v_{i_1} \cdots v_{i_{n-1}}) \otimes (v_{i_n}-v_{h_n}) \in I_{n-1}\otimes V
\]
and it is nonzero. Therefore, 
\[
\smashoperator{\sum \limits_{\substack{i=(i_1,\cdots ,i_n)
\\
(i_1,\cdots ,i_n)\ne (h_1,\cdots ,h_n)}}} \lambda_i v_{i_1}\cdots v_{i_{n-1}} \in I_{n-1}.
\]
Thus,
\[
\smashoperator{\sum \limits_{i=(i_1,\cdots ,i_n)}} \lambda_i v_{i_1} \cdots v_{i_{n-1}} - \smashoperator{\sum \limits_{\substack{
i=(i_1,\cdots ,i_n)
\\
(i_1,\cdots ,i_n)\ne (h_1,\cdots ,h_n)}}} \lambda_i v_{i_1} \cdots v_{i_{n-1}}= \lambda_{h} v_{h_1} \cdots v_{h_{n-1}} \in I_{n-1}.
\]
\noindent Since $h$ is arbitrary, the result follows. \qed
\end{prof}

\begin{rmk}\label{rmk:belongtoideal}
Consider a nonzero element $x=\sum_{i=(i_1,\cdots ,i_n)} \lambda_i v_{i_1} \cdots v_{i_{n-1}}\otimes v_{i_n} \in I_{n-1}\otimes V$ and suppose that there exist $i$ and $i'$ such that $i\ne i'$ and $v_{i_n}=v_{i'_n}$. Reordering the sum, we may see that $x=\sum_{t=1}^s (\sum_{\{ i / v_{i_n} = v_t\}} \lambda_i v_{i_1} \cdots v_{i_{n-1}})\otimes v_t$. By Lemma \ref{lem:belongtoauxideal}, $\sum_{\{ i/ v_{i_n}=v_t\}} \lambda_i v_{i_1} \cdots v_{i_{n-1}} \in I_{n-1}$ for all $t$. It follows that $\sum_{i=(i_1,\cdots ,i_n)} \lambda_i v_{i_1} \cdots v_{i_{n-1}} \in I_{n-1}$, since $I$ is an ideal.
\end{rmk}


\section{\texorpdfstring{Pure resolutions of $k$}{sec:pureresolutions}}

In this section, after some preliminaries, we construct a minimal graded projective resolution of the trivial $A$-module $k$. The description of the succesive kernels of the differentials in this resolution will be related to some specific conditions. The situation is summarize in Theorem \ref{thm:kerdeltai}.

We recall the following definition.

\begin{defi}
Let $M$ be a $\mathbb{Z}$-graded left $A$-module provided of a minimal graded projective resolution
\begin{equation}\label{eq:resolution}
\cdots \longrightarrow P_n\overset{d_n}{\longrightarrow} P_{n-1}\longrightarrow \cdots
\longrightarrow P_0\overset{d_0}{\longrightarrow} M \longrightarrow 0.
\end{equation}

\noindent We say that \eqref{eq:resolution} is a \textbf{pure projective resolution} if for each $n\in \mathbb{N}_0$ there exists $s_n\in \mathbb{Z}$ such that $P_n$ is $s_n$-pure.
\end{defi}

\begin{rmk}
Every module admits a pure projective resolution. In the graded category this resolution may be chosen minimal (see \cite{GMV}.)
\end{rmk}

\begin{defi}
Let $S\subseteq V^{(n)}$ and $T\subseteq V^{(m)}$ be two subspaces with $n\le m$. We say that $S$ and $T$ are \textbf{exclusive} if the intersection of the two-sided ideal generated by $S$ with $T$ is trivial. Explicitly, $(\sum_{i+j+n=m} V^{(i)}\otimes S\otimes V^{(j)})\cap \hskip 0.1cm T=0$. 
\end{defi}

\medskip
 
From now on, we fix two integers $a$ and $b$ such that $2\le a<b$, $R_a$ and $R_b$ subspaces of $V^{(a)}$ and $V^{(b)}$ respectively and $R=R_a \oplus R_b$. We also assume that $R_a$ and $R_b$ are exclusive.

The two-sided ideal $I=I(R)$ generated by $R$ in the tensor algebra $T(V)$ is $\mathbb{Z}$-graded, i.e. $I=\bigoplus_{n\in \mathbb{Z}} I_n$ where:
\begin{align*}
I_n &= 0 & \hskip 0.2cm & \text{if} \hskip 0.2cm n<a,
\\
I_n &=\smashoperator{\sum \limits_{i+j+a=n}} V^{(i)}\otimes R_a \otimes V^{(j)} & \hskip 0.2cm & \text{if} \hskip 0.2cm a\le n<b,
\\
I_n &= \smashoperator{\sum \limits_{i+j+a=n}} V^{(i)}\otimes R_a \otimes V^{(j)}+ \smashoperator{\sum \limits_{h+l+b=n}} V^{(h)}\otimes R_b \otimes V^{(l)} & \hskip 0.2cm & \text{if} \hskip 0.2cm b\le n.
\end{align*}

\medskip

$R$ is called a {\bf space of relations} for the $\mathbb{Z}$-graded algebra $A=T(V)/I$. The homogeneous components of $A$ are the vector subspaces $A_n=V^{(n)}/I_n$ for $n\in \mathbb{N}_0$ and zero for $n<0$. The algebra $A$ is generated in degree $1$ and $A_n=V^{(n)}$ for $0\le n<a$.

We fix the following bases $\mathcal{B}_a=\{ r_1,\cdots ,r_p\}$ of $R_a$ and $\mathcal{B}_b=\{ s_1,\cdots ,s_{p'}\}$ of $R_b$ where $r_i$ and $s_t$ are uniquely written as $r_i= \sum_{j=(j_1,\cdots ,j_a)} \lambda_j^i v_{j_1} \cdots v_{j_a}$ with $v_{j_h}\in \mathcal{B}$ and $\lambda_j^i\in k$ for $1\le i\le p$, and $s_t= \sum_{j=(j_1,\cdots ,j_b)} \mu_j^t v_{j_1} \cdots v_{j_b}$ with $v_{j_h}\in \mathcal{B}$ y $\mu_j^t\in k$ for $1\le t\le p'$ respectively.

\begin{rmk}
The notion of space of relations is already present in \cite{B2} and it involves a minimality idea. Suppose $R$ is a space of relations for the two-sided ideal $I$, and assume that any basis of $R$ is a minimal set of generators of $I$. The minimality condition implies that $R_a$ and $R_b$ are exclusive. In fact, suppose that there exists $v\ne 0$ such that $v\in (\sum_{i+j+a=b} V^{(i)}\otimes R_a \otimes V^{(j)})\cap R_b$. Then $v=\sum_{h=1}^{p'} \lambda_h s_h$ with $\lambda_h\in k$, $1\le h\le p'$, not all zero. Let $h'$ be such that $\lambda_{h'}\ne 0$, then $s_{h'}=v-\frac{1}{\lambda_{h'}} \sum_{h=1,h\ne h'}^{p'} \lambda_h s_h$, and since $v$ belongs to the two-sided ideal generated by $R_a$, $R$ is not a minimal set of generators.
\end{rmk}

Our purpose now is to obtain a pure projective resolution of the trivial graded $A$-module $k$. The natural projection $\epsilon : A\rightarrow k$ is a projective cover of $k$ and $\ker \epsilon =\bigoplus_{n\ge 1} A_n$, which is pure in degree $1$ and $(\ker \epsilon )_1=V$.

Next we briefly describe the kernels of the morphisms of a resolution of the $A$-module $k$. For further details we refer to \cite{R}.


\subsection{\texorpdfstring{Description of $\ker \delta_1$}{sec:kerdelta1}}

The injection $\tilde{g}_1: V\rightarrow \ker \epsilon$ induces the linear map $g_1:A\otimes V\rightarrow \ker \epsilon$, $\overline{\alpha}\otimes v\mapsto \overline{\alpha v}$, which is a projective cover of $\ker \epsilon$. Let the inclusion $inc_0:\ker \epsilon \rightarrow A\otimes k$, we consider the map $\delta_1=inc_0\circ g_1$. One can easily check (see \cite{R}) that
\begin{itemize}
\item[$\bullet$] $(\ker \delta_1)_n=0$ if $n<a$,
\item[$\bullet$] $(\ker \delta_1)_a=R_a$,
\item[$\bullet$] $(\ker \delta_1)_n \simeq \frac{I_{n-1}\otimes V + V^{(n-a)} \otimes R_a}{I_{n-1}\otimes V}$ if $a<n<b$,
\item[$\bullet$] $(\ker \delta_1)_n \simeq \frac{I_{n-1}\otimes V + V^{(n-a)}\otimes R_a+ V^{(n-b)}\otimes R_b}{I_{n-1}\otimes V}$ if $n\ge b$.
\end{itemize}

The following lemma is a graded version of a well-known result proved in \cite{B5} for the non graded case. It provides a description of $\ker \delta_1$.

\begin{lema}
The $A$-module $\ker \delta_1$ is $2$-pure in degrees $a$ and $b$.
\end{lema}

\begin{prof}
In order to prove the statement, note that
\begin{itemize}
\item[$\bullet$] if $n<a$, $(\ker \delta_1)_n=0$,
\item[$\bullet$] if $n=a$, $(\ker \delta_1)_a=R_a$,
\item[$\bullet$] if $a<n<b$,
\[
(\ker \delta_1)_n \simeq \frac{I_{n-1}\otimes V + V^{(n-a)}\otimes R_a}{I_{n-1}\otimes V} \simeq \frac{V^{(n-a)}\otimes R_a}{(I_{n-1}\otimes V)\cap (V^{(n-a)}\otimes R_a)}=A_{n-a}\otimes R_a,
\]

\item[$\bullet$] if $n\ge b$
\begin{align*}
& (\ker \delta_1)_n \simeq \frac{I_{n-1}\otimes V + V^{(n-a)}\otimes
R_a + V^{(n-b)}\otimes R_b}{I_{n-1}\otimes V}
\\
& \simeq \frac{V^{(n-a)}\otimes R_a \oplus V^{(n-b)}\otimes R_b}{(I_{n-1}\otimes V)\cap (V^{(n-a)}\otimes R_a \oplus V^{(n-b)}\otimes R_b)}=A_{n-a}\otimes R_a\oplus A_{n-b}\otimes R_b.
\end{align*}

\noindent Since $R_a$ and $R_b$ are exclusive the sums are direct. 
\end{itemize}

Consequently, by Proposition \ref{prop:essentialineachdegree}, $\ker \delta_1$ is $2$-pure in degrees $a$ and $b$. \qed
\end{prof}


\subsection{\texorpdfstring{Description of $\ker \delta_2$}{sec:kerdelta2}}

Let $g_2:A\otimes R\rightarrow \ker \delta_1$ be the morphism induced by the injection $\tilde{g}_2:R\rightarrow \ker \delta_1$, then $g_2$ is a projective cover of $\ker \delta_1$. Let also $inc_1:\ker \delta_1\rightarrow A\otimes V$ and consider $\delta_2=inc_1\circ g_2$. We have thus obtained the beginning of a resolution
\[
A\otimes R\overset{\delta_2}{\longrightarrow} A\otimes V \overset{\delta_1}{\longrightarrow} A\otimes k \overset{\epsilon}{\longrightarrow} k \longrightarrow 0.
\]

\noindent Note that $k$ and $V$ are respectively concentrated in degrees $0$ and $1$ and that $R$ is $2$-concentrated in degrees $a$ and $b$. We define the following vector spaces, both concentrated in degree $n$,
\begin{align*}
J_n^a & = \smashoperator{\bigcap \limits_{i+j+a=n}} V^{(i)}\otimes R_a \otimes V^{(j)} & \text{for} \hskip 0.2cm n\ge a,
\\
J_n^b & = \smashoperator{\bigcap \limits_{s+t+b=n}} V^{(s)}\otimes R_b \otimes V^{(t)} & \text{for} \hskip 0.2cm n\ge b.
\end{align*}

\medskip

Since
\[
(\ker \delta_2)_n =0 \hskip 0.2cm \text{if} \hskip 0.2cm n\le a \hskip 0.2cm \text{and} \hskip 0.2cm (\ker \delta_2)_{a+1} = J_{a+1}^a,
\]

\noindent the projective cover of $\ker \delta_2$ must include $A\otimes J_{a+1}^a$ but may also include $s$-pure modules with $s>a+1$. We want to describe this projective cover.

It is straightforward to see that for $n=a+m$ with $2\le m\le \min \{ a-1, b-a\}$
\begin{equation}
(\ker \delta_2)_n = (V^{(m)}\otimes R_a) \cap \big(\smashoperator{\sum \limits_{i+j=m-1}} V^{(i)}\otimes R_a\otimes V^{(j)}\otimes V\big)\supseteq V^{(m-1)}\otimes J_{a+1}^a.
\end{equation}

\noindent Then,
\begin{equation}\label{eq:econa}
(\ker \delta_2)_n =A_{m-1} \otimes J_{a+1}^a \Longleftrightarrow
(V^{(m)}\otimes R_a)\cap \big(\smashoperator{\sum \limits_{i+j=m-1}} V^{(i)}\otimes R_a\otimes V^{(j)}\otimes V\big) = V^{(m-1)}\otimes J_{a+1}^a .
\end{equation}

\noindent Note that if the right hand side in \eqref{eq:econa} holds for $m=a-1$, then for $2\le m\le a-1$ and $2\le t\le m-2$,
\begin{equation}
(V^{(m)}\otimes R_a)\cap (V^{(m-t)}\otimes R_a\otimes V^{(t)}+\cdots +V^{(m-1)}\otimes R_a\otimes V)\subseteq V^{(m-1)}\otimes J_{a+1}^a.
\end{equation}

\noindent Since $V$ is $k$-faithfully flat, we get \eqref{eq:econa} for $2\le m\le a-1$.

\medskip

For $n=a+t=b+h$ with $1\le h\le 2a-b-1$, $(\ker \delta_2)_n$ is equal to
\[
[(V^{(t)}\otimes R_a) \oplus (V^{(h)}\otimes R_b)]\cap \big(\smashoperator{\sum \limits_{i+j=t-1}} V^{(i)}\otimes R_a\otimes V^{(j)}\otimes V+\smashoperator{\sum \limits_{s+t=h-1}} V^{(s)}\otimes R_b\otimes V^{(t)}\otimes V\big),
\]

\noindent and contains $(V^{(t-1)}\otimes J_{a+1}^a)\oplus (V^{(h-1)}\otimes J_{b+1}^b)$. 

\medskip 

Consider now the following relations
\begin{equation}\label{eq:ec}
\begin{split}
(V^{(s)}\otimes R_a)\cap \big(\smashoperator{\sum \limits_{i+j=s-1}} V^{(i)}\otimes R_a \otimes V^{(j)}\otimes V\big) & = V^{(s-1)}\otimes J_{a+1}^a,
\\
(V^{(l)}\otimes R_b)\cap \big(\smashoperator{\sum \limits_{s+t=l-1}} V^{(s)}\otimes R_b \otimes V^{(t)}\otimes V\big) & = V^{(l-1)}\otimes J_{b+1}^b.
\end{split}
\end{equation}

\noindent With similar arguments as before, \eqref{eq:ec} holds for $l\le b-2$ whenever it holds for $l=b-1$. Moreover, \eqref{eq:ec} implies the right hand side equality in \eqref{eq:econa}.

The relations \eqref{eq:ec} for $s=a-1$ and $l=b-1$ will be called \textbf{extra conditions} and will be denoted by \textit{(e.c.)}.

\bigskip

Next we recall a definition that will be generalized afterwards.

\begin{defi}\cite{B2, J, O}
A triple $(E,F,G)$ of subspaces of a given vector space is said to be \textbf{distributive} if
\begin{equation}
E\cap (F+G)= (E\cap F) + (E\cap G).
\end{equation}
\end{defi}

\begin{prop}\label{prop:eciffdistributivity}
The \textit{(e.c.)} hold if and only if for $2\le m\le a-1$ and $2\le l\le b-1$ the triples
\begin{equation}
(V^{(m)}\otimes R_a, R_a\otimes V^{(m)}, \smashoperator[r]{\sum \limits_{i+j=m-2}} V\otimes V^{(i)}\otimes R_a \otimes V^{(j)}\otimes V) \text{ and}
\end{equation}
\[
(V^{(l)}\otimes R_b, R_b\otimes V^{(l)}, \smashoperator[r]{\sum \limits_{s+t=l-2}} V\otimes V^{(s)}\otimes R_b \otimes V^{(t)}\otimes V)
\]

\noindent are distributive and there are inclusions
\begin{align*}
(V^{(m)}\otimes R_a)\cap (R_a \otimes V^{(m)}) & \subseteq V^{(m-1)}\otimes R_a\otimes V,
\\
(V^{(l)}\otimes R_b)\cap (R_b \otimes V^{(l)}) & \subseteq V^{(l-1)}\otimes R_b \otimes V.
\end{align*}
\end{prop}

\begin{prof}
Let $E=V^{(m)}\otimes R_a$, $F=R_a\otimes V^{(m)}$ and $G=\sum_{i+j=m-2} V\otimes V^{(i)}\otimes R_a \otimes V^{(j)}\otimes V$. Recall that if the \textit{(e.c.)} are satisfied, then we have \eqref{eq:ec} for $2\le s\le a-1$ and $1\le l\le b-2$. If we assume that the \textit{(e.c.)} hold, then in \eqref{eq:ec}, $E\cap (F+G)= V^{(m-1)}\otimes J_{a+1}^a\subseteq E\cap G \subseteq (E\cap F)+(E\cap G)$. The argument is analogous for $b$. The inclusions are trivial.

Conversely, because of the distributivity of the first triple and the first inclusion for $m$ such that $2\le m\le a-1$, 
\begin{align*}
& (V^{(a-2)}\otimes R_a)\cap \big(\smashoperator{\sum \limits_{i+j=a-3}} V^{(i)}\otimes R_a\otimes V^{(j)}\otimes V\big)
\\
& = (V^{(a-2)}\otimes R_a)\cap (R_a\otimes V^{(a-2)})+(V^{(a-2)}\otimes R_a)\cap \big(\smashoperator{\sum \limits_{i+j=a-4}} V\otimes V^{(i)}\otimes R_a\otimes V^{(j)}\otimes V\big)
\\
& \subseteq (V^{(a-3)}\otimes R_a \otimes V)+ V\otimes [(V^{(a-3)} \otimes R_a)\cap \big(\sum_{i+j=a-4} V^{(i)}\otimes R_a\otimes V^{(j)}\otimes V\big)]
\\
& = (V^{(a-3)}\otimes R_a \otimes V)+ V\otimes [(V^{(a-3)}\otimes R_a) \cap (R_a\otimes V^{(a-3)})+(V^{(a-3)}\otimes R_a) \cap
\\
& \big(\smashoperator{\sum \limits_{i+j=a-5}} V\otimes V^{(i)}\otimes R_a\otimes V^{(j)}\otimes V\big)]
\\
& \subseteq (V^{(a-3)}\otimes R_a\otimes V)+ V\otimes [(V^{(a-4)}\otimes R_a\otimes V)+(V^{(a-3)}\otimes R_a)\cap
\\
& \big(\smashoperator{\sum \limits_{i+j=a-5}} V\otimes V^{(i)}\otimes R_a\otimes V^{(j)}\otimes V\big)]
\\
& =(V^{(a-3)}\otimes R_a\otimes V)+V^{(2)}\otimes [V^{(a-5)}\otimes R_a\otimes V+(V^{(a-4)}\otimes R_a)\cap \big(\smashoperator{\sum \limits_{i+j=a-5}} V^{(i)}\otimes R_a\otimes V^{(j)}\otimes V\big)]
\\
& =(V^{(a-3)}\otimes R_a\otimes V)+V^{(2)}\otimes [(V^{(a-4)}\otimes R_a)\cap \big(\smashoperator{\sum \limits_{i+j=a-5}} V^{(i)}\otimes R_a\otimes V^{(j)}\otimes V\big)].
\end{align*}

\noindent In general, $(V^{(a-3)}\otimes R_a\otimes V)+V^{(s-3)}\otimes [(V^{(a-s+1)}\otimes R_a)\cap (\sum_{i+j=a-s} V^{(i)}\otimes R_a\otimes V^{(j)}\otimes V)]\subseteq (V^{(a-3)}\otimes R_a\otimes V)+V^{(s-3)}\otimes [V^{(a-s)}\otimes R_a\otimes V +(V^{(a-s+1)}\otimes R_a)\cap (\sum_{i+j=a-s} V^{(i)}\otimes R_a\otimes V^{(j)}\otimes V)]=(V^{(a-3)}\otimes R_a \otimes V)+V^{(s-3)}\otimes [(V^{(a-s+1)}\otimes R_a)\cap (\sum_{i+j=a-s} V^{(i)}\otimes R_a\otimes V^{(j)}\otimes V)]$ with $3\le s\le a$ and therefore
\[
V^{(a-3)}\otimes R_a \otimes V \subseteq V^{(a-3)}\otimes J_{a+1}^a.
\]

\noindent Hence, the first relation in \eqref{eq:ec} holds. A similar argument can be applied for the second relation in \eqref{eq:ec}. \qed
\end{prof}

\begin{lema}\label{lem:equivalenceforincl}
For $2\le m\le a-1$ and $2\le l\le b-1$, the following inclusions hold
\begin{equation}\label{eq:auxinclusion}
(V^{(m)}\otimes R_a)\cap (R_a\otimes V^{(m)}) \subseteq V^{(m-1)}\otimes R_a\otimes V,
\end{equation}
\[
(V^{(l)}\otimes R_b)\cap (R_b\otimes V^{(l)}) \subseteq V^{(l-1)}\otimes R_b\otimes V,
\]

\noindent if and only if the following equalities are satisfied
\begin{equation*}
(V^{(m)}\otimes R_a)\cap (R_a\otimes V^{(m)})=J_{a+m}^a, \hskip 0.2cm (V^{(l)}\otimes R_b)\cap (R_b\otimes V^{(l)})=J_{b+l}^b.
\end{equation*}
\end{lema}

\begin{prof}
For the necessity we use that $J_{a+m}^a \subseteq (V^{(m)}\otimes R_a)\cap (R_a\otimes V^{(m)})$. We shall immediately prove the other inclusion by induction on $m$. If $m=2$, the equality becomes, using \eqref{eq:auxinclusion}
\[
J_{a+2}^a=(V^{(2)}\otimes R_a)\cap (V\otimes R_a \otimes V) \cap (R_a\otimes V^{(2)})=(V^{(2)}\otimes R_a)\cap (R_a\otimes V^{(2)}).
\]

\noindent We suppose now that the inclusion we need is valid for some $m<a-1$; let $\sum_{i=(i_1,\cdots ,i_{a+m+1})} \lambda_i v_{i_1}\cdots$  $v_{i_{a+m+1}}$ $\in (V^{(m+1)}\otimes R_a)\cap (R_a\otimes V^{(m+1)})$ with $v_{i_j}\in \mathcal{B}$, then $\sum_{i=(i_1,\cdots ,i_{a+m+1})} \lambda_i v_{i_1}\cdots v_{i_a}$ $\in R_a$ and $\sum_{i=(i_1,\cdots ,i_{a+m+1})} \lambda_i v_{i_{m+2}}\cdots v_{i_{a+m+1}} \in R_a$, but using \eqref{eq:auxinclusion}, $\sum_{i=(i_1,\cdots , i_{a+m+1})} \lambda_i v_{i_{m+1}}\cdots v_{i_{a+m}} \in R_a$ if $m+1\le a-1$. Hence, $\sum_{i=(i_1,\cdots ,i_{a+m+1})} \lambda_i v_{i_1}\cdots v_{i_{a+m}} \in (R_a\otimes V^{(m)})\cap (V^{(m)}\otimes R_a) = J_{a+m}^a$ by inductive hypothesis. We can then conclude that $\sum_{i=(i_1,\cdots ,i_{a+m+1})} \lambda_i v_{i_1}\cdots v_{i_{a+m+1}} \in \bigcap_{h+j=m} V^{(h)}\otimes R_a \otimes V^{(j)}\otimes V$ and since it also belongs to $R_a\otimes V^{(m+1)}$, it belongs to $J_{a+m+1}^a$ and the desired equality is satisfied.

The proof for the equality involving $J_{b+l}^b$ is analogous.

The converse is trivial since the following inclusions are satisfied:
\[
J_{a+m}^a \subseteq V^{(m-1)}\otimes R_a \otimes V, \hskip 0.2cm J_{b+l}^b \subseteq V^{(l-1)}\otimes R_b \otimes V.
\] \qed
\end{prof}

\bigskip

Next we state two technical lemmas concerning $A$-modules (not necessarily graded) that will be used afterwards. The proofs are straightforward.

\begin{lema}\label{lem:samemod}
Let $I$, $N$ and $M$ be $A$-modules such that $I\subseteq N\subseteq M$, $N/I\simeq M/I$ and the following diagram is commutative
\[
\xymatrix@R-20pt
{
N
\ar@{^(->}[r]
\ar[dd]^{\pi}
& M
\ar[dd]^{\pi'}
\\
&
\\
N/I
\ar[r]^{\sim}
& M/I
.}
\]

\noindent Then $N=M$.
\end{lema}

\begin{lema}\label{lem:intersection}
Let $I$, $J$ and $M$ be subspaces of a $k$-vector space $N$ such that $I\subseteq M$, and $I\subseteq J$. Suppose also that
\begin{align*}
\phi :M/I & \longrightarrow N/J
\\
\overline{m} & \longmapsto [m]
\end{align*}

\noindent is an isomorphism. Then $M\cap J=I$.
\end{lema}

\begin{defi}\label{def:multdisttriple}
A 4-tuple of subspaces of a given vector space of the form $(E, E', F_1+\cdots +F_t, G_1+\cdots +G_{t'})$ is said to be \textbf{multidistributive} if $E\cap E'=0$ and
\begin{multline}
(E\oplus E')\cap (F_1+\cdots +F_t+G_1+\cdots G_{t'})
\\
= (E\cap F_1)+\cdots +(E\cap F_t)\oplus (E'\cap G_1)+\cdots +(E'\cap G_{t'}).
\end{multline}
\end{defi}

\medskip

The following proposition gives an equivalence in terms of distributivity and multidistributivity and the \textit{(e.c.)}, to decide when $\ker \delta_2$ is $2$-pure. The purpose will be then to generalize this equivalence for the other morphisms in the resolution.

\begin{prop}\label{prop:kerdelta2pure}
If we assume that the second equality of the \textit{(e.c.)} is satisfied; i.e. 
\[
(V^{(b-1)}\otimes R_b)\cap \big(\smashoperator{\sum \limits_{s+t=b-2}} V^{(s)}\otimes R_b \otimes V^{(t)}\otimes V\big) = V^{(b-2)}\otimes J_{b+1}^b,
\]

\noindent then $\ker \delta_2$ is $2$-pure in degrees $a+1$ and $b+1$ if and only if \eqref{eq:econa} holds and, for all $n\ge 2a$, $(E,F,G)$ and $(E',E'',F'+G',F''+G'')$ are respectively distributive and multidistributive where:
\begin{align*}
E & = E'= V^{(n-a)}\otimes R_a,
\\
F & = R_a\otimes V^{(n-a)}+ \cdots + V^{(n-2a)}\otimes R_a \otimes V^{(a)},
\\
G & = G' = V^{(n-2a+1)}\otimes R_a\otimes V^{(a-1)}+ \cdots + V^{(n-a-1)}\otimes R_a \otimes V,
\\
F' & = R_a\otimes V^{(n-a)}+ \cdots + V^{(n-2a)}\otimes R_a \otimes V^{(a)} + R_b\otimes V^{(n-b)} + \cdots + V^{(n-a-b)}\otimes R_b \otimes V^{(a)},
\\
E'' & = V^{(n-b)}\otimes R_b,
\\
F'' & = R_a\otimes V^{(n-a)}+ \cdots + V^{(n-a-b)}\otimes R_a \otimes V^{(b)}+ R_b\otimes V^{(n-b)} + \cdots + V^{(n-2b)}\otimes R_b \otimes V^{(b)},
\\
G'' & = V^{(n-2b+1)}\otimes R_b\otimes V^{(b-1)}+ \cdots + V^{(n-b-1)}\otimes R_b \otimes V.
\end{align*}
\end{prop}

\begin{prof}
Fix $V^{(i)}=0$ for $i<0$. We have already proved that
\[
(\ker \delta_2)_n =A_{m-1} \otimes J_{a+1}^a \hskip 0.2cm \text{for} \hskip 0.2cm n=m+a, \hskip 0.2cm \text{with} \hskip 0.2cm 2\le m\le \min \{a-1,b-a\},
\]

\noindent if and only if \eqref{eq:econa} is satisfied. We shall first prove the if part. Suppose now that the \textit{(e.c.)} hold and fix an integer $n\ge 2a$. We want to describe $(\ker \delta_2)_n$ for arbitrary $n$.
\begin{itemize}
\item[(i)] If $n\le b$, $(\ker \delta_2)_n\subseteq A_{n-a}\otimes R_a$. Note that if $n=b$, $A_{n-b}=k$ and if $s\in R_b$, $\delta_2 (s)$ is not zero, if not,
\[
0=\delta_2(s)=\delta_2\big(\smashoperator{\sum \limits_{\substack{l=1
\\
h=(h_1,\cdots ,h_b)}}^{p'}} \eta_l \mu_h^l v_{h_1}\cdots v_{h_b}\big)=\smashoperator{\sum \limits_{\substack{l=1
\\
h=(h_1,\cdots ,h_b)
}}^{p'}} \overline{\eta_l \mu_h^l v_{h_1}\cdots v_{h_{b-1}}}\otimes v_{h_b},
\]

\noindent and then $\sum_h \sum_{l=1}^{p'} \eta_l \mu_h^l v_{h_1}\cdots v_{h_{b-1}}\in I_{b-1}$ but in this case $R_a$ and $R_b$ would not be exclusive.

Let $\overline{\alpha_i} \in A_{n-a}$, $\rho_i\in R_a$, such that $0= \delta_2 (\sum_i \overline{\alpha_i} \otimes \rho_i)=\delta_2 (\sum_{i,h,j=(j_1,\cdots ,j_a)} \gamma_h^i \lambda_j^h \overline{\alpha_i} \otimes v_{j_1}\cdots v_{j_a})= \sum_{i,h,j=(j_1,\cdots ,j_a)} \overline{\gamma_h^i \lambda_j^h \alpha_i v_{j_1}\cdots v_{j_{a-1}}}\otimes v_{j_a}$. So, $\sum_{i,h,j=(j_1,\cdots ,j_a)} \gamma_h^i \lambda_j^h \alpha_i v_{j_1}\cdots v_{j_{a-1}}$ belongs to $I_{n-1}$.

Consider the following subspace of $T(V)$
\[
N_n=(V^{(n-a)}\otimes R_a)\cap (I_{n-1}\otimes V)=(V^{(n-a)}\otimes R_a)\cap (V^{(n-a-1)}\otimes R_a \otimes V + \cdots +R_a\otimes V^{(n-a)}).
\]

\noindent We get that $(\ker \delta_2)_n\subseteq \frac{N_n}{I_{n-a}\otimes R_a}$. The other inclusion is trivial.
\item[(ii)] If $n>b$, $(\ker \delta_2)_n\subseteq (A_{n-a}\otimes R_a)\oplus (A_{n-b}\otimes R_b)$. Let $\overline{\alpha_i} \in A_{n-a}$, $\overline{\beta_h}\in A_{n-b}$, $\rho_i\in R_a$ and $\theta_h\in R_b$, if
\begin{align*}
0 & =\delta_2 \big(\smashoperator{\sum \limits_i} \overline{\alpha_i}\otimes \rho_i +\smashoperator{\sum \limits_h} \overline{\beta_h}\otimes \theta_h\big)
\\
& =\delta_2 \big(\smashoperator{\sum \limits_{\substack{i,t
\\
j=(j_1,\cdots ,j_a)}}} \gamma_t^i \lambda_j^t \overline{\alpha_i} \otimes v_{j_1}\cdots v_{j_a} + \smashoperator{\sum \limits_{\substack{h,m
\\
l=(l_1,\cdots ,l_b)}}} \eta_m^h \mu_l^m \overline{\beta_h} \otimes v_{l_1}\cdots v_{l_b}\big)
\\
& = \smashoperator{\sum \limits_{\substack{i,t
\\
j=(j_1,\cdots ,j_a)}}} \gamma_t^i\lambda_j^t \overline{\alpha_i v_{j_1}\cdots v_{j_{a-1}}}\otimes v_{j_a} + \smashoperator{\sum \limits_{\substack{h,m
\\
l=(l_1,\cdots ,l_b)}}} \eta_m^h \mu_l^m \overline{\beta_h v_{l_1}\cdots v_{l_{b-1}}}\otimes v_{l_b},
\end{align*}

\noindent then $\sum_{i,t,j} \gamma_t^i \lambda_j^t \alpha_i v_{j_1}\cdots v_{j_{a-1}} +\sum_{h,m,l} \eta_m^h \mu_l^m \beta_h v_{l_1}\cdots v_{l_{b-1}} \in I_{n-1}$.

The subspace
\begin{align*}
N_n = & (V^{(n-a)}\otimes R_a \oplus V^{(n-b)}\otimes R_b)\cap (V^{(n-a-1)} \otimes R_a\otimes V +
\\
& V^{(n-a-2)}\otimes R_a \otimes V^{(2)} + \cdots +R_a\otimes V^{(n-a)}+
\\
& V^{(n-b-1)}\otimes R_b\otimes V + V^{(n-b-2)}\otimes R_b \otimes V^{(2)} + \cdots + R_b \otimes V^{(n-b)})
\end{align*}

\noindent is included in $I_{n-1}\otimes V$. As a consequence we obtain that
\[
(\ker \delta_2)_n = \frac{N_n}{I_{n-a}\otimes R_a \oplus I_{n-b}\otimes R_b}.
\]
\end{itemize}

In order to prove the distributivity and multidistributivity we note that:
\begin{align*}
E\cap (F+G) & = N_n & \hskip 0.2cm & \text{in case \textit{(i)}},
\\
(E'\oplus E'')\cap (F'+F''+G'+G'') & =N_n & \hskip 0.2cm & \text{in case \textit{(ii)}}.
\end{align*}

We shall analyse separately both cases for $(E,F,G)$ and $(E', E'',F'+F'',G'+G'')$, which are respectively distributive and multidistributive, using that the \textit{(e.c.)} are satisfied.
\begin{itemize}
\item[$\bullet$] If $n\le b$, $N_n = E\cap (F+G) = (E\cap F) + (E\cap G) =_{\text{\textit{(e.c.)}}} (I_{n-a}\otimes R_a)+(V^{(n-a-1)}\otimes J_{a+1}^a).$

\noindent Hence, 
\[
(\ker \delta_2)_n = \frac{N_n}{I_{n-a}\otimes R_a} = \frac{V^{(n-a-1)}\otimes J_{a+1}^a}{(V^{(n-a-1)}\otimes J_{a+1}^a)\cap (I_{n-a}\otimes R_a)}
\]

\noindent and it is surjectively mapped to $\frac{V^{(n-a-1)}\otimes J_{a+1}^a}{I_{n-a-1}\otimes J_{a+1}^a} = A_{n-a-1}\otimes J_{a+1}^a$ by a morphism $\varphi$ induced by the identity on $V^{(n-a-1)}\otimes J_{a+1}^a$. Moreover, the domain and the image of $\varphi$ are isomorphic as $k$-vector spaces (see \cite{B2}). For dimensional reasons, $\varphi$ is an isomorphism.

\item[$\bullet$] If $n>b$,
\[
E'\cap F' = I_{n-a}\otimes R_a, \hskip 0.2cm E''\cap F'' = I_{n-b}\otimes R_b, \hskip 0.2cm E'\cap G' = V^{(n-a-1)}\otimes J_{a+1}^a, \hskip 0.2cm E''\cap G'' = V^{(n-b-1)}\otimes J_{b+1}^b,
\]

\noindent where the last two hold because of the \textit{(e.c.)}. The kernel is then characterized as follows,
\begin{align*}
(\ker \delta_2)_n & = \frac{N_n}{I_{n-a}\otimes R_a \oplus I_{n-b} \otimes R_b}
\\
& =\frac{(I_{n-a}\otimes R_a + V^{(n-a-1)}\otimes J_{a+1}^a) \oplus (I_{n-b}\otimes R_b + V^{(n-b-1)}\otimes J_{b+1}^b)}{I_{n-a}\otimes R_a\oplus I_{n-b}\otimes R_b}
\\
& \simeq \frac{I_{n-a}\otimes R_a + V^{(n-a-1)}\otimes J_{a+1}^a}{I_{n-a}\otimes R_a} \oplus \frac{I_{(n-b)}\otimes R_b + V^{(n-b-1)}\otimes J_{b+1}^b}{I_{n-b}\otimes R_b}
\\
& \simeq \frac{V^{(n-a-1)}\otimes J_{a+1}^a}{(V^{(n-a-1)}\otimes J_{a+1}^a)\cap (I_{n-a}\otimes R_a)}\oplus \frac{V^{(n-b-1)}\otimes J_{b+1}^b}{(V^{n-b-1}\otimes J_{b+1}^b) \cap (I_{n-b}\otimes R_b)}
\\
& \simeq \frac{V^{(n-a-1)}\otimes J_{a+1}^a}{I_{n-a-1}\otimes J_{a+1}^a} \oplus \frac{V^{(n-b-1)}\otimes J_{b+1}^b}{I_{n-b-1}\otimes J_{b+1}^b}\simeq A_{n-a-1} \otimes J_{a+1}^a \oplus A_{n-b-1} \otimes J_{b+1}^b,
\end{align*}

\noindent where the isomorphism for the term involving $R_b$ is analogous to the previous case.
\end{itemize}

We then conclude that $\ker \delta_2$ is $2$-pure in degrees $a+1$ and $b+1$.

Conversely, if $\ker \delta_2$ is $2$-pure in degrees $a+1$ and $b+1$, then $(\ker \delta_2)_n= A_{n-a-1} \otimes J_{a+1}^a \oplus A_{n-b-1} \otimes J_{b+1}^b$ and using the equivalence mentioned at the beginning of the proof, \eqref{eq:econa} holds.

Assume first that $n>b$, then $(\ker \delta_2)_n \subseteq A_{n-a}\otimes R_a \oplus A_{n-b}\otimes R_b$. We know that
\begin{equation}\label{eq:Nn}
\begin{split}
\frac{N_n}{I_{n-a}\otimes R_a\oplus I_{n-b}\otimes R_b} & \simeq A_{n-a-1}\otimes J_{a+1}^a\oplus A_{n-b-1} \otimes J_{b+1}^b
\\
& \simeq \frac{V^{(n-a-1)}\otimes J_{a+1}^a}{I_{n-a-1}\otimes J_{a+1}^a} \oplus \frac{V^{(n-b-1)}\otimes J_{b+1}^b}{I_{n-b-1}\otimes J_{b+1}^b}.
\end{split}
\end{equation}

On the other hand, the following inclusions are trivial:
\begin{align*}
V^{(n-a-1)}\otimes J_{a+1}^a & \subseteq I_{n-a}\otimes R_a+ V^{(n-a-1)}\otimes J_{a+1}^a,
\\
I_{n-a-1}\otimes J_{a+1}^a & \subseteq I_{n-a}\otimes R_a,
\\
V^{(n-b-1)}\otimes J_{b+1}^b & \subseteq I_{n-b}\otimes R_b+ V^{(n-b-1)}\otimes J_{b+1}^b,
\\
I_{n-b-1}\otimes J_{b+1}^b & \subseteq I_{n-b}\otimes R_b.
\end{align*}

Hence, by Lemma \ref{lem:intersection} we get the equalities:
\begin{align*}
(V^{(n-a-1)}\otimes J_{a+1}^a)\cap (I_{n-a}\otimes R_a) & = I_{n-a-1}\otimes J_{a+1}^a,
\\
(V^{(n-b-1)}\otimes J_{b+1}^b)\cap (I_{n-b}\otimes R_b) & = I_{n-b-1}\otimes J_{b+1}^b.
\end{align*}

\noindent So, the last expression in \eqref{eq:Nn} is equivalent to $\frac{V^{(n-a-1)}\otimes J_{a+1}^a}{(V^{(n-a-1)}\otimes J_{a+1}^a) \cap (I_{n-a}\otimes R_a)}\oplus \frac{V^{(n-b-1)}\otimes J_{b+1}^b}{(V^{(n-b-1)}\otimes J_{b+1}^b)\cap (I_{n-b}\otimes R_b)}$, which is in term isomorphic to $\frac{(V^{(n-a-1)}\otimes J_{a+1}^a)+(I_{n-a}\otimes R_a)}{I_{n-a}\otimes R_a}\oplus \frac{(V^{(n-b-1)}\otimes J_{b+1}^b)+(I_{n-b}\otimes R_b)}{I_{n-b}\otimes R_b}$. By Lemma \ref{lem:samemod},
\begin{align*}
N_n & = V^{(n-a-1)}\otimes J_{a+1}^a + I_{n-a}\otimes R_a \oplus V^{(n-b-1)}\otimes J_{b+1}^b + I_{n-b}\otimes R_b
\\
& = (E'\cap F') + (E'\cap G') \oplus (E''\cap F'') + (E''\cap G'').
\end{align*}

\noindent As a consequence, the 4-tuple is multidistributive and, analogously, the triple is distributive. \qed
\end{prof}


\subsection{\texorpdfstring{Description of $\ker \delta_3$}{sec:kerdelta3}}

From now on we assume that $\ker \delta_2$ is $2$-pure in degrees $a+1$ and $b+1$. The canonical injection $\tilde{g}_3:J_{a+1}^a\oplus J_{b+1}^b\rightarrow \ker \delta_2$ induces a projective cover $g_3:(A\otimes J_{a+1}^a)\oplus (A\otimes J_{b+1}^b) \rightarrow \ker \delta_2$. Let $\delta_3$ be the composition of this map with $inc_2:\ker \delta_2 \rightarrow A\otimes R$.

Now, our aim is to describe $\ker \delta_3$. Note first that if $\overline{\alpha_i}, \overline{\alpha '_j}\in A$, $x_i, x'_j\in V$, $y_i\in R_a$ and $y'_j\in R_b$ are such that $\delta_3 (\sum_i \overline{\alpha_i}\otimes x_i\otimes y_i+ \sum_j \overline{\alpha '_j}\otimes x'_j\otimes y'_j)=0$, then  $\sum_i \overline{\alpha_i x_i}\otimes y_i=\sum_j \overline{\alpha '_j x'_j}\otimes y'_j=0$, since $R_a$ and $R_b$ are exclusive. So we can describe each direct summand of $\ker \delta_3$.

Note that $\delta_3$ can vanish only on elements such that $\sum_i \overline{\alpha_i x_i}=0$ and $\gr(\sum_i \overline{\alpha_i x_i})$ $\ge a$, hence $\gr (\sum_i \overline{\alpha_i} \otimes x_i \otimes y_i)\ge 2a$. As a consequence $(\ker \delta_3)_n=0$ if $n<2a$.

It is straightforward (see \cite[\S 4.2.3]{R}) that
\[
(\ker \delta_3)_{2a} = (V^{(a-1)}\otimes J_{a+1}^a) \cap (R_a\otimes V^{(a)})=J_{2a}^a.
\]

For $n$ such that $2a<n<a+b$, $(\ker \delta_3)_n\subseteq A_{n-a-1}\otimes J_{a+1}^a\oplus A_{n-b-1}\otimes J_{b+1}^b$. In fact, $(\ker \delta_3)_n\cap (A_{n-b-1}\otimes J_{b+1}^b)=0$ since $R_a$ and $R_b$ are exclusive and $n-b-1<a$. Therefore
\[
(\ker \delta_3)_n = \frac{(V^{(n-a-1)}\otimes J_{a+1}^a)\cap (I_{n-a-1}\otimes V^{(a+1)}+V^{(n-2a)}\otimes R_a\otimes V^{(a)})}{I_{n-a-1}\otimes J_{a+1}^a}.
\]

\medskip

If $a+b\le n<2b$, an inspection of the kernel gives
\[
(\ker \delta_3)_n = \frac{(V^{(n-a-1)}\otimes J_{a+1}^a\oplus V^{(n-b-1)}\otimes J_{b+1}^b)\cap (I_{n-a}\otimes V^{(a)}+I_{n-b}\otimes V^{(b)})}{I_{n-a-1}\otimes J_{a+1}^a\oplus I_{n-b-1}\otimes J_{b+1}^b}.
\]

\noindent Since $n-b<b$, $R_b$ does not appear in $I_{n-b}$ and the numerator equals
\begin{align*}
& (V^{(n-a-1)}\otimes J_{a+1}^a\oplus V^{(n-b-1)}\otimes J_{b+1}^b) \cap(I_{n-a-1}\otimes V^{(a+1)}+ V^{(n-2a)}\otimes R_a\otimes V^{(a)}+
\\
& V^{(n-a-b)}\otimes R_b\otimes V^{(a)}+ I_{n-b-1}\otimes V^{(b+1)}+V^{(n-a-b)}\otimes R_a\otimes V^{(b)}).
\end{align*}

\medskip

Finally, for $n\ge 2b$, we have
\[
(\ker \delta_3)_n = \frac{(V^{(n-a-1)}\otimes J_{a+1}^a\oplus V^{(n-b-1)}\otimes J_{b+1}^b)\cap (I_{n-a}\otimes V^{(a)}+ I_{n-b}\otimes V^{(b)})}{I_{n-a-1}\otimes J_{a+1}^a\oplus I_{n-b-1}\otimes J_{b+1}^b}.
\]

We shall denote by $N_n$ the numerator of $(\ker \delta_3)_n$.
\medskip

\begin{defi}
A 4-tuple $(E,F,G,H)$ is \textbf{distributive} if 
\[
E\cap (F+G+H) = (E\cap F) + (E\cap G) + (E\cap H).
\]
\end{defi}

\medskip

Whenever the following equalities hold we will say that the \textbf{extra vanishing conditions} \textit{(e.v.c.)} are satisfied
\begin{align}
(V^{(b-1)}\otimes R_a \otimes V) \cap (V^{(b)}\otimes R_a) \cap (R_b\otimes V^{(a)}) & = 0,
\\
(V^{(a-1)}\otimes R_b \otimes V) \cap (V^{(a)}\otimes R_b) \cap (R_a\otimes V^{(b)}) & = 0.
\end{align}

\bigskip

In case the \textit{(e.c.)} and the \textit{(e.v.c.)} are satisfied, we are able to derive properties of $\ker \delta_3$ from those of $\ker \delta_2$ as it is shown in the following proposition.

\begin{prop}\label{prop:kerdelta3dist}
Assume that $\ker \delta_2$ is $2$-pure in degrees $a+1$ and $b+1$, and that the \textit{(e.v.c.)} and the second relation of the \textit{(e.c.)} \eqref{eq:ec} are satisfied. Then, $\ker \delta_3$ is $2$-pure in degrees $2a$ and $2b$ if and only if for all $n\ge 2a+1$, $(E,F,G,H)$ and $(E,E',F+G+H,F'+G'+H')$ are respectively distributive and multidistributive, where
\begin{align*}
E & = V^{(n-a-1)} \otimes J_{a+1}^a, & \hskip 2cm E' & = V^{(n-b-1)} \otimes J_{b+1}^b,
\\
F & = I_{n-a-1}\otimes V^{(a+1)}, & \hskip 2cm F' & = I_{n-b-1}\otimes V^{(b+1)},
\\
G & = V^{(n-2a)}\otimes R_a \otimes V^{(a)}, & \hskip 2cm G'& = V^{(n-a-b)}\otimes R_a \otimes V^{(b)},
\\
H & = V^{(n-a-b)}\otimes R_b\otimes V^{(a)}, & \hskip 2cm H' & = V^{(n-2b)}\otimes R_b\otimes V^{(b)}.
\end{align*}
\end{prop}

\begin{prof}
In order to prove the ``if'' part suppose that the distributivity and multidistributivity of the tuples are satisfied. We have already seen that
\[
(\ker \delta_3)_n = 0 \hskip 0.2cm \text{if} \hskip 0.2cm n<2a, \text{ and } (\ker \delta_3)_{2a} = J_{2a}^a.
\]

Moreover, if $2a<n<a+b$, 
\begin{align*}
N_n & = E\cap (F+G)= (E\cap F)+(E\cap G)
\\
& = I_{n-a-1}\otimes J_{a+1}^a+ V^{(n-2a)}\otimes [(V^{(a-1)}\otimes R_a\otimes V) \cap (V^{(a)}\otimes R_a)\cap (R_a\otimes V^{(a)})]
\\
& = I_{n-a-1}\otimes J_{a+1}^a+V^{(n-2a)}\otimes J_{2a}^a.
\end{align*}

\noindent The equality $(\ker \delta_3)_n = A_{n-2a} \otimes J_{2a}^a$ follows using arguments similar to the proof of Proposition \ref{prop:kerdelta2pure}.

If $a+b\le n<2b$,
\begin{align*}
N_n = & (E\oplus E')\cap (F+G+H+F'+G')
\\
= & (E\cap F)+(E\cap G)+(E\cap H)\oplus (E'\cap F')+(E'\cap G')
\\
= & I_{n-a-1}\otimes J_{a+1}^a + V^{(n-2a)}\otimes J_{2a}^a+
\\
& V^{(n-a-b)}\otimes \underbrace{[(V^{(b-1)}\otimes R_a\otimes V)\cap (V^{(b)}\otimes R_a)\cap (R_b\otimes V^{(a)})]}_{=0 \textit{ (e.v.c.)}} \oplus
\\
& I_{n-b-1}\otimes J_{b+1}^b + V^{(n-b-a)}\otimes \underbrace{[(V^{(a-1)}\otimes R_b\otimes V)\cap (V^{(a)}\otimes R_b) \cap (R_a\otimes V^{(b)})]}_{=0 \textit{ (e.v.c.)}}
\\
= & I_{n-a-1}\otimes J_{a+1}^a + V^{(n-2a)}\otimes J_{2a}^a \oplus I_{n-b-1}\otimes J_{b+1}^b.
\end{align*}

\noindent Again, as in Proposition \ref{prop:kerdelta2pure}, $(\ker \delta_3)_n = A_{n-2a}\otimes J_{2a}^a$. Note that the class of any element in $I_{n-b-1}\otimes J_{b+1}^b$ vanishes in $(\ker \delta_3)_n$.

Finally, if $n\ge 2b$ we proceed in the same way and we use Lemma \ref{lem:equivalenceforincl}. Then,
\[
N_n = (I_{n-a-1}\otimes J_{a+1}^a+V^{(n-2a)}\otimes J_{2a}^a) \oplus (I_{n-b-1}\otimes J_{b+1}^b+V^{(n-2b)}\otimes J_{2b}^b).
\]

\noindent Hence, $\ker \delta_3$ is $2$-pure in degrees $2a$ and $2b$.

Conversely, if $\ker \delta_3$ is $2$-pure in degrees $2a$ and $2b$, 
\begin{align*}
(\ker \delta_3)_n & = 0 & \hskip 0.2cm & \text{if} \hskip 0.2cm n<2a,
\\
(\ker \delta_3)_n & = A_{n-2a}\otimes J_{2a}^a & \hskip 0.2cm & \text{if}\hskip 0.2cm 2a\le n<2b,
\\
(\ker \delta_3)_n & = A_{n-2a}\otimes J_{2a}^a\oplus A_{n-2b}\otimes J_{2b}^b & \hskip 0.2cm & \text{if} \hskip 0.2cm n\ge 2b.
\end{align*}

Mimicking the description of $\ker \delta_2$ and taking into account that $E\cap H=0=E'\cap H'$,
\begin{align*}
(\ker \delta_3)_n & = \frac{(E\oplus E')\cap (F+G+H+F'+G'+H')} {I_{n-a-1}\otimes J_{a+1}^a\oplus I_{n-b-1}\otimes J_{b+1}^b}
\\
& = A_{n-2a}\otimes J_{2a}^a \oplus A_{n-2b}\otimes J_{2b}^b = \frac{V^{(n-2a)}\otimes J_{2a}^a}{I_{n-2a}\otimes J_{2a}^a}\oplus \frac{V^{(n-2b)}\otimes J_{2b}^b}
{I_{n-2b}\otimes J_{2b}^b}
\\
& = \frac{E\cap G}{I_{n-2a}\otimes J_{2a}^a} \oplus \frac{E\cap H}{I_{n-a-b}\otimes R_b\otimes R_a} \oplus \frac{E'\cap H'}{I_{n-a-b}\otimes R_a\otimes R_b} \oplus \frac{E'\cap G'}{I_{n-2b}\otimes J_{2b}^b}
\\
& = \frac{E\cap G}{I_{n-2a}\otimes J_{2a}^a} \oplus \frac{E'\cap G'}{I_{n-2b}\otimes J_{2b}^b}.
\end{align*}

By Lemma \ref{lem:intersection} and since $E\cap F = I_{n-a-1}\otimes J_{a+1}^a$ and $I_{n-2a}\otimes J_{2a}^a \subseteq I_{n-a-1}\otimes J_{a+1}^a$, $(E\cap G)\cap (I_{n-a-1}\otimes J_{a+1}^a) = I_{n-2a}\otimes J_{2a}^a$. Thus,
\[
\frac{E\cap G}{I_{n-2a}\otimes J_{2a}^a}= \frac{E\cap G}{(E\cap G) \cap (I_{n-a-1}\otimes J_{a+1}^a)} \simeq \frac{(E\cap G)+(E\cap F)}{I_{n-a-1}\otimes J_{a+1}^a}.
\]

\noindent The other cases are analogous. It follows from Lemma \ref{lem:samemod} that the tuples are respectively distributive and multidistributive. \qed
\end{prof}


\subsection{\texorpdfstring{Description of $\ker \delta_i$ for $i>3$}{sec:kerdeltai}}\label{sec:kerdeltai}

The study of the subspace $\ker \delta_i$ for $i>3$ will be related to the following conditions.
\begin{equation}\label{eq:firstecc}
(V^{(b-1)}\otimes R_a) \cap (R_b\otimes V^{(a-1)}+\cdots +V^{(a-2)}\otimes R_b\otimes V)=0,
\end{equation}
\begin{equation}\label{eq:secondecc}
(V^{(a-1)}\otimes R_b)\cap (R_a\otimes V^{(b-1)}+ \cdots + V^{(b-2)}\otimes R_a\otimes V)=0.
\end{equation}

\noindent We shall call them \textbf{extra crossed conditions}, denoted \textit{(e.c.c.)}. Note that if \eqref{eq:firstecc} holds, then $(V^{(t)}\otimes R_a)\cap (R_b\otimes V^{(t+a-b)}+\cdots +V^{(t+a-b-1)}\otimes R_b\otimes V)=0$ for $b-a+1\le t\le b-2$. Analogously, if \eqref{eq:secondecc} is satisfied, then $(V^{(t)}\otimes R_b)\cap (R_a\otimes V^{(t-a+b)}+\cdots +V^{(t-a+b-1)}\otimes R_a\otimes V)=0$ for $b-a+1\le t\le a-2$.
 
\bigskip

The arguments used to compute $\ker \delta_i$ will be similar to those of previous subsections. Recall that we have fixed a basis $\mathcal{B}$ of $V$ which induces bases of the spaces $V^{(m)}$ for $m\ge 2$. 

For each $s\in \mathbb{N}$, we define the map $n_s:\mathbb{Z}_{\ge 0}\rightarrow \mathbb{Z}_{\ge 0}$ as follows
\[
n_s(2l)=la, \hskip 0.4cm n_s(2l+1)=la+1.
\]

Given $i>3$, suppose that $\delta_0, \cdots ,\delta_{i-1}$ have been defined in such a way that there are injections
\begin{equation}
\tilde{g}_j:J_{n_a(j)}^a\oplus J_{n_b(j)}^b\longrightarrow \ker \delta_{j-1} \text{ for all } 0<j\le i,
\end{equation}
inducing $g_j:A\otimes J_{n_a(j)}^a\oplus A\otimes J_{n_b(j)}^b\longrightarrow \ker \delta_{j-1}$ and there are inclusions
\begin{equation}
inc_j:\ker \delta_j \longrightarrow A\otimes J_{n_a(j)}^a\oplus A\otimes J_{n_b(j)}^b \text{ for all } 0\le j\le i-1.
\end{equation}

\noindent Then $\delta_i$ may be defined as the composition of $g_i$ and $inc_{i-1}$. Explicitly,
\[
\delta_i:A\otimes J_{n_a(i)}^a\oplus A\otimes J_{n_b(i)}^b \longrightarrow A\otimes J_{n_a(i-1)}^a\oplus A\otimes J_{n_b(i-1)}^b
\]
is such that, for $s=a,b$
\begin{equation}
\delta_i(\overline{\alpha}\otimes v_{j_1}\cdots v_{j_{n_s(i)}})=\begin{cases}
\overline{\alpha v_{j_1}\cdots v_{j_{s-1}}}\otimes v_{j_s}\cdots v_{j_{n_s(i)}} & \text{if $i$ is even,}
\\
\overline{\alpha v_{j_1}}\otimes v_{j_2}\cdots v_{j_{n_s(i)}} & \text{if $i$ is odd.}
\end{cases}
\end{equation}

The following statements are straightforward (see \cite{R}, \S 4.2.4, for details).

\medskip

$\bullet$ \textbf{Case $i$ even, $i\ge 4$}

By inspection, $(\ker \delta_i)_n=0$ if $n\le n_a(i)$ and $(\ker \delta_i)_{n_a(i)+1}= (V\otimes J_{n_a(i)}^a)\cap (R_a\otimes J_{n_a(i-1)}^a)=J_{n_a(i)+1}^a$.

Let $n=n_a(i)+m$ with $2\le m\le a-1$. We consider the cases:
\begin{itemize}
\item[$-$] $m$ such that $m<b-a+1$, then $R_b$ does not appear in $I_{m+a-1}$ and $(\ker \delta_i)_n=V^{(m-1)}\otimes J_{n_a(i)+1}^a$. 
\item[$-$] $m$ such that $b-a+1\ge m<n_{b-a}(i)$, then
\begin{align*}
(\ker \delta_i)_n = & (V^{(m)}\otimes J_{n_a(i)}^a) \cap [(R_a\otimes V^{(m-1)}+\cdots + V^{(m-1)}\otimes R_a)\otimes V^{(n_a(i-1))}+
\\
& (R_b\otimes V^{(m+a-b-1)}+\cdots +V^{(m+a-b-1)}\otimes R_b)\otimes V^{(n_a(i-1))}].
\end{align*}
\item[$-$] If $m=n_{b-a}(i)$, then $n=n_b(i)$ and
\begin{align*}
(\ker \delta_i)_{n_b(i)} = & (V^{(n_{b-a}(i))}\otimes J_{n_a(i)}^a \oplus J_{n_b(i)}^b) \cap [(R_a\otimes V^{(n_{b-a}(i)-1)}+\cdots + V^{(n_{b-a}(i)-1)}\otimes R_a)\otimes 
\\
& V^{(n_a(i-1))}+(R_b\otimes V^{(n_{b-a}(i-2))}+\cdots +V^{(n_{b-a}(i-2))}\otimes R_b)\otimes V^{(n_a(i-1))}].
\end{align*}
\item[$-$] For $m>n_{b-a}(i)$, 
\[
(\ker \delta_i)_n=(V^{(m)}\otimes J_{n_a(i)}^a\oplus V^{(n-n_b(i))}\otimes J_{n_b(i)}^b) \cap (I_{n-n_a(i-1)}\otimes V^{(n_a(i-1))}+I_{n-n_b(i-1)}\otimes V^{(n_b(i-1))}).
\]
\end{itemize}

Finally, for $n\ge n_a(i)+a$, then
\begin{align*}
(\ker \delta_i)_n & = \frac{N_n}{I_{n-n_a(i)}\otimes J_{n_a(i)}^a} & \text{if} \hskip 0.2cm n<a+n_b(i),
\\
(\ker \delta_i)_n & = \frac{N_n}{I_{n-n_a(i)}\otimes J_{n_a(i)}^a\oplus I_{n-n_b(i)}\otimes J_{n_b(i)}^b} & \text{if} \hskip 0.2cm n\ge a+n_b(i),
\end{align*}

\noindent where
\[
N_n = (V^{(n-n_a(i))}\otimes J_{n_a(i)}^a\oplus V^{(n-n_b(i))}\otimes J_{n_b(i)}^b) \cap (I_{n-n_a(i-1)}\otimes V^{(n_a(i-1))} +I_{n-n_b(i-1)}\otimes V^{(n_b(i-1))}).
\]

\medskip

$\bullet$ \textbf{Case $i$ odd}

In this case $(\ker \delta_i)_n=0$ for $n<n_a(i+1)$ and $(\ker \delta_i)_{n_a(i+1)}=J_{n_a(i+1)}^a$.

For $n\ge n_a(i+1)+1=n_a(i+2)$,
\begin{align*}
(\ker \delta_i)_n & = \frac{N_n}{I_{n-n_a(i)}\otimes J_{n_a(i)}^a} & \text{if} \hskip 0.2cm n<a+n_b(i),
\\
(\ker \delta_i)_n & = \frac{N_n}{I_{n-n_a(i)}\otimes J_{n_a(i)}^a\oplus I_{n-n_b(i)}\otimes J_{n_b(i)}^b} & \text{if} \hskip 0.2cm n\ge a+n_b(i),
\end{align*}

\noindent where
\begin{itemize}
\item[$-$] $N_n = (V^{(n-n_a(i))}\otimes J_{n_a(i)}^a) \cap (\sum_{j+h=n-n_a(i+2)}V^{(j)}\otimes R_a\otimes V^{(h)}\otimes V^{(n_a(i))}+ V^{(n-n_a(i+1))}\otimes R_a\otimes V^{(n_a(i-1))})$, if $n<\min \{n_b(i),n_a(i-1)+b\}$.
\item[$-$] $N_n = (V^{(n-n_a(i))}\otimes J_{n_a(i)}^a) \cap (\sum_{j+h=n-n_a(i+2)}V^{(j)}\otimes R_a\otimes V^{(h)}\otimes V^{(n_a(i))} +\sum_{p+q=n-n_a(i)-b}V^{(p)}\otimes R_b\otimes V^{(q)}\otimes V^{(n_a(i))}+V^{(n-n_a(i+1))}\otimes R_a\otimes V^{(n_a(i-1))}+ V^{(n-n_a(i-1)-b)}\otimes R_b\otimes V^{(n_a(i-1))})$, if $n_a(i-1)+b\le n<n_b(i)$.

Note that if $n=n_a(i-1)+b$ then the term $\sum_{p+q=n-n_a(i)-b}V^{(p)}\otimes R_b\otimes V^{(q)}\otimes V^{(n_a(i))}$ vanishes.

\item[$-$] $N_n =(V^{(n-n_a(i))}\otimes J_{n_a(i)}^a \oplus V^{(n-n_b(i))}\otimes J_{n_b(i)}^b)\cap (\sum_{j+h=n-n_a(i+2)}V^{(j)}\otimes R_a\otimes V^{(h)}\otimes V^{(n_a(i))}+\sum_{p+q=n-n_a(i)-b}V^{(p)}\otimes R_b\otimes V^{(q)}\otimes V^{(n_a(i))}+ V^{(n-n_a(i+1))}\otimes R_a\otimes V^{(n_a(i-1))}+V^{(n-n_a(i-1)-b)}\otimes R_b\otimes V^{(n_a(i-1))}+\sum_{l+m=n-a-n_b(i)}V^{(l)}\otimes R_a\otimes V^{(m)}\otimes V^{(n_b(i))}+\sum_{s+t=n-n_b(i+2)}V^{(s)}\otimes R_b\otimes V^{(t)} \otimes V^{(n_b(i))}+V^{(n-a-n_b(i-1))}\otimes R_a \otimes V^{(n_b(i-1))}+V^{(n-n_b(i+1))}\otimes R_b\otimes V^{(n_b(i-1))})$, if $n\ge n_b(i)$.
\end{itemize}

\medskip

We will consider the spaces:
\begin{itemize}
\item[$\bullet$] $E_a= V^{(n-n_a(i))}\otimes J_{n_a(i)}^a$;
\item[$\bullet$] $E_b= V^{(n-n_b(i))}\otimes J_{n_b(i)}^b$;
\item[$\bullet$] $F_1= (R_a\otimes V^{(n-n_a(i+2))}+\cdots +V^{(n-n_a(i+2))}\otimes R_a)\otimes V^{(n_a(i))}$;
\item[$\bullet$] $F_2= (R_b\otimes V^{(n-n_a(i)-b)}+\cdots +V^{(n-n_a(i)-b)}\otimes R_b)\otimes V^{(n_a(i))}$;
\item[$\bullet$] $F_3= (R_a\otimes V^{(n-a-n_b(i))}+\cdots +V^{(n-a-n_b(i))}\otimes R_a)\otimes V^{(n_b(i))}$;
\item[$\bullet$] $F_4= (R_b\otimes V^{(n-n_b(i+2))}+\cdots +V^{(n-n_b(i+2))}\otimes R_b)\otimes V^{(n_b(i))}$;
\item[$\bullet$] $F_5= V^{(n-n_a(i+2))}\otimes R_a\otimes V^{(n_a(i)-1)}+\cdots +V^{(n-n_a(i)-1)}\otimes R_a\otimes V^{(n_a(i-2))}$ if $i$ is even;
\item[$\bullet$] $F_6= V^{(n-n_b(i+2)+1)}\otimes R_b\otimes V^{(n_b(i)-1)}+\cdots +V^{(n-n_b(i)-1)}\otimes R_b\otimes V^{(n_b(i-2)+1)}$ if $i$ is even;
\item[$\bullet$] $F_7= V^{(n-n_a(i)-b+1)}\otimes R_b\otimes V^{(n_a(i)-1)}+\cdots +V^{(n-n_a(i-2)-b-1)}\otimes R_b\otimes V^{(n_a(i-2)+1)}$ if $i$ is even;
\item[$\bullet$] $F_8= V^{(n-a-n_b(i)+1)}\otimes R_a\otimes V^{(n_b(i)-1)}+\cdots +V^{(n-a-n_b(i-2)-1)}\otimes R_a\otimes V^{(n_b(i-2)+1)}$ if $i$ is even;
\item[$\bullet$] $F_9=(V^{(n-n_a(i)-1)}\otimes R_a+\cdots +R_a\otimes V^{(n-n_a(i)-1)})\otimes V^{(n_a(i-1))}$ if $i$ is even;
\item[$\bullet$] $F_{10}=(V^{(n-n_a(i-1)-b)}\otimes R_b+\cdots +R_b\otimes V^{(n-n_a(i-1)-b)})\otimes V^{(n_a(i-1))}$ if $i$ is even;
\item[$\bullet$] $F_{11}=(V^{(n-a-n_b(i-1))}\otimes R_a+\cdots +R_a\otimes V^{(n-a-n_b(i-1))})\otimes V^{(n_b(i-1))}$ if $i$ is even;
\item[$\bullet$] $F_{12}=(V^{(n-n_b(i)-1)}\otimes R_b+\cdots +R_b\otimes V^{(n-n_b(i)-1)})\otimes V^{(n_b(i-1))}$ if $i$ is even;
\item[$\bullet$] $F_{13}= V^{(n-n_a(i+1))}\otimes R_a\otimes V^{(n_a(i-1))}$ if $i$ is odd;
\item[$\bullet$] $F_{14}= V^{(n-n_a(i-1)-b)}\otimes R_b\otimes V^{(n_a(i-1))}$ if $i$ is odd;
\item[$\bullet$] $F_{15}= V^{(n-a-n_b(i-1))}\otimes R_a\otimes V^{(n_b(i-1))}$ if $i$ is odd;
\item[$\bullet$] $F_{16}= V^{(n-n_b(i+1))}\otimes R_b\otimes V^{(n_b(i-1))}$ if $i$ is odd.
\end{itemize}

Then for $n\ge n_a(i)$, $N_n$ equals:
\begin{itemize}
\item[] $(E_a\oplus E_b)\cap (F_9+F_{10}+F_{11}+F_{12})$ if $i$ is even and $n\le n_a(i+2)-1$,
\item[] $(E_a\oplus E_b)\cap (F_1+F_2+F_3+F_4+F_5+F_6+F_7+F_8)$ if $i$ is even and $n\ge n_a(i+2)$,
\item[] $(E_a\oplus E_b)\cap (F_1+F_2+F_3+F_4+F_{13}+F_{14}+F_{15}+F_{16})$ if $i$ is odd.
\end{itemize}

\bigskip

\begin{thm}\label{thm:kerdeltai}
Given $E_a$, $E_b$ and $F_l$, $1\le l \le 16$, as before and $i\ge 4$, assume that for all $j$ such that $j<i$, $\ker \delta_j$ is $2$-pure in degrees:
\[
\begin{cases}
n_a(j)+1 \hskip 0.2cm \text{and}\hskip 0.2cm n_b(j)+1 & \text{if}\hskip 0.2cm j\hskip 0.2cm \text{is even},
\\
n_a(j+1) \hskip 0.2cm \text{and} \hskip 0.2cm n_b(j+1) & \text{if}\hskip 0.2cm j \hskip 0.2cm \text{is odd},
\end{cases}
\]

\noindent and that the \textit{(e.c.)}, the \textit{(e.v.c.)} and the \textit{(e.c.c.)} hold. Then, $\ker \delta_i$ is $2$-pure in degrees:
\[
\begin{cases}
n_a(i)+1 \hskip 0.2cm \text{and} \hskip 0.2cm n_b(i)+1 & \text{if} \hskip 0.2cm i \hskip 0.2cm \text{is even},
\\
n_a(i+1) \hskip 0.2cm \text{and} \hskip 0.2cm n_b(i+1) & \text{if} \hskip 0.2cm i \hskip 0.2cm \text{is odd},
\end{cases}
\]

\noindent if and only if the following conditions are satisfied:
\begin{itemize}
\item[$\bullet$] If $i$ is even:
\begin{itemize}
\item[\textit{(i)}] $(E_a, E_b, F_9+F_{10}, F_{11}+F_{12})$ is multidistributive for $n=n_a(i)+m$, with $1\le m\le a-1$,
\item[\textit{(ii)}] $(E_a, E_b, F_1+F_2+F_5+F_7,F_3+F_4+F_6+F_8)$ is multidistributive for all $n\ge n_a(i+2)$.
\end{itemize}
\item[$\bullet$] If $i$ is odd:
\begin{itemize}
\item[] $(E_a, E_b,F_1+F_2+F_{13}+F_{14}, F_3+F_4+F_{15}+F_{16})$ is multidistributive for all $n\ge n_a(i+2)$.
\end{itemize}
\end{itemize}
\end{thm}

\begin{prof}

\textbf{Case 1:} $i$ even. We know that
\[ 
(\ker \delta_i)_n=0 \hskip 0.2cm \text{if} \hskip 0.2cm n\le n_a(i), \hskip 1.0cm (\ker \delta_i)_{n_a(i)+1}=J_{n_a(i)+1}^a.
\]

For $n=n_a(i)+m$, with $2\le m\le a-1$, the following holds:
\begin{itemize}
\item[$\bullet$] If $m<b-a+1$, then
\[
(\ker \delta_i)_n = V^{(m-1)}\otimes J_{n_a(i)+1}^a = A_{m-1} \otimes J_{n_a(i)+1}^a.
\]
\item[$\bullet$] If $b-a+1\le m<n_{b-a}(i)$, then
\begin{equation}
(\ker \delta_i)_n = E_a\cap (F_9+F_{10})=_{\text{dist.}} (E_a\cap F_9)+(E_a\cap F_{10}),
\end{equation}

\noindent but
\begin{align*}
E_a\cap F_{10} & = (V^{(m)}\otimes J_{n_a(i)}^a)\cap [(V^{(m+a-b-1)}\otimes R_b+\cdots +R_b\otimes V^{(m+a-b-1)})\otimes V^{(n_a(i-1))}]
\\
& \subseteq [(V^{(m)}\otimes R_a)\cap (R_b\otimes V^{(m+a-b)}+\cdots +V^{(m+a-b-1)}\otimes R_b\otimes V)]\otimes V^{(n_a(i-2))}
\\
& =_{\text{\textit{(e.c.c.)}}} 0,
\end{align*}

\noindent so $(\ker \delta_i)_n=E_a\cap F_9=V^{m-1} \otimes J_{n_a(i)+1}^a$.
\item[$\bullet$] If $m\ge n_{b-a}(i)$, it follows from the \textit{(e.c.c.)} that $E_b\cap F_{11}=0$, so:
\begin{align*}
(\ker \delta_i)_n & = (E_a\oplus E_b)\cap (F_9+F_{10}+F_{11}+F_{12})
\\
& =_{\text{multidist.}} (E_a\cap F_9+E_a\cap F_{10}) \oplus (E_b\cap F_{11} + E_b\cap F_{12})
\\
& = V^{(m-1)}\otimes J_{n_a(i)+1}^a\oplus V^{(n-n_b(i)-1)}\otimes J_{n_b(i)+1}^b
\\
& = A_{m-1} \otimes J_{n_a(i)+1}^a\oplus A_{n-n_b(i)-1} \otimes J_{n_b(i)+1}^b.
\end{align*}

For $n\ge n_a(i)+a$, note first that
\begin{itemize}
\item[$-$] $E_a\cap (F_1+F_2) =I_{n-n_a(i)}\otimes J_{n_a(i)}^a$.
\item[$-$] $E_b\cap (F_3+F_4) = I_{n-n_b(i)}\otimes J_{n_b(i)}^b$.
\item[$-$] $F_1+F_2+F_5+F_7 = I_{n-n_a(i-1)}\otimes V^{(n_a(i-1))}$.
\item[$-$] $F_3+F_4+F_6+F_8 = I_{n-n_b(i-1)}\otimes V^{(n_b(i-1))}$.
\item[$-$] $E_a\cap F_5 =(V^{(n-n_a(i))}\otimes J_{n_a(i)}^a)\cap (V^{(n-n_a(i+2)+1)}\otimes R_a\otimes V^{(n_a(i)-1)}+\cdots + V^{(n-n_a(i)-1)}\otimes R_a\otimes V^{(n_a(i-1))})=V^{(n-n_a(i+2)+1)}\otimes [(V^{(a-1)}\otimes J_{n_a(i)}^a)\cap (R_a\otimes V^{(n_a(i)-1)}+\cdots + V^{(a-2)}\otimes R_a\otimes V^{(n_a(i-1))})]\subseteq V^{(n-n_a(i+2)+1)}\otimes [(V^{(a-1)}\otimes R_a)\cap (R_a\otimes V^{(a-1)}+\cdots +V^{(a-2)}\otimes R_a\otimes V)]\otimes J_{n_a(i-2)}^a=_{\textit{(e.c.)}} V^{(n-n_a(i+2)+1)}\otimes V^{(a-2)}\otimes J_{a+1}^a\otimes J_{n_a(i-2)}^a= V^{(n-n_a(i)-1)}\otimes J_{n_a(i)+1}^a$.

\noindent Since the other inclusion is trivial, $E_a\cap F_5 = V^{(n-n_a(i)-1)}\otimes J_{n_a(i)+1}^a$.
\item[$-$] In the same way, $E_b\cap F_6 = V^{(n-n_b(i)-1)}\otimes J_{n_b(i)+1}^b$.
\item[$-$] $E_a\cap F_7 = (V^{(n-n_a(i))}\otimes J_{n_a(i)}^a)\cap (V^{(n-n_a(i)-b+1)}\otimes R_b\otimes V^{(n_a(i)-1)}+\cdots + V^{(n-n_a(i-1)-b)}\otimes R_b\otimes V^{(n_a(i-1))}) \subseteq V^{(n-n_a(i)-b+1)}\otimes [(V^{(b-1)}\otimes R_a)\cap (R_b\otimes V^{(a-1)}+\cdots +V^{(a-2)}\otimes R_b \otimes V)]\otimes V^{(n_a(i-2))} =0$, by the \textit{(e.c.c.)}.
\item[$-$] Analogously, $E_b\cap F_8 = 0$.
\end{itemize}
\end{itemize}

We then obtain the following equalities
\begin{align*}
N_n & = (E_a\oplus E_b)\cap (F_1+F_2+F_5+F_7+F_3+F_4+F_6+F_8)
\\
& =_{\text{multdist.}} E_a\cap (F_1+F_2)+ E_a\cap F_5 + E_a\cap F_7 + E_b\cap (F_3+F_4)+E_b\cap F_6+ E_b\cap F_8
\\
& =I_{n-n_a(i)}\otimes J_{n_a(i)}^a + V^{(n-n_a(i)-1)}\otimes J_{n_a(i)+1}^a+I_{n-n_b(i)}\otimes J_{n_b(i)}^b+V^{(n-n_b(i)-1)}\otimes J_{n_b(i)+1}^b.
\end{align*}

Thus, $(\ker \delta_i)_n=\frac{N_n}{I_{n-n_a(i)}\otimes J_{n_a(i)}^a\oplus I_{n-n_b(i)}\otimes J_{n_b(i)}^b}$ is isomorphic to
\begin{align*}
& \frac{I_{n-n_a(i)}\otimes J_{n_a(i)}^a+V^{(n-n_a(i)-1)}\otimes J_{n_a(i)+1}^a}{I_{n-n_a(i)}\otimes J_{n_a(i)}^a} \oplus \frac{I_{n-n_b(i)}\otimes J_{n_b(i)}^b+V^{(n-n_b(i)-1)}\otimes J_{n_b(i)+1}^b}{I_{n-n_b(i)}\otimes J_{n_b(i)}^b} \simeq
\\
& \frac{V^{(n-n_a(i)-1)}\otimes J_{n_a(i)+1}^a}{(I_{n-n_a(i)}\otimes J_{n_a(i)}^a)\cap (V^{(n-n_a(i)-1)}\otimes J_{n_a(i)+1}^a)} \oplus \frac{V^{(n-n_b(i)-1)}\otimes J_{n_b(i)+1}^b}{(I_{n-n_b(i)}\otimes J_{n_b(i)}^b)\cap (V^{(n-n_b(i)-1)}\otimes J_{n_b(i)+1}^b)}.
\end{align*}

Arguments similar to those used in the proof of Proposition \ref{prop:kerdelta2pure} imply that
\[
(\ker \delta_i)_n \simeq A_{n-n_a(i)-1}\otimes J_{n_a(i)+1}^a\oplus A_{n-n_b(i)-1}\otimes J_{n_b(i)+1}^b,
\]

\noindent and then $\ker \delta_i$ is $2$-pure in degrees $n_a(i)+1$ and $n_b(i)+1$.

Conversely, suppose that $\ker \delta_i$ is $2$-pure in degrees $n_a(i)+1$ and $n_b(i)+1$, i.e.
\[
(\ker \delta_i)_n = A_{n-n_a(i)-1} \otimes J_{n_a(i)+1}^a\oplus A_{n-n_b(i)-1} \otimes J_{n_b(i)+1}^b.
\]

\noindent Using the description of $(\ker \delta_i)_n$ and by Lemma \ref{lem:samemod}, \textit{(i)} and \textit{(ii)} are satisfied.

\medskip

\textbf{Case 2:} $i$ odd.

We already know that:
\[
(\ker \delta_i)_n = 0 \hskip 0.2cm \text{if} \hskip 0.2cm n\le n_a(i+1)-1, \hskip 0.6cm (\ker \delta_i)_{n_a(i+1)} = J_{n_a(i+1)}^a.
\]

If $n>n_a(i+1)$, note first that
\begin{itemize}
\item[$-$] $F_1+F_2+F_{13}+F_{14}+F_3+F_4+F_{15}+F_{16} = I_{n-n_a(i-1)}\otimes V^{(n_a(i-1))}+I_{n-n_b(i-1)}\otimes V^{(n_b(i-1))}$.
\item[$-$] $E_a\cap (F_1+F_2) = I_{n-n_a(i)}\otimes J_{n_a(i)}^a$.
\item[$-$] $E_b\cap (F_3+F_4) = I_{n-n_b(i)}\otimes J_{n_b(i)}^b$.
\item[$-$] $E_a\cap F_{13} = (V^{(n-n_a(i))}\otimes J_{n_a(i)}^a)\cap (V^{(n-n_a(i+1))}\otimes R_a\otimes V^{(n_a(i-1))}) = V^{(n-n_a(i+1))}\otimes [(V^{(a-1)}\otimes J_{n_a(i)}^a)\cap (R_a\otimes V^{(n_a(i-1))})] =_{\textit{(e.c.)}} V^{(n-n_a(i+1))}\otimes J_{n_a(i+1)}^a$.  
\item[$-$] $E_a\cap F_{14} =(V^{(n-n_a(i))}\otimes J_{n_a(i)}^a)\cap (V^{(n-n_a(i-1)-b)}\otimes R_b\otimes V^{(n_a(i-1))})= V^{(n-n_a(i-1)-b)}\otimes [(V^{(b-1)}\otimes J_{n_a(i)}^a)\cap (R_b\otimes V^{(n_a(i-1))})]\subseteq V^{(n-n_a(i-1)-b)}\otimes [(V^{(b-1)}\otimes R_a)\cap (R_b\otimes V^{(a-1)})]\otimes V^{(n_a(i-2))}\subseteq  V^{(n-n_a(i-1)-b)}\otimes [(V^{(b-1)}\otimes R_a)\cap
(R_b\otimes V^{(a-1)}+\cdots +V^{(a-1)}\otimes R_b)]\otimes V^{(n_a(i-2))}$ $=_{\text{\textit{(e.c.c.)}}} 0$.
\item[$-$] Analogously, $E_b\cap F_{15}=0$. 
\item[$-$] $E_b\cap F_{16}=V^{(n-n_b(i+1))}\otimes J_{n_b(i+1)}^b$ by the \textit{(e.c.)}.
\end{itemize}

Thus,
\begin{align*}
N_n = & (E_a\oplus E_b)\cap (F_1+F_2+F_{13}+F_{14}+F_3+F_4+F_{15}+F_{16})
\\
=& E_a\cap (F_1+F_2) + E_a\cap F_{13} + E_a\cap F_{14} + E_b\cap (F_3+F_4) + E_b\cap F_{15} + E_b\cap F_{16}
\\
=& I_{n-n_a(i)}\otimes J_{n_a(i)}^a+ V^{(n-n_a(i+1))}\otimes J_{n_a(i+1)}^a+I_{n-n_b(i)}\otimes J_{n_b(i)}^b+ V^{(n-n_b(i+1))}\otimes J_{n_b(i+1)}^b.
\end{align*}

From now on the argument is analogous to the case $i$ even. \qed
\end{prof}


\section{\texorpdfstring{$(a,b)$-Koszul algebras}{sec:abkoszulalg}}

In this section we define $(a,b)$-Koszul algebras.

\medskip

For $i\ge 0$, let $K_i= A\otimes J_{n_a(i)}^a + A\otimes J_{n_b(i)}^b$, where $J_0^a=J_0^b= k$, $J_1^a=J_1^b=V$ and the sum is direct if $i\ge 2$.

\begin{defi}
An algebra $A=T(V)/I$ is \textbf{$(a,b)$-homogeneous} if $I$ admits a set of homogeneous generators, which are in degrees $a$ and $b$.
\end{defi}

\begin{defi}\label{def:abKoszulalg}
Let $a,b\in \mathbb{Z}$ be such that $2\le a<b$ and let $R=R_a\oplus R_b$ be a space of relations for the two-sided ideal $I$ with $R_a$ and $R_b$ exclusive. We will say that an $(a,b)$-homogeneous algebra $A$ is \textbf{$(a,b)$-Koszul} if the graded vector space $\tor_i^A (k,k)$ is $2$-pure in degrees $n_a(i)$ and $n_b(i)$ for all $i\ge 2$.
\end{defi}

\begin{rmk}
For $a=b$ - omitting of course in this case the condition that $R_a$ and $R_b$ are exclusive - the previous definition coincides with the definition of $a$-Koszul given in \cite{B2}.
\end{rmk}

\medskip

We recall the following lemma without proof.

\begin{lema}\cite[Corollaire 2.7]{B4}\label{lem:torbound}
Let $W$ be a graded $k$-vector space concentrated in degree $1$. Consider $A=T(W)/J$ and $R$ a space of relations for $J$ such that $R=\sum_{n\ge s} R_n$ for some integer $s\ge 2$. Then, for all $i\ge 2$, the graded vector spaces $\tor_i^A (k,k)$ and $\Ext_A^i (k,k)$ are zero up to degrees $n_s(i)-1$ and $-(n_s(i)-1)$ respectively.
\end{lema}

\begin{thm}\label{thm:Kozulequiv}
Let $A=T(V)/I$ be an $(a,b)$-homogeneous algebra ($2\le a<b$) and let $R=R_a\oplus R_b$ be a space of relations for $I$ such that $R_a$ and $R_b$ are exclusive. The following statements are equivalent:
\begin{itemize}
\item[\textit{(i)}] $A$ is $(a,b)$-Koszul.
\item[\textit{(ii)}] The procedure described in \S \ref{sec:kerdeltai} gives a minimal pure projective resolution of $k$
\begin{equation}\label{eq:Koszulres}
\cdots \longrightarrow K_i \overset{\delta_i}{\longrightarrow} K_{i-1}\longrightarrow \cdots \longrightarrow K_1 \overset{\delta_1}{\longrightarrow} K_0 \overset{\epsilon}{\longrightarrow} k \longrightarrow 0
\end{equation}

\noindent in the category of left bounded graded left $A$-modules.
\item[\textit{(iii)}] The \textit{(e.c.)}, the \textit{(e.v.c.)} and the \textit{(e.c.c.)} are satisfied and for all $j\ge 1$
\begin{itemize}
\item[$\bullet$] the 4-tuple $(E_a, E_b, F_9+F_{10}, F_{11}+F_{12})$ is multidistributive for all $n<(j+1)a$,
\item[$\bullet$] the 4-tuple $(E_1, E_2, D_1+G_1+H_1,D_2+G_2+H_2)$ is multidistributive for all $n\ge (j+1)a$,
\item[$\bullet$] the 4-tuple $(E'_1,E'_2, D'_1+G'_1+H'_1,D'_2+G'_2+H'_2)$ is multidistributive for all $n\ge (j+1)a+1$,
\end{itemize}

\noindent where $E_a$, $E_b$, $F_h$ ($9\le h\le 12$) are those of Theorem \ref{thm:kerdeltai} and
\begin{align*}
E_1 & = V^{(n-ja)}\otimes J_{ja}^a, & \hskip 0.2cm E'_1 & = V^{(n-ja-1)}\otimes J_{ja+1}^a,
\\
E_2 & = V^{(n-jb)}\otimes J_{jb}^b, & \hskip 0.2cm E'_2 & = V^{(n-jb-1)}\otimes J_{jb+1}^b,
\\
D_1 & = I_{n-ja}\otimes V^{(ja)}, & \hskip 0.2cm D'_1 & = I_{n-ja-1}\otimes V^{(ja+1)},
\\
D_2 & = I_{n-jb}\otimes V^{(jb)}, & \hskip 0.2cm D'_2 & = I_{n-jb-1}\otimes V^{(jb+1)},
\\
G_1 & = V^{(n-(j+1)a+1)}\otimes I_{2a-2}^a\otimes V^{((j-1)a+1)}, & \hskip 0.2cm G'_1 & = V^{(n-(j+1)a)}\otimes R_a\otimes V^{(ja)},
\\
G_2 & = V^{(n-(j+1)b+1)}\otimes I_{2b-2}^b\otimes V^{((j-1)b+1)}, & \hskip 0.2cm G'_2 & = V^{(n-(j+1)b)}\otimes R_b\otimes V^{(jb)},
\\
H_1 & = V^{(n-ja-b+1)}\otimes I_{a+b-2}^b\otimes V^{((j-1)a+1)}, & \hskip 0.2cm H'_1 & = V^{(n-ja-b)}\otimes R_b\otimes V^{(ja)},
\\
H_2 & = V^{(n-a-jb+1)}\otimes I_{a+b-2}^a\otimes V^{((j-1)b+1)}, & \hskip 0.2cm H'_2 & = V^{(n-a-jb)}\otimes R_a\otimes V^{(jb)}.
\end{align*}
\end{itemize}
\end{thm}

\begin{prof}
Notice first that if
\[
\cdots \longrightarrow A\otimes Q_i\overset{d_i}{\longrightarrow} A\otimes Q_{i-1}\longrightarrow \cdots \longrightarrow A\otimes R\overset{d_2}{\longrightarrow} A\otimes V\overset{d_1}{\longrightarrow} A \overset{d_0}{\longrightarrow} k \longrightarrow 0
\]
is a minimal projective resolution, we apply the functor $k\otimes_A \place$ and obtain the complex
\[
\cdots \longrightarrow Q_i\overset{1_k\otimes_A d_i}{\longrightarrow} Q_{i-1}\longrightarrow \cdots \longrightarrow R\overset{1_k\otimes_A d_2}{\longrightarrow} V\overset{1_k\otimes_A d_1}{\longrightarrow} k \longrightarrow 0.
\]
It is easy to see that $1_k\otimes_A d_i=0$ and then $\tor_0^A(k,k)=k$, $\tor_1^A(k,k)=V$, $\tor_2^A(k,k)=R$ and $\tor_i^A(k,k)=Q_i$ for $i\ge 3$.

Since $J_{n_s(i)}^s \subseteq R_s\otimes V^{(n_s(i-2))}$ for $s=a,b$, it is clear that $\delta_{i-1}\delta_i=0$.

The $k$-modules $J_n^a$ and $J_n^b$ are projective and $A$ is $k$-flat, so the $A$-modules $K_i$ are projective. Moreover, they are $2$-pure except for $K_0$ and $K_1$ which are pure.

Given $i\ge 2$, by Lemma \ref{lem:torbound}, $\tor_i^A (k,k)$ may not vanish only in degrees greater than or equal to $n_a(i)$. It is then clear that there is a resolution of the form \eqref{eq:Koszulres} if and only if $\tor_i^A(k,k)$ is $2$-pure in degrees $n_a(i)$ and $n_b(i)$, since $\tor^A_i(k,k)= J_{n_a(i)}^a\oplus J_{n_b(i)}^b$.

Now, we consider the spaces $F_l$ , $1\le l\le 16$ as in Theorem \ref{thm:kerdeltai} and we note that
\begin{itemize}
\item[$\bullet$] if $i$ is even:
\begin{align*}
& E_1 = E_a, & \hskip 0.2cm & D_1 = F_1+F_2, & \hskip 0.2cm  & G_1 = F_5, & \hskip 0.2cm  & H_1 = F_7,
\\
& E_2 = E_b, & \hskip 0.2cm & D_2 = F_3+F_4, & \hskip 0.2cm & G_2 = F_6, & \hskip 0.2cm & H_2 = F_8;
\end{align*}
\item[$\bullet$] if $i$ is odd:
\begin{align*}
& E'_1 = E_a, & \hskip 0.2cm & D'_1 = F_1+F_2, & \hskip 0.2cm & G'_1 = F_{13}, & \hskip 0.2cm & H'_1 = F_{14},
\\
& E'_2 = E_b, & \hskip 0.2cm & D'_2 = F_3+F_4, & \hskip 0.2cm & G'_2 = F_{16}, & \hskip 0.2cm & H'_2 = F_{15}.
\end{align*}
\end{itemize}

\noindent It is straightforward that condition \textit{(iii)} is equivalent to Theorem \ref{thm:kerdeltai}. \qed
\end{prof}

\begin{rmk}\label{rmk:left=right}
In Definition \ref{def:abKoszulalg}, $k$ is considered as the trivial left $A$-module and $A$ is said to be left $(a,b)$-Koszul. If $k$ is the trivial right $A$-module and we take $K_i= J_{n_a(i)}^a \otimes A+ J_{n_b(i)}^b\otimes A$ and differentials analogous to $\delta_i$, $A$ is said to be right $(a,b)$-Koszul. It is clear by definition that $A$ is left $(a,b)$-Koszul if and only if it is right $(a,b)$-Koszul.
\end{rmk}

\medskip

Given graded $A$-modules $N$ and $N'$, we recall the following notation:
\begin{align*}
\Hom_A (N,N')_d & = \{ f\in \Hom_A (N,N') \hskip 0.2cm / \hskip 0.2cm f(N_i) \subseteq N'_{i+d}, \hskip 0.2cm i\in \mathbb{Z}\},
\\
\underline{\Hom}_A (N,N') & = \smashoperator{\bigoplus \limits_{d\in \mathbb{Z}}} \Hom_A (N,N')_d,
\\
\hom_A (N,N') & = \Hom_A (N,N')_0.
\end{align*}

If $N$ is left bounded, let $\mathcal{P}$ be a minimal projective resolution of $N$. We denote
\begin{align*}
\underline{\Ext}_A^i (N,N') & = H^i (\underline{\Hom}_A (\mathcal{P},N')),
\\
\ext_A^i (N,N') & = H^i (\hom_A (\mathcal{P}, N')).
\end{align*}

\medskip

Next we prove some equivalent conditions to $(a,b)$-Koszulity.

\begin{prop}
For a $\mathbb{Z}$-graded $k$-algebra $A$, the following conditions are equivalent:
\begin{itemize}
\item[\textit{(i)}] $A$ is $(a,b)$-Koszul.
\item[\textit{(ii)}] $\ext_A^i (k,k[-n])=0$ if $n\ne n_a(i)$ and $n\ne n_b(i)$.
\item[\textit{(iii)}] Let $M$ and $N$ be two graded $A$-modules concentrated in degrees $m$ and $n$ respectively. If $n\ne m+n_a(i)$ and $n\ne m+n_b(i)$, then $\ext_A^i(M,N)=0$.
\end{itemize}
\end{prop}

\begin{prof}
Since $k$ is a finitely generated $A$-module, $\Ext_A^{i,-n} (k,k) = \underline{\Ext}_A^{i,-n}(k,k) = \ext_A^i (k,k[-n])$ and also $\Ext_A^{i,-n} (k,k) = (\tor_{i,n}^A (k,k))^*$ (see \cite{B4}). It is then clear that \textit{(i)} is equivalent to \textit{(ii)}.

Since $k$ is concentrated in degree $0$, $k[-n]$ is concentrated in degree $-n$, so it is straightforward that \textit{(iii)} implies \textit{(ii)}.

If \textit{(i)} holds, then $k$ admits a graded projective resolution $\mathcal{P}=(P_i)_{i\ge 0}$ such that for all $i$, $P_i=A\otimes Q_i$ where $Q_0$ and $Q_1$ are concentrated in degree $0$ and $1$ respectively and for $i\ge 2$, $Q_i$ is a $k$-vector space $2$-concentrated in degrees $n_a(i)$ and $n_b(i)$. In order to prove \textit{(iii)},we can suppose that $m=0$ since $\ext_A^i (M,N) = \ext_A^i (M[-m],N[-m])$ and $M$, being concentrated, is a direct sum of copies of $k$. Without loss of generality, we may also assume that $M=k$. The following isomorphisms are natural
\[
\hom_A (P_j,N)=\hom_A (A\otimes Q_j,N) \simeq \hom_k (Q_j,N),
\]

\noindent and $\hom_k (Q_j,N)$ vanishes in degrees different from $n_a(j)$ or $n_b(j)$. Thus, the complex $\hom_A (\mathcal{P},N)$ may have nonzero terms only in degrees $n_a(j)$ and $n_b(j)$. So
\[
\ext_A^i (k,N) = H^i (\hom_A (\mathcal{P},N)) =\frac{\ker d_i}{\im d_{i-1}}=\begin{cases}
\hom_A(P_i,N) & \text{if $i=n_a(j)$ or $i=n_b(j)$,}
\\
0 & \text{otherwise.}
\end{cases}
\] \qed
\end{prof}

\bigskip

Given an $(a,b)$-Koszul algebra $A$, the complex \eqref{eq:Koszulres} is called the \textbf{Koszul resolution} of the left $A$-module $k$. Using Proposition \ref{prop:essentialineachdegree}, this resolution is minimal and projective.

In general, for an $(a,b)$-homogeneous algebra $A$, the complex \eqref{eq:Koszulres} is called the (left) \textbf{Koszul complex} of $A$. It is a generalization of the Koszul complex defined by Priddy in the quadratic case (\cite{M, P}).

\begin{prop}\label{prop:Koszulcomplexequiv}
Let $A$ be an $(a,b)$-homogeneous algebra ($2\le a<b$) such that $I$ admits a space of relations $R=R_a\oplus R_b$ where $R_a$ and $R_b$ are exclusive. Then $A$ is $(a,b)$-Koszul if and only if its left (or right) Koszul complex is exact in positive degrees.
\end{prop}

\begin{prof}
Assume first that the left Koszul complex is exact. Applying the functor $k\otimes_A \place$, $\tor_i^A (k,k)\simeq k\otimes_A K_i$, which is $2$-pure in degrees $n_a(i)$ and $n_b(i)$ for $i\ge 2$. Then $A$ is $(a,b)$-Koszul. The converse is a consequence of Theorem \ref{thm:Kozulequiv}. \qed
\end{prof}

\bigskip

It is well-known (see \cite{B4}) that the beginning of a minimal graded projective resolution of $k$ is 
\[
A\otimes R \overset{\delta_2}{\longrightarrow} A\otimes V \overset{\delta_1}{\longrightarrow} A \overset{\delta_0}{\longrightarrow} k \longrightarrow 0.
\]

Since the length of a minimal projective resolution of $k$ gives the global dimension of $A$, the following proposition is immediate:
 
\begin{prop}
Let $A=T(V)/I$ be an $(a,b)$-homogeneous algebra ($2\le a<b$) such that it admits a space of relations $R=R_a\oplus R_b$ where $R_a$ and $R_b$ are exclusive. If the global dimension of $A$ is $2$, then $A$ is $(a,b)$-Koszul.
\end{prop}


\subsection{\texorpdfstring{Distributive lattices}{sec:distributivelattices}}

The aim of this subsection is to give some criteria which will allow us to decide whether certain triples and tuples are distributive or multidistributive. The proofs are rather technical.

The set of all the $k$-subspaces of $V$ will be denoted by $\mathcal{L}(V)$. It is known that $(\mathcal{L}(V),\subseteq ,+, \cap)$ is a modular lattice; i.e. given $W_1,W_2,W_3\in \mathcal{L} (V)$, if $W_1\subseteq W_3$, then $W_1+(W_2\cap W_3)=(W_1+W_2)\cap W_3$.

We say that a sublattice $\mathcal{S}\subseteq \mathcal{L}(V)$ is \textbf{distributive} if $E\cap (F+G)=(E\cap F)+(E\cap G)$ for all $E,F,G\in \mathcal{S}$.

\begin{rmk}\label{rmk:zerointer}
It is immediate that if a sublattice $\mathcal{S}\subseteq \mathcal{L}(V)$ is distributive, then for all $E,F_1,\cdots ,F_t$ $\in \mathcal{S}$, $E\cap (F_1+\cdots +F_t)=(E\cap F_1)+\cdots +(E\cap F_t)$. Also, if $\{ E_1,E_2,E'\} \subseteq \mathcal{S}$ are such that $E_1\cap E'=E_2\cap E'=0$, then $(E_1+E_2)\cap E'=0$. 
\end{rmk}

\bigskip

Given $W_1,\cdots ,W_n$ subspaces of $V$, we shall denote by $\mathcal{T}$ the sublattice generated by $W_1,\cdots,W_n$. The following result is a first criterion for distributivity. 

\begin{prop}\cite{B2, BF, BGS1}\label{prop:distbasis}
$\mathcal{T}$ is distributive if and only if there exists a basis $\mathcal{B}$ of $V$ such that $\mathcal{B}_i=\mathcal{B}\cap W_i$ is a basis of $W_i$ for all $1\le i\le n$. In this case, we say that $\mathcal{B}$ distributes with respect to $W_1,\cdots ,W_n$.
\end{prop}

\medskip

\begin{lema}\label{lem:crossedincl}
If $(E, E',F_1+\cdots +F_t,G_1+\cdots +G_{t'})$ is multidistributive, then 
\begin{align*}
E\cap (G_1+\cdots +G_{t'}) & \subseteq E\cap (F_1+\cdots +F_t) \text{ and }
\\
E'\cap (F_1+\cdots +F_t) & \subseteq E'\cap (G_1+\cdots +G_{t'}).
\end{align*}
\end{lema}

\begin{prof}
Since the tuple is multidistributive,
\[
(E\oplus E')\cap (F_1+\cdots +F_t+G_1+\cdots +G_{t'}) = (E\cap F_1)+\cdots +(E\cap F_t)+(E'\cap G_1)+\cdots +(E'\cap G_{t'}).
\]

If $v\in E\cap (G_1+\cdots +G_{t'})$, then $v\in (E\oplus E')\cap (F_1+\cdots +F_t+G_1+\cdots +G_{t'})$. Since $E\cap E'=0$, $v\in (E\cap F_1)+\cdots +(E\cap F_t)\subseteq E\cap (F_1+\cdots +F_t)$. The other inclusion is analogous. \qed
\end{prof}

\medskip

\begin{rmk}\label{rmk:avoidcasesdist}
$(E, E', F_1+\cdots +F_t, G_1+\cdots G_{t'})$ is multidistributive if and only if $(E', E, G_1+\cdots +G_{t'}, F_1+\cdots +F_t)$ is multidistributive in $\mathcal{L}(V)$.
\end{rmk}

Next we prove the following lemma.

\begin{lema}\label{lem:generators}
Let $\mathcal{S}$ be a distributive sublattice of $\mathcal{L}(V)$. The following statements hold for $E,E',E_1,$ $E_2,E'_1, E'_2,F_i,G_j\in \mathcal{S}$, $1\le i\le t$, $1\le j \le t'$.
\begin{itemize}
\item[\textit{(i)}] If $(E,E',F_1,G_1+\cdots +G_{t'})$ and $(E,E',F_2,G_1+\cdots +G_{t'})$ are multidistributive, then $(E,E',F_1+F_2,G_1+\cdots +G_{t'})$ and $(E,E',F_1\cap F_2,G_1+\cdots +G_{t'})$ are also multidistributive.
\item[\textit{(ii)}] If $(E_1,E',F_1+\cdots +F_t,G_1+\cdots + G_{t'})$ and $(E_2,E',F_1+\cdots +F_t,G_1+\cdots +G_{t'})$ are multidistributive, then $(E_1+E_2,E',F_1+\cdots +F_t,G_1+\cdots + G_{t'})$ and $(E_1\cap E_2,E',F_1+\cdots +F_t,G_1+\cdots +G_{t'})$ are also multidistributive.
\end{itemize}
\end{lema}

\begin{prof}
\begin{itemize}
\item[\textit{(i)}] By hypothesis,
\[
(E\oplus E')\cap (F_i+G_1+\cdots +G_{t'}) = (E\cap F_i)+(E'\cap G_1)+\cdots + (E'\cap G_{t'}) \text{ for $i=1,2$.}
\]

\noindent Thus,
\begin{align*}
& (E\oplus E')\cap (F_1+F_2+G_1+\cdots +G_{t'})
\\
& =_{\text{dist.}} (E\oplus E') \cap (F_1+G_1+\cdots +G_{t'})+(E\oplus E')\cap (F_2+G_1+\cdots +G_{t'})
\\
& =(E\cap F_1)+(E'\cap G_1)+\cdots (E'\cap G_{t'}) + (E\cap F_2)+(E'\cap G_1)+\cdots +(E'\cap G_{t'})
\\
& = (E\cap F_1) + (E\cap F_2)+(E'\cap G_1)+\cdots +(E'\cap G_{t'}).
\end{align*}

It is clear that for $i=1,2$:
\begin{equation}\label{eq:F1}
\begin{split}
(E\oplus E')\cap ((F_1\cap F_2)+G_1+\cdots +G_{t'}) & \subseteq (E\oplus E')\cap (F_i+G_1+\cdots +G_{t'})
\\
&= (E\cap F_i)+(E'\cap G_1)+\cdots +(E'\cap G_{t'}).
\end{split}
\end{equation}

\noindent Let $v\in (E\oplus E')\cap ((F_1\cap F_2)+G_1+\cdots + G_{t'})$. By \eqref{eq:F1} there exist $x\in E\cap F_1$, $y\in E\cap F_2$ and $z,w\in (E'\cap G_1)+\cdots +(E'\cap G_{t'})$ such that $v=x+z=y+w$, then $x-y=w-z\in (E'\cap G_1)+\cdots +(E'\cap G_{t'})$. Since $x,y\in E$, $x-y\in E$ and so $x=y\in F_2$. Therefore, $x\in E\cap F_1\cap F_2$ and $(E,E',F_1\cap F_2,G_1+\cdots +G_{t'})$ is multidistributive.
\item[\textit{(ii)}] Note first that $(E_1+E_2)\cap E'=0$ by Remark \ref{rmk:zerointer} and that
\begin{align*}
&((E_1+E_2)\oplus E')\cap (F_1+\cdots +F_t+G_1+\cdots +G_{t'})
\\
& =_{\text{dist.}} (E_1\oplus E')\cap (F_1+\cdots +F_t+G_1+\cdots +G_{t'})+ (E_2\oplus E')\cap (F_1+\cdots +F_t+G_1+\cdots +G_{t'})
\\
& =_{\text{hip}}(E_1\cap F_1)+\cdots +(E_1\cap F_t)+(E'\cap G_1)+\cdots +(E'\cap G_{t'})+
\\
& (E_2\cap F_1)+\cdots +(E_2\cap F_t)+(E'\cap G_1)+\cdots +(E'\cap G_{t'})
\\
& =_{\text{dist.}} ((E_1+E_2)\cap F_1)+\cdots +((E_1+E_2)\cap F_t)+(E'\cap G_1)+\cdots +(E'\cap G_{t'}).
\end{align*}

For the other tuple, recall that $((E_1\cap E_2)\oplus E')\cap (F_1+\cdots +F_t+G_1+\cdots +G_{t'}) \subseteq (E_i\oplus E')\cap (F_1+\cdots +F_t+G_1+\cdots +G_{t'})= (E_i\cap F_1)+\cdots +(E_i\cap F_t)+(E'\cap G_1)+\cdots +(E'\cap G_{t'})$ for $i=1,2$. Thus,
\begin{align*}
& ((E_1\cap E_2)\oplus E')\cap (F_1+\cdots +F_t+G_1+\cdots +G_{t'})
\\
& \subseteq [(E_1\cap F_1)+\cdots +(E_1\cap F_t)+(E'\cap G_1)+\cdots +(E'\cap G_{t'})]\cap
\\
& [(E_2\cap F_1)+\cdots +(E_2\cap F_t)+(E'\cap G_1)+\cdots +(E'\cap G_{t'})]
\\
& =_{\text{dist.}} (E_1\cap E_2\cap F_1)+\cdots +(E_1\cap E_2\cap F_t) +\smashoperator{\sum \limits_{\substack{1\le i,j\le t
\\
i\ne j}}}(E_1\cap E_2 \cap F_i \cap F_j)+
\\
& \smashoperator{\sum \limits_{\substack{1\le i\le t
\\
1\le h\le t'}}}(E_1\cap F_i \cap E' \cap G_h)+\smashoperator{\sum \limits_{\substack{1\le h \le t'
\\
1\le i \le t}}}(E'\cap G_h \cap E_2 \cap F_i)+
\\
& (E'\cap G_1)+\cdots +(E'\cap G_{t'}) + \smashoperator{\sum \limits_{\substack{1\le h,l\le t'
\\
h\ne l}}}(E' \cap G_h \cap G_l)
\\
& \subseteq (E_1\cap E_2\cap F_1)+\cdots +(E_1\cap E_2\cap F_t)+(E'\cap G_1)+\cdots +(E'\cap G_{t'}).
\end{align*}
\end{itemize}
\end{prof}

Using Remark \ref{rmk:avoidcasesdist}, we obtain the following

\begin{coro}
\item[\textit{(i)}] If $(E,E',F_1+\cdots +F_t,G_1)$ and $(E,E',F_1+\cdots +F_t,G_2)$ are multidistributive, then $(E,E',F_1+\cdots +F_t,G_1+G_2)$ and $(E,E',F_1+\cdots +F_t,G_1\cap G_2)$ are also multidistributive.
\item[\textit{(ii)}] If $(E,E'_1,F_1+\cdots +F_t,G_1+\cdots + G_{t'})$ and $(E,E'_2,F_1+\cdots +F_t,G_1+\cdots +G_{t'})$ are multidistributive, then $(E,E'_1+E'_2,F_1+\cdots +F_t,G_1+\cdots +G_{t'})$ and $(E,E'_1\cap E'_2,F_1+\cdots +F_t,G_1+\cdots + G_{t'})$ are also multidistributive.
\end{coro}

\bigskip

From now on, we will consider an $(a,b)$-homogeneous algebra $A=T(V)/I$ with space of relations $R=R_a\oplus R_b$ such that $R_a$ and $R_b$ are exclusive. For $n\in \mathbb{N}$, denote by $\mathcal{T}_n$ the sublattice of $\mathcal{L}(V^{(n)})$ generated by:
\begin{equation}
\begin{cases}
V^{(i)}\otimes R_a\otimes V^{(n-i-a)} \text{ with $0\le i\le n-a$,} & \text{if $a\le n<b$},
\\
V^{(i)}\otimes R_a\otimes V^{(j)} \text{ and } V^{(h)}\otimes R_b \otimes V^{(l)}, \text{ with $i+j+a=h+l+b=n$},  & \text{if $n\ge b$}.
\end{cases}
\end{equation}


\subsection{\texorpdfstring{Opposite algebra and Koszul dual algebra}{sec:oppositealganddualalg}}

In this section we prove that the opposite $(a,b)$-homogeneous algebra $A^{\circ}$ of an algebra $A=T(V)/I$ is $(a,b)$-Koszul if and only if $A$ is $(a,b)$-Koszul. Moreover, we define the dual $(a,b)$-homogeneous algebra of $A$ and exhibit an $(a,b)$-Koszul algebra whose dual algebra is not $(a,b)$-Koszul.

As always, $R=R_a\oplus R_b$ is a space of relations for $I$, where $R_a$ and $R_b$ are exclusive. The $k$-linear endomorphism of $T(V)$ given by
\begin{itemize}
\item[$\bullet$] $\tau (1)=1$,
\item[$\bullet$] $\tau(w_1\otimes w_2\otimes \cdots \otimes w_n)=w_n\otimes \cdots \otimes w_2\otimes w_1$ for $w_1,\cdots ,w_n\in V$ and $n\ge 1$,
\end{itemize}
is an anti-isomorphism and it induces an algebra anti-isomorphism:
\[
\overline{\tau}: A\longrightarrow \frac{T(V)}{\langle\tau (R)\rangle}=A^{\circ}.
\]

It is clear that if $R_a$ and $R_b$ are exclusive, then $\tau (R)=\tau (R_a)\oplus \tau (R_b)$ and $\tau (R_a)$ and $\tau (R_b)$ are also exclusive. We say that $A^{\circ}$ is the \textbf{opposite $(a,b)$-homogeneous algebra} of $A$.

\begin{prop}
The algebra $A^{\circ}$ is $(a,b)$-Koszul if and only if $A$ is $(a,b)$-Koszul.
\end{prop}

\begin{prof}
The map $\overline{\tau}$ induces an isomorphism from $\Mod\text{-}A$ to $A^{\circ}\text{-}\Mod$ preserving the objects, where the left $A^{\circ}$-module structure of a right $A$-module $M$ is given by
\begin{align*}
A^{\circ} \otimes M & \rightarrow  M
\\
a\otimes m & \mapsto  am:=m\overline{\tau}^{-1}(a).
\end{align*}

Suppose that $A$ is $(a,b)$-Koszul. Applying $\overline{\tau}$ to the Koszul resolution of the right $A$-module $k$, we obtain a Koszul resolution of the left $A^{\circ}$-module $k$ since $\overline{\tau}$ is an homogeneous anti-isomorphism. In particular, $A$ is right $(a,b)$-Koszul if and only if $A^{\circ}$ is left $(a,b)$-Koszul. \qed
\end{prof}

\bigskip

Even if $A$ is not $(a,b)$-Koszul, we may still prove the following.

\begin{prop}
The \textit{(e.c.)} hold for $A^{\circ}$ if and only if $A$ satisfies its \textit{(e.c.)}.
\end{prop}

\begin{prof}
Note first that for $W_1,W_2\subseteq V^{(n)}$, 
\begin{itemize}
\item[\textit{(i)}] $\tau(W_1+W_2)=\tau (W_1)+\tau (W_2)$, \hskip 0.6cm \textit{(iii)} \hskip 0.1cm $\tau (V^{(n)})=V^{(n)}$,
\item[\textit{(ii)}] $\tau (W_1\cap W_2)=\tau (W_1)\cap \tau (W_2)$, \hskip 0.6cm \textit{(iv)} \hskip 0.1cm $\tau (W_1\otimes W_2)= \tau(W_2)\otimes \tau(W_1)$.
\end{itemize}

Suppose that the \textit{(e.c.)} are satisfied for $A$, then by Proposition \ref{prop:eciffdistributivity} and for $2\le m\le a-1$, $2\le l\le b-1$ the triples
\begin{align*}
(E,F,G) = & \hskip 0.2cm  (V^{(m)}\otimes R_a,R_a\otimes V^{(m)}, V\otimes R_a\otimes V^{(m-1)}+\cdots +V^{(m-1)}\otimes R_a\otimes V) \hskip 0.2cm \text{and}
\\
(\hat E,\hat F,\hat G) = & \hskip 0.2cm (V^{(l)}\otimes R_b,R_b\otimes V^{(l)}, V\otimes R_b\otimes V^{(l-1)}+\cdots +V^{(l-1)}\otimes R_b\otimes V)
\end{align*}

\noindent are distributive.

Let $F'=\tau (E)$, $E'=\tau (F)$ and $G'=\tau (G)$. If $v\in E'\cap (F'+G')$, then there exist $v_1\in F'$ and $v_2\in G'$ such that $v=v_1+v_2$. Thus, $\tau^{-1}(v)=\tau^{-1}(v_1)+\tau^{-1}(v_2)$, and then $\tau^{-1}(v_1)=\tau^{-1}(v)-\tau^{-1}(v_2)$. Since $(E,F,G)$ is distributive, there exist $w_1\in E\cap F$ and $w_2\in E\cap G$ such that $\tau^{-1}(v_1)=w_1+w_2$. Therefore, $v_1=\tau (w_1)+\tau (w_2)$ and $v=\tau (w_1)+\tau (w_2)+v_2$, where $\tau (w_1)\in E'\cap F'$, $\tau (w_2)+v_2\in G'$ and $\tau(w_2)+v_2=v-\tau (w_1)\in E'$. Hence, $v\in (E'\cap F')+(E'\cap G')$ and then $(E',F',G')$ is distributive. The other case is analogous.

By Lemma \ref{lem:equivalenceforincl} the inclusions $(V^{(m)}\otimes R_a)\cap (R_a\otimes V^{(m)}) \subseteq V^{(m-1)}\otimes R_a\otimes V$ and $(V^{(l)}\otimes R_b)\cap (R_b\otimes V^{(l)}) \subseteq V^{(l-1)}\otimes R_b\otimes V$ are equivalent to
\[
(V^{(m)}\otimes R_a)\cap (R_a\otimes V^{(m)})  = J_{a+m}^a \hskip 0.2cm \text{and} \hskip 0.2cm (V^{(l)}\otimes R_b)\cap (R_b\otimes V^{(l)}) = J_{b+l}^b.
\]

\noindent Applying again $\tau$ we verify that the \textit{(e.c.)} hold for $A^{\circ}$. The proof of the converse statement is similar. \qed
\end{prof}

\begin{rmk}
There are no analogous equivalences for the \textit{(e.v.c.)} or for the \textit{(e.c.c.)}. For instance, consider $A=\frac{k\langle x,y\rangle}{\langle x^3,xy^3\rangle}$. It is easy to see that the \textit{(e.v.c.)} are true for $A$ but not for $A^{\circ}$.

\medskip
For $A=\frac{k\langle x,y\rangle}{\langle xy^2,x^4+x^3y\rangle}$, the \textit{(e.c.c.)} hold, but they are not satisfied for $A^{\circ}$.
\end{rmk}

\medskip

We have already said that $(a,b)$-Koszul algebras are generalizations of $N$-Koszul algebras. It is well-known (\cite{BDW}) that the dual algebra $A^!$ of an $N$-Koszul algebra $A$ may not be Koszul for $N\ge 3$. Of course, the same happens for $(a,b)$-Koszul algebras.

\bigskip

Given $s\in \mathbb{N}$, let $R_s^{\perp}=\{ f\in (V^{(s)})^* / f(R_s)=0\}$. We define the \textbf{dual $(a,b)$-homogeneous algebra of $A$} as $A^{!}=\frac{T(V^*)}{I(R_a^{\perp}\oplus R_b^{\perp})}$.

Next we introduce some notation. For $s=a,b$
\[
J(R_s)_n = \begin{cases}
V^{(n)} & \text{if} \hskip 0.2cm 0\le n\le s-1,
\\
J_n^s & \text{if} \hskip 0.2cm n\ge s.
\end{cases}
\]

\begin{prop}
$(I(R_a^{\perp}\oplus R_b^{\perp}))_n =(J(R_a)^{\perp})_n+(J(R_b)^{\perp})_n$.
\end{prop}

\begin{prof}
It follows from the isomorphism $(V^{(i)}\otimes R\otimes V^{(j)})^{\perp} \simeq (V^*)^{(i)}\otimes R^{\perp}\otimes (V^*)^{(j)}$, where $R$ is any subspace of $V^{(n)}$. \qed
\end{prof}

\medskip

The map
\[
\begin{array}{cccc}
\rho : & \mathcal{L}(V^{(n)}) & \rightarrow & \mathcal{L}((V^{(n)})^*)
\\
& W & \mapsto & W^{\perp}.
\end{array}
\]
is bijective and $\rho (W+W')=\rho (W)\cap \rho (W')$ and $\rho (W\cap W') =\rho (W)+ \rho (W')$ for all $W,W'\in \mathcal{L}(V^{(n)})$.

Moreover, $\rho$ transforms $\mathcal{T}_n$ into the sublattice generated by $(V^*)^{(i)}\otimes R_a^{\perp}\otimes (V^*)^{(j)}$ and $(V^*)^{(h)}\otimes R_b^{\perp}\otimes (V^*)^{(l)}$ for $i+j+a=h+l+b=n$.

Suppose now that $R_a^{\perp}$ and $R_b^{\perp}$ are exclusive. Note that, applying $\rho$,
\begin{itemize}
\item[\textit{(i)}] The \textit{(e.c.)} hold for $A^{!}$ if and only if
\begin{align*}
V^{(a-2)}\otimes R_a\otimes V & \subseteq (V^{(a-1)}\otimes R_a) +(R_a\otimes V^{(a-1)}\cap \cdots \cap V^{(a-2)}\otimes R_a\otimes V) \text{ and}
\\
V^{(b-2)}\otimes R_b\otimes V & \subseteq (V^{(b-1)}\otimes R_b) +(R_b\otimes V^{(b-1)}\cap \cdots \cap V^{(b-2)}\otimes R_b\otimes V).
\end{align*}
\item[\textit{(ii)}] The \textit{(e.v.c.)} hold for $A^{!}$ if and only if
\begin{align*}
(V^{(b-1)}\otimes R_a\otimes V)+(V^{(b)}\otimes R_a)+(R_b\otimes V^{(a)}) & = V^{(a+b)} \text{ and}
\\
(V^{(a-1)}\otimes R_b\otimes V)+(V^{(a)}\otimes R_b)+(R_a\otimes V^{(b)}) & = V^{(a+b)}.
\end{align*}
\item[\textit{(iii)}] The \textit{(e.c.c.)} hold for $A^{!}$ if and only if
\begin{align*}
(V^{(b-1)}\otimes R_a)+(R_b\otimes V^{(a-1)}\cap \cdots \cap V^{(a-2)}\otimes R_b\otimes V) & = V^{(a+b-1)} \text{ and}
\\
(V^{(a-1)}\otimes R_b)+(R_a\otimes V^{(b-1)}\cap \cdots \cap V^{(b-2)}\otimes R_a\otimes V) & = V^{(a+b-1)}.
\end{align*}
\end{itemize}

\bigskip

\begin{example}
If $R_a$ and $R_b$ are exclusive, $R_a^{\perp}$ and $R_b^{\perp}$ are not necessarily exclusive. For instance, $R_3=\langle x^3\rangle$ and $R_4=\langle y^4\rangle$ are exclusive but $(x^*)^2y^*\in R_3^{\perp}$ and $(x^*)^3y^*\in R_4^{\perp}$, and so $(V^*\otimes R_3^{\perp})\cap R_4^{\perp}\ne 0$.
\end{example}


\section{\texorpdfstring{Examples}{sec:example}}\label{sec:example}

In this section we exhibit some examples of $(a,b)$-Koszul algebras.

Given $(a,b)$ with $4\le a<b$, we define $\tilde{A}_{a,b}=k\langle x,y\rangle/I$, with $I$ generated by a space of relations $R=R_a\oplus R_b$ where
\begin{itemize}
\item[$\bullet$] $R_a=\langle x^2w_1\cdots w_{a-4}y^2\rangle$ such that if $t=\left[ \frac{a+1}{2}\right]-1$, for $2\le j\le t$, $w_{a-2j}=x$ and $w_{a-2j-1}=y$, and ,
\item[$\bullet$] $R_b=\langle x^2y^{b-4}xy\rangle$.
\end{itemize}

\noindent Our aim is to decide whether $\tilde{A}_{a,b}$ is $(a,b)$-Koszul or not.

\begin{rmk}\label{rmk:avoidcase}
If $a=4$ and $b\ge 6$, then $R_a$ and $R_b$ are not exclusive. We restrict to $(a,b)=(4,5)$ if $a=4$.
\end{rmk}

To check the \textit{(e.c.)}, note that $\{ v_{i_1}\cdots v_{i_{a-1}}x^2w_1\cdots w_{a-4}y^2 \hskip 0.2cm / \hskip 0.2cm v_{i_j}=x,y\}$ is a basis of $V^{(a-1)}\otimes R_a$. It is clear that $(V^{(a-1)}\otimes R_a)\cap (R_a\otimes V^{(a-1)}+\cdots +V^{(a-2)}\otimes R_a\otimes V)$ vanishes. Moreover, $J_{a+1}^a=(V\otimes R_a)\cap (R_a\otimes V)=0$. Thus, $(V^{(a-1)}\otimes R_a)\cap (R_a\otimes V^{(a-1)}+\cdots +V^{(a-2)}\otimes R_a\otimes V) = V^{(a-2)}\otimes J_{a+1}^a$. The argument is similar for the other case of the \textit{(e.c.)}.

\medskip

We now verify the \textit{(e.v.c.)}. Observe that $(V^{(b-1)}\otimes R_a\otimes V)\cap (V^{(b)}\otimes R_a) \cap (R_b\otimes V^{(a)})\subseteq (V^{(b)}\otimes R_a)\cap (R_b\otimes V^{(a)})=R_b\otimes R_a=\langle x^2y^{b-4}xyx^2w_1\cdots w_{a-4}y^2\rangle$; however, $x^2y^{b-4}xyx^2w_1\cdots w_{a-4}y^2\not \in V^{(b-1)}\otimes R_a\otimes V$. Thus, the intersection is zero. Analogously, $(V^{(a-1)}\otimes R_b\otimes V)\cap (V^{(a)}\otimes R_b) \cap (R_a\otimes V^{(b)})=0$.

\medskip

For the \textit{(e.c.c.)}, given $v\in (V^{(b-1)}\otimes R_a)\cap (R_b\otimes V^{(a-1)}+\cdots +V^{(a-2)}\otimes R_b\otimes V)$, there exist $\lambda_i \in k$ and $v_{i_j}\in \{x,y\}$ for $1\le j\le a+b-1$ such that $v=\sum_{i=(i_1,\cdots ,i_{a+b-1})} \lambda_i v_{i_1}\cdots v_{i_{a+b-1}}$. It must be $v_{i_1}\cdots v_{i_{a+b-1}} = v_{i_1}\cdots v_{i_{b-1}}x^2w_1\cdots w_{a-4}y^2$. Also, $v$ belongs to the sum only if there exists $s$ such that $0\le s\le a-2$ and $v_{i_1}\cdots v_{i_{a+b-1}}= v_{i_1}\cdots v_{i_s}x^2y^{b-4}xyv_{i_{b+s+1}}\cdots v_{i_{a+b-1}}$. It is easy to see that this is a contradiction and hence the intersection vanishes.

Similarly, $(V^{(a-1)}\otimes R_b)\cap (R_a\otimes V^{(b-1)}+\cdots +V^{(b-2)}\otimes R_a\otimes V)=0$.

\medskip

Given $\mathcal{B}_n=\{ v_{i_1}\cdots v_{i_n} \hskip 0.2cm / \hskip 0.2cm v_{i_j}=x\hskip 0.2cm \text{or} \hskip 0.2cm v_{i_j}=y \hskip 0.2cm \text{for} \hskip 0.2cm 1\le j\le n\}$ which is a basis of $V^{(n)}$, it is clear that $\mathcal{B}_n$ distributes with $\mathcal{T}_n$.

\medskip

Next, considering a minimal resolution of $\tilde{A}_{a,b}$ we shall prove that $\ker \delta_2$ is $2$-pure in degrees $a+1$ and $b+1$. By previous arguments and Proposition \ref{prop:distbasis}, we know that the sublattice $\mathcal{T}_n$ is distributive and hence the triples defined in Proposition \ref{prop:kerdelta2pure} are also distributive. Moreover, following the notations of Proposition \ref{prop:kerdelta2pure}, we know that $(E'\oplus E'')\cap (F'+G'+F''+G'') = E'\cap (F'+G')+E'\cap (F''+G'')+E''\cap (F'+G')+E''\cap (F''+G'')$.

The set $\{ v=v_{i_1}\cdots v_{i_{n-a}}x^2w_1\cdots w_{a-4}y^2 \hskip 0.2cm / \hskip 0.2cm v_{i_j}=x\hskip 0.2cm \text{or} \hskip 0.2cm v_{i_j}=y\}$, appearing in the definition of $R_a$, is a basis of $E'\cap (F''+G'')$. 
\begin{itemize}
\item[$\bullet$] If $v\in F''$ then $v\in F'$. 
\item[$\bullet$] If $v\in G''=(V^{(n-2b+1)}\otimes R_b\otimes V^{(b-1)}+\cdots +V^{(n-a-b)}\otimes R_b\otimes V^{(a)})+(V^{(n-a-b+1)}\otimes R_b\otimes V^{(a-1)}+\cdots +V^{(n-b-1)}\otimes R_b\otimes V)$ and there exists $s$ such that $1\le s\le a-1$ and $v=v_{j_1}\cdots v_{j_{n-b-s}}x^2y^{b-4}xyv_{j_{n-s+1}}\cdots v_{j_n}$, then $x^2y^{b-4}xyv_{j_{n-s+1}}\cdots v_{j_n}=z_1x^2w_1\cdots w_{a-4}y^2z_2$ where $z_1$ and $z_2$ belong respectively to $\langle x,y\rangle^{(s)}$ and $\langle x,y\rangle^{(t)}$ for some $s,t\in \mathbb{N}_0$. This is possible only if $(a,b)=(4,6)$, but this case is excluded. Otherwise, $v\in F'$.
\end{itemize}

Therefore, $E'\cap (F''+G'')\subseteq E'\cap (F'+G')$.

We now want a basis of $E''\cap (F'+G')$. Given $v=v_{i_1}\cdots v_{i_{n-b}}x^2y^{b-4}xy$ in a basis, let us first suppose that $v\in F'$. Recall that $F'=(R_a\otimes V^{(n-a)}+\cdots +V^{(n-a-b)}\otimes R_a\otimes V^{(b)})+(V^{(n-a-b+1)}\otimes R_a\otimes V^{(b-1)}+\cdots +V^{(n-2a)}\otimes R_a\otimes V^{(a)})
+(R_b\otimes V^{(n-b)}+\cdots +V^{(n-2b)}\otimes R_b\otimes V^{(b)})+(V^{(n-2b+1)}\otimes R_b\otimes V^{(b-1)}+ \cdots +V^{(n-a-b)}\otimes R_b\otimes V^{(a)})= F''+(V^{(n-a-b+1)}\otimes R_a\otimes V^{(b-1)}+\cdots +V^{(n-2a)}\otimes R_a\otimes V^{(a)})+V^{(n-2b+1)}\otimes R_b\otimes V^{(b-1)}+\cdots +V^{(n-a-b)}\otimes R_b\otimes V^{(a)}$.

\noindent If $v\in D=(V^{(n-a-b+1)}\otimes R_a\otimes V^{(b-1)}+\cdots +V^{(n-2a)}\otimes R_a\otimes V^{(a)})$, there exists $s$ such that $a\le s\le b-1$ and $v=v_{j_1}\cdots v_{j_{n-a-s}}x^2w_1\cdots w_{a-4}y^2v_{j_{n-s+1}}\cdots v_{j_n}$. Note that $v_{j_{n-s-2}}=x$ and $v_{j_{n-s-1}}=v_{j_{n-s}}=y$ and using that $v\in E''$,
\begin{itemize}
\item[$\bullet$] there exists $v_h$ such that $h\le n-b-2$, so $n-s-2\le n-b-2$, but this is a contradiction.
\item[$\bullet$] there exists $v_h$ such that $h\le n-b+2$ and so $s=b-4$. As a consequence, $v_{i_1}\cdots v_{i_{n-b}}x^2y^{b-4}xy=v_{j_1}\cdots v_{j_{n-a-b+4}}x^2w_1\cdots w_{a-4}y^2v_{j_{n-b+5}}\cdots v_{j_n}$ and $v_{j_{n-b+1}}=x=w_{a-3}=y$ giving a contradiction if $a\ge 6$.

If $a=4$, then $v=v_{j_1}\cdots v_{j_{n-b}}x^2y^2v_{j_{n-b+5}}\cdots v_{j_n}$ and this happens only if $b\ge 6$. 

If $a=5$, then $v=v_{j_1}\cdots v_{j_{n-b-1}}x^3y^2v_{j_{n-b+5}}\cdots v_{j_n}$ which implies $b\ge 6$, $v_{i_{n-b}}=x$ and $v_{j_{n-b+5}}\cdots v_{j_n}=y^{b-6}xy$. Note that $v_{j_1}\cdots v_{j_{n-b-1}}x^3y^2v_{j_{n-b+5}}\cdots v_{j_n}\in V^{(n-b-1)}\otimes R_5\otimes V^{(b-4)}$ and that $V^{(n-b-1)}\otimes R_5\otimes V^{(b-4)}\cap (F''+G'')=0$.
\end{itemize}

For $n\ge a+b-1$, 
\begin{align*}
E''\cap G' & =E''\cap G'=(V^{(n-b)}\otimes R_b)\cap (V^{(n-2a+1)}\otimes R_a\otimes V^{(a-1)}+\cdots +V^{(n-a-1)}\otimes R_a\otimes V)
\\
& = V^{(n-a-b+1)}\otimes [(V^{(a-1)}\otimes R_b)\cap (V^{(b-a)}\otimes R_a\otimes V^{(a-1)}+\cdots +V^{(b-2)}\otimes R_a\otimes V)]
\\
& \subseteq V^{(n-a-b+1)}\otimes [(V^{(a-1)}\otimes R_b)\cap (R_a\otimes V^{(b-1)}+\cdots +V^{(b-2)}\otimes R_a\otimes V)]=0.
\end{align*}

\noindent If $2a\le n<a+b-1$ and $v\in E''\cap G'$ we argue as for $E''\cap D$. Hence, $E''\cap (F'+G')\subseteq E''\cap (F''+G'')$.

We conclude that $(E', E'',F'+G',F''+G'')$ is multidistributive and by Proposition \ref{prop:kerdelta2pure}, $\ker \delta_2$ is $2$-pure in degrees $a+1$ and $b+1$.

\medskip

Given $m\in \mathbb{N}$ such that $m\ge s+1$ for $s=a$ or $s=b$ respectively,
\begin{align*}
J_m^s & =(R_s\otimes V^{(m-s)})\cap (V\otimes R_s\otimes V^{(m-s-1)})\cap \cdots \cap (V^{(m-s)}\otimes R_s)
\\
& \subseteq V^{(m-s-1)}\otimes [(R_s\otimes V)\cap (V\otimes R_s)]=V^{(m-s-1)}\otimes J_{s+1}^s=0.
\end{align*}

\noindent Following Theorem \ref{thm:kerdeltai}, there are two cases:
\begin{itemize}
\item[$\bullet$] if $i>2$ is even, since $\frac{i}{2}s>s+1$, $E_s=V^{(n-\frac{i}{2}s)}\otimes J_{\frac{i}{2}s}^s=0$ for $s=a,b$;
\item[$\bullet$] if $i>3$ is odd,  since $\frac{i-1}{2}s+1\ge s+1$, $E_s=V^{(n-\frac{i-1}{2}s-1)}\otimes J_{\frac{i-1}{2}s+1}^s=0$ for $s=a,b$.
\end{itemize}

\noindent As a consequence, the distributivities and multidistributivities in Theorem \ref{thm:kerdeltai} are satisfied.

By Theorems \ref{thm:kerdeltai} and \ref{thm:Kozulequiv}, we obtain

\begin{prop}\label{prop:abKoszulexample}
The algebras $\tilde{A}_{a,b}$ are $(a,b)$-Koszul if and only if $(a,b)=(4,5)$ or $6\le a<b$.
\end{prop}

\begin{rmk}
Given $a,b\in \mathbb{Z}$ such that $2\le a<b$, the algebra $A=\frac{k\langle x,y\rangle}{\langle x^a,y^b\rangle}$ is not $(a,b)$-Koszul. For this, consider $(E_a, E_b, F_1+F_2+F_{13}+F_{14}, F_3+F_4+F_{15}+F_{16})$ as in Theorem \ref{thm:kerdeltai} and $i$ odd.

The basis $\mathcal{B}_n=\{ v_{i_1}\cdots v_{i_n} \hskip 0.2cm / \hskip 0.2cm v_{i_j}=x\hskip 0.2cm \text{or} \hskip 0.2cm v_{i_j}=y \hskip 0.2cm \text{for} \hskip 0.2cm 1\le j\le n\}$ of $V^{(n)}$ clearly distributes with $\mathcal{T}_n$ and by Proposition \ref{prop:distbasis}, $\mathcal{T}_n$ is distributive.

Let $i\ge 5$, $n\ge \frac{i+1}{2}b+1$ and $v=v_{j_1}\cdots v_{j_{n-\frac{i+1}{2}b+1}}y^{\frac{i+1}{2}b-1}$ where $h\ge 0$, $t\ge 1$ and 
\[
v_{j_t}= \begin{cases}
x & \text{if} \hskip 0.2cm t=n-\frac{i+1}{2}b-2h+1,
\\
y & \text{if} \hskip 0.2cm t=n-\frac{i+1}{2}b-2h.
\end{cases}
\]

\noindent Since $\frac{i-1}{2}b+1<\frac{i+1}{2}b-1$, $v\in E_b\cap F_{13}$. It is easy to see that $v$ belong neither to $E_a$ nor to $F_3+F_4+F_{15}+F_{16}$. Hence, the tuple is not multidistributive.
\end{rmk}

\bigskip

The following example of $(a,b)$-Koszul algebra is a quotient of a down-up algebra. We refer to \cite{B1} for more details on down-up algebras.

In \cite{B1}, the author considers a $3$-dimensional Lie algebra $\mathfrak{g}$ over $\mathbb{C}$, with basis $\{ x,y,[x,y]\}$ such that $[x,[x,y]]=\gamma x$ and $[[x,y],y]=\gamma y$, $\gamma \in \mathbb{C}$. The relations in the universal enveloping algebra $\mathcal{U}(\mathfrak{g})$ of $\mathfrak{g}$ become:
\begin{align*}
x^2y-2xyx+yx^2 & =\gamma x,
\\
xy^2-2yxy+y^2x & =\gamma y,
\end{align*}

\noindent and $\mathcal{U}(\mathfrak{g})\simeq A(2,-1,\gamma )$. For $\gamma =0$, they are homogeneous.

Let $\mathfrak{h}$ be the Heisenberg algebra generated by $x$ and $y$ and satisfying $[x,[x,y]]=0=[y,[x,y]]$. Then, $\mathcal{U}(\mathfrak{h})\simeq A(2,-1,0)$. 

Consider the algebra $A=\frac{k\langle x,y\rangle}{\langle R \rangle}$ where $R=R_3\oplus R_4$, $R_3=\langle x^2y-2xyx+yx^2,y^2x-2yxy+xy^2\rangle$ and $R_4=\langle x^4,y^4\rangle$. It is clear that $R_3$ and $R_4$ are exclusive.

The following statements are easy:

\begin{itemize}
\item[$\bullet$] $(V^{(2)}\otimes R_3)\cap (R_3\otimes V^{(2)}+V\otimes R_3 \otimes V)=\langle xyx^2y-2xyxyx+xy^2x^2+x^2y^2x-2x^2yxy+x^3y^2,y^2x^2y-2y^2xyx+y^3x^2+yxy^2x-2yxyxy+yx^2y^2\rangle =V\otimes J_4^3$. Also, $(V^{(3)}\otimes R_4)\cap (R_4\otimes V^{(3)}+V\otimes R_4\otimes V^{(2)}+V^{(2)}\otimes R_4\otimes V)=\langle x^7,xyx^5,yx^6,y^2x^5,x^2y^5,xy^6,yxy^5,y^7\rangle$ $=V^{(2)}\otimes J_5^4$. Thus the \textit{(e.c.)} are satistied.
\item[$\bullet$] $(V^{(3)}\otimes R_3\otimes V)\cap (V^{(4)}\otimes R_3)\cap (R_4\otimes V^{(3)})=0=(V^{(2)}\otimes R_4\otimes V)\cap (V^{(3)}\otimes R_4)\cap (R_3\otimes V^{(4)})$. Hence, the \textit{(e.v.c.)} hold.
\item[$\bullet$] $(V^{(3)}\otimes R_3)\cap (R_4\otimes V^{(2)}+V\otimes R_4\otimes V)=0=(V^{(2)}\otimes R_4)\cap (R_3\otimes V^{(3)}+V\otimes R_3\otimes V^{(2)}+V^{(2)}\otimes R_3\otimes V)$. Therefore, the \textit{(e.c.c.)} are satisfied.
\end{itemize} 

\begin{lema}\label{lem:Jn4}
For all $n\ge 4$, $J_n^4=\langle x^n,y^n\rangle$.
\end{lema}

\begin{prof}
If $n=4$, $J_4^4=R_4=\langle x^4,y^4\rangle$. Suppose now that the statement is true for some $n\ge 4$. Then,
\begin{align*}
J_{n+1}^4 & =(R_4\otimes V^{(n-3)})\cap (V\otimes R_4\otimes V^{(n-4)})\cap \cdots \cap (V^{(n-3)}\otimes R_4)
\\
& =(R_4\otimes V^{(n-3)})\cap V\otimes [(R_4\otimes V^{(n-4)})\cap \cdots \cap (V^{(n-4)}\otimes R_4)] 
\\
& =(R_4\otimes V^{(n-3)})\cap (V\otimes J_n^4)
\\
& =\langle x^4v_1\cdots v_{n-3},y^4w_1\cdots w_{n-3} \hskip 0.2cm / \hskip 0.2cm v_i,w_j\in \{ x,y\}\rangle \cap \langle x^{n+1},yx^n,xy^n,y^{n+1}\rangle.
\end{align*}
Thus, it is clear that $J_{n+1}^4$ is generated by $x^{n+1}$ and $y^{n+1}$. \qed
\end{prof}

\begin{rmk}\label{rmk:Jn3}
Notice that $J_4^3=\langle x^2y^2-2xyxy+yx^2y+y^2x^2-2yxyx+xy^2x\rangle$ and that $J_5^3=(R_3\otimes V^{(2)})\cap (V\otimes R_3\otimes V)\cap (V^{(2)}\otimes R_3)=(R_3\otimes V^{(2)})\cap (V\otimes J_4^3)=0$. For $n\ge 6$, the arguments of Lemma \ref{lem:Jn4} prove that $J_n^3=(R_3\otimes V^{(n-3)})\cap (V\otimes J_{n-1}^3)$. Since $J_5^3=0$, then $J_n^3=0$ for all $n\ge 5$.
\end{rmk}

\bigskip

The Koszul complex of $A$ is
\[
\cdots \longrightarrow A\otimes J_{12}^4 \overset{\delta_6}{\longrightarrow} A\otimes J_9^4 \overset{\delta_5}{\longrightarrow} A\otimes J_8^4 \overset{\delta_4}{\longrightarrow} A\otimes J_4^3\oplus A\otimes J_5^4 \overset{\delta_3}{\longrightarrow} A\otimes R \overset{\delta_2}{\longrightarrow} A\otimes V \overset{\delta_1}{\longrightarrow} A \longrightarrow 0.
\]

Next we study the exactness of this complex.

\begin{itemize}
\item[$\bullet$] Let $\mathcal{T}_n$ be the subspace generated by
\[
\begin{cases}
V^{(i)}\otimes R_3\otimes V^{(n-i-3)} & \text{if $n=3$},
\\
V^{(i)}\otimes R_3\otimes V^{(j)} \text{ and } V^{(h)}\otimes R_4 \otimes V^{(l)}, \text{ with $i+j+3=h+l+4=n$},  & \text{if $n\ge 4$}.
\end{cases}
\]

It is clear that the basis $\mathcal{B}_n=\{ v_{i_1}\cdots v_{i_n} \hskip 0.2cm / \hskip 0.2cm v_{i_j}=x \text{ or } v_{i_j}=y \text{ for } 1\le j\le n\}$ of $V^{(n)}$ distributes with $\mathcal{T}_n$. By Proposition \ref{prop:distbasis} $\mathcal{T}_n$ is distributive.

With the notations of Proposition \ref{prop:kerdelta2pure}, $(E,F,G)$ is distributive. Using the \textit{(e.c.c.)} it is easy to see that $E'\cap (F''+G'')\subseteq E'\cap (F'+G')+(V^{(n-3)}\otimes R_3)\cap (V^{(n-6)}\otimes R_4\otimes V^{(2)}+V^{(n-5)}\otimes R_4\otimes V)=E'\cap (F'+G')$ and that $E''\cap (F'+G')\subseteq E''\cap (F''+G'')+(V^{(n-4)}\otimes R_4)\cap (V^{(n-6)}\otimes R_3\otimes V^{(3)}+V^{(n-5)}\otimes R_3\otimes V^{(2)}+V^{(n-4)}\otimes R_3\otimes V)=E''\cap (F''+G'')$. Since $(E'\oplus E'')\cap (F'+G'+F''+G'')=E'\cap (F'+G')+E'\cap (F''+G'')+E''\cap (F'+G')+E''\cap (F''+G'')$, $(E', E'',F'+G',F''+G'')$ is multidistributive.

By Proposition \ref{prop:kerdelta2pure}, $\ker \delta_2$ is $2$-pure in degrees $4$ and $5$.
\item[$\bullet$] Since $\ker \delta_3=A\overline{x^3}\otimes x^5+A\overline{y^3}\otimes y^5$ and $\im \delta_4$ is generated by $\delta_4(\overline{1}\otimes x^8)=\overline{x^3}\otimes x^5$ and $\delta_4(\overline{1}\otimes y^8)=\overline{y^3}\otimes y^5$, $\im \delta_4=\ker \delta_3$.
\item[$\bullet$] For $i\ge 4$ even,
\[
\begin{array}{cccc}
\delta_i : & A\otimes J_{2i}^4 & \longrightarrow & A\otimes J_{2i-3}^4
\\
 & \overline{\alpha}\otimes x^{2i} & \mapsto & \overline{\alpha x^3}\otimes x^{2i-3},
\\
 & \overline{\beta}\otimes y^{2i} & \mapsto & \overline{\beta y^3}\otimes y^{2i-3},
\end{array}
\]
and $\ker \delta_i=(A\overline{x}+A\overline{y})\otimes J_{2i}^4$.

Since
\[
\begin{array}{cccc}
\delta_{i+1} : & A\otimes J_{2i+1}^4 & \longrightarrow & A\otimes J_{2i}^4
\\
 & \overline{\alpha}\otimes x^{2i+1} & \mapsto & \overline{\alpha x}\otimes x^{2i},
\\
 & \overline{\beta}\otimes y^{2i+1} & \mapsto & \overline{\beta y}\otimes y^{2i},
\end{array}
\]
it is clear that $\im \delta_{i+1}=(A\overline{x}+A\overline{y})\otimes J_{2i}^4$ and then $\im \delta_{i+1}=\ker \delta_i$.
\item[$\bullet$] For $i\ge 5$ odd,
\[
\begin{array}{cccc}
\delta_i : & A\otimes J_{2i-1}^4 & \longrightarrow & A\otimes J_{2i-2}^4
\\
 & \overline{\alpha}\otimes x^{2i-1} & \mapsto & \overline{\alpha x}\otimes x^{2i-2},
\\
 & \overline{\beta}\otimes y^{2i-1} & \mapsto & \overline{\beta y}\otimes y^{2i-2},
\end{array}
\]
and $\ker \delta_i=(A\overline{x^3}+A\overline{y^3})\otimes J_{2i-1}^4$.

Moreover,
\[
\begin{array}{cccc}
\delta_{i+1} : & A\otimes J_{2i+2}^4 & \longrightarrow & A\otimes J_{2i-1}^4
\\
 & \overline{\alpha}\otimes x^{2i+2} & \mapsto & \overline{\alpha x^3}\otimes x^{2i-1},
\\
 & \overline{\beta}\otimes y^{2i+2} & \mapsto & \overline{\beta y^3}\otimes y^{2i-1},
\end{array}
\]
and then $\im \delta_{i+1}=(A\overline{x^3}+A\overline{y^3})\otimes J_{2i-1}^4$. Thus, $\im \delta_{i+1}=\ker \delta_i$.
\end{itemize}

As a consequence, the Koszul complex of $A$ is exact in positive degrees and by Proposition \ref{prop:Koszulcomplexequiv}, $A$ is $(3,4)$-Koszul. Note that $A\simeq \frac{\mathcal{U}(\mathfrak{h})}{\langle x^4,y^4\rangle}$.


\section{\texorpdfstring{Hochschild homology}{sec:hochschildhomology}}

In this section, $A$ denotes an $(a,b)$-homogeneous algebra. When $A$ is $(a,b)$-Koszul, we obtain a bimodule resolution of $A$ which is useful to compute its Hochschild homology. The abelian category of $\mathbb{Z}$-graded left bounded $A$-bimodules with bimodule morphisms preserving the grading will be denoted by $\mathcal{C}$. Let $A^{e}=A\otimes_k A^{op}$ be the enveloping algebra of $A$. It is well-known that $\mathcal{C}$ is naturally isomorphic to the category of $\mathbb{Z}$-graded left bounded left $A^{e}$-modules, so we can use the results and notations of \S \ref{sec:preliminaries} for the graded $k$-algebra $A^e$. In particular, the following statements hold (see \cite{B2}):
\begin{itemize}
\item[\textit{(i)}] An $A$-bimodule $M$ is projective in $\mathcal{C}$ if and only if it is graded-free.
\item[\textit{(ii)}] Every $A$-bimodule $M$ admits a projective cover, unique up to isomorphism.
\item[\textit{(iii)}] Every $A$-bimodule $M$ admits a minimal projective resolution, unique up to isomorphism.
\item[\textit{(iv)}] Every projective resolution of an $A$-bimodule $M$ contains a minimal projective resolution as a direct summand.
\end{itemize}

We recall the following result from \cite{B4}.

\begin{lema}\cite[Lemme 1.6]{B4}\label{lem:tensorexactness}
Suppose that $M,M',M''\in \mathcal{C}$ are graded-free. If the complex of graded left $A$-modules $M''\otimes_A k \overset{g\otimes 1_k}{\longrightarrow} M'\otimes_A k \overset{f\otimes 1_k}{\longrightarrow} M\otimes_A k$ is exact, then $M'' \overset{g}{\longrightarrow} M' \overset{f}{\longrightarrow} M$ is exact in $\mathcal{C}$.
\end{lema}

\bigskip

The following ideas are based on \cite{B3}. For $s=a,b$ and $i\ge 0$ we consider the left $A$-module $(\overline{K}_{L,s})_i=A\otimes J_{n_s(i)}^s$ and the $A$-linear morphism $(\delta_{L,s})_i: (\overline{K}_{L,s})_i\rightarrow (\overline{K}_{L,s})_{i-1}$ defined extending the natural inclusion $J_{n_s(i)}^s \hookrightarrow  A\otimes J_{n_s(i-1)}^s$. It is clear that $(\delta_{L,s})^s=0$ and then $(\overline{K}_{L,s},\delta_{L,s})$ is an $s$-complex. Analogously, $(\overline{K}_{R,s},\delta_{R,s})$ is an $s$-complex where $(\delta_{R,s})_i: (\overline{K}_{R,s})_i\rightarrow (\overline{K}_{R,s})_{i-1}$ is defined by the restriction of the $k$-linear map $V^{(n_s(i))}\rightarrow A\otimes V^{(n_s(i-1))}$, $ v_{j_1}\cdots v_{j_{n_s(i)}}\mapsto \overline{v_{j_{n_s(i-1)-1}}\cdots v_{j_{n_s(i)}}}\otimes v_{j_1}\cdots v_{j_{n_s(i-1)}}$. Thus, $\overline{K}_{L-R,s}=\overline{K}_{L,s}\otimes A=A\otimes \overline{K}_{R,s}$ is an $s$-complex of bimodules with differential $\delta'_{L,s}=\delta_{L,s}\otimes 1_A$ (respectively, $\delta '_{R,s}=1_A \otimes \delta_{R,s}$). To show that $\delta '_{L,s}$ and $\delta '_{R,s}$ commute, we look at the composition
\[
(\overline{K}_{L,s})_i\otimes_k A \overset{(\delta_{L,s})_i\otimes 1_A}{\longrightarrow} (\overline{K}_{L,s})_{i-1}\otimes_k A = A\otimes_k (\overline{K}_{R,s})_{i-1} \overset{1_A\otimes (\delta_{R,s})_{i-1}}{\longrightarrow} A\otimes_k (\overline{K}_{R,s})_{i-2}.
\]

\noindent If $i$ is even,
\begin{align*}
& \overline{\alpha}\otimes v_{j_1}\cdots v_{j_{\frac{i}{2}s}}\otimes \overline{\beta} \mapsto (\overline{\alpha v_{j_1}\cdots v_{j_{s-1}}}\otimes v_{j_s}\cdots v_{j_{\frac{i}{2}s}})\otimes \overline{\beta} = \overline{\alpha v_{j_1}\cdots v_{j_{s-1}}}\otimes (v_{j_s}\cdots v_{j_{\frac{i}{2}s}}\otimes \overline{\beta})
\\
& \mapsto \overline{\alpha v_{j_1}\cdots v_{j_{s-1}}}\otimes v_{j_s}\cdots v_{j_{\frac{i-2}{2}s+1}}\otimes \overline{v_{j_{\frac{i-2}{2}s+2}}\cdots v_{j_{\frac{i}{2}s}}\otimes \beta}.
\end{align*}

On the other hand, $(\delta_{L,s})_{i-1}\otimes 1_A\circ 1_A\otimes (\delta_{R,s})_i$ is in this case:
\begin{align*}
& \overline{\alpha}\otimes v_{j_1}\cdots v_{j_{\frac{i}{2}s}}\otimes \overline{\beta} \mapsto \overline{\alpha}\otimes (v_{j_1}\cdots v_{j_{\frac{i-2}{2}s+1}}\otimes \overline{v_{j_{\frac{i-2}{2}s+2}}\cdots v_{j_{\frac{i}{2}s}}\beta}) \mapsto
\\
& (\overline{\alpha}\otimes v_{j_1}\cdots v_{j_{\frac{i-2}{2}s+1}})\otimes \overline{v_{j_{\frac{i-2}{2}s+2}}\cdots v_{j_{\frac{i}{2}s}}\beta} \mapsto \overline{\alpha v_{j_1}\cdots v_{j_{s-1}}}\otimes v_{j_s}\cdots v_{j_{\frac{i-2}{2}s+1}}\otimes \overline{v_{j_{\frac{i-2}{2}s+2}}\cdots v_{j_{\frac{i}{2}s}}\beta}.
\end{align*}

\noindent The case $i$ odd is similar.

Setting $(K_{L-R})_i=(\overline{K}_{L,a})_i\otimes A \oplus (\overline{K}_{L,b})_i\otimes A$ and $d'=d'_a\oplus d'_b$ where 
\[
(d'_s)_i=\begin{cases}
\delta'_{L,s}-\delta'_{R,s} & \text{if $i$ is odd},
\\ 
(\delta'_{L,s})^{s-1}+(\delta'_{L,s})^{s-2}\delta'_{R,s}+\cdots +\delta'_{L,s}(\delta'_{R,s})^{s-2}+(\delta'_{R,s})^{s-1} & \text{if $i$ is even},
\end{cases}
\]

\noindent it is clear that $(d'_s)^2=0$ and thus $(K_{L-R},d')$ is a complex of pure projective modules. It is called the \textbf{Koszul complex of the $A$-bimodule $A$}.

\begin{thm}
Let $A=T(V)/I$ be an $(a,b)$-homogeneous algebra such that $I$ is generated by $R=R_a\oplus R_b$ with $R_a$ and $R_b$ exclusive. The augmented Koszul complex 
\begin{equation}\label{eq:bimodkoszulcomplex}
\cdots \longrightarrow (K_{L-R})_2 \overset{d'_2}{\longrightarrow} (K_{L-R})_1 \overset{d'_1}{\longrightarrow} (K_{L-R})_0 \overset{\mu}{\longrightarrow} A \longrightarrow 0
\end{equation}

\noindent is exact if and only if $A$ is $(a,b)$-Koszul.
\end{thm}

\begin{prof}
Applying the functor $\place\otimes_A k$ to \eqref{eq:bimodkoszulcomplex}, we obtain the complex $K_L$ with augmentation $\varepsilon$, which is exact when $A$ is $(a,b)$-Koszul. Since the $A$-bimodules $(K_{L-R})_i$ are graded-free and left bounded for all $i\in \mathbb{Z}$, Lemma \ref{lem:tensorexactness} implies that the complex \eqref{eq:bimodkoszulcomplex} is exact.

Let us denote by $\mathcal{C}_R$ the abelian category of $\mathbb{Z}$-graded left bounded right $A$-modules. If \eqref{eq:bimodkoszulcomplex} is exact then it is a projective resolution of $A$ in $\mathcal{C}_R$.

Since $A$ is clearly projective in $\mathcal{C}_R$, the complex \eqref{eq:bimodkoszulcomplex} is homotopically trivial. Finally, applying the functor $\place \otimes_A k$ to \eqref{eq:bimodkoszulcomplex}, the left Koszul complex is exact in positive degrees and by Proposition \ref{prop:Koszulcomplexequiv}, $A$ is $(a,b)$-Koszul. \qed
\end{prof}

\begin{rmk}
If $A$ is $(a,b)$-Koszul, using Proposition \ref{prop:essentialineachdegree}, Lemmas 1.3 and 1.4 in \cite{B4} and Proposition 7 in \cite{C}, it follows that the complex \eqref{eq:bimodkoszulcomplex} is a minimal projective resolution of $A$ in $\mathcal{C}$.
\end{rmk}

\bigskip

From now on $A$ will be an $(a,b)$-Koszul algebra.

Recall that since $A$ is $k$-flat, the Hochschild homology $HH_*(A)$ is isomorphic to $\tor_{*}^{A^{e}} (A,A)$, so $HH_*(A) \simeq H_* (A\otimes _{A^e} K_{L-R}, 1_A\otimes d')$.

Fix $i\ge 0$. We identify the complex $(A\otimes_{A^e} K_{L-R},1_A\otimes d')$ with $(\overline{K}_{L,a}+\overline{K}_{L,b},\overline{d})$ using the $k$-linear isomorphisms,
\[
\begin{array}{cccc}
f_i^s: & A\otimes_{A^e} (A\otimes_k J_{n_s(i)}^s \otimes_k A) & \rightarrow & A\otimes_k J_{n_s(i)}^s
\\
& \overline{\alpha} \otimes (\overline{\beta}\otimes m\otimes \overline{\gamma}) & \mapsto & \overline{\gamma \alpha \beta}\otimes m
\end{array} \hskip 0.2cm \text{and}
\]
\[
\begin{array}{cccc}
g_i^s: & A\otimes_k J_{n_s(i)}^s & \rightarrow & A\otimes_{A^e} (A\otimes_k J_{n_s(i)}^s\otimes_k A)
\\
& \overline{\alpha}\otimes m & \mapsto & \overline{\alpha}\otimes (1\otimes m\otimes 1),
\end{array}
\]

\noindent where the map $\overline{d}$ is defined as follows:
\begin{itemize}
\item[$\bullet$] if $v\in V$ and $\overline{\alpha}\in A$, then $\overline{d_1} (\overline{\alpha} \otimes v)=f_0^a((1_A\otimes d'_1)(\overline{\alpha}\otimes 1\otimes v\otimes 1))=f_0^a(\overline{\alpha}\otimes (\overline{v}\otimes 1-1\otimes \overline{v})) = \overline{\alpha v}-\overline{v\alpha}$;
\item[$\bullet$] if $i\ge 3$ is odd, $\overline{\alpha}\in A$, $v,v'\in V$ and $w\in V^{(\frac{i-1}{2}s-1)}$ for $s=a,b$, then
\[
\overline{d}_i(\overline{\alpha}\otimes vwv')=\overline{\alpha v}\otimes wv'-\overline{v'\alpha}\otimes vw;
\]
\item[$\bullet$] if $i$ is even and $v_{j_1},\cdots ,v_{j_{\frac{i}{2}s}}\in V$ for $s=a,b$, then
\begin{align*}
& \overline{d}_i(\overline{\alpha}\otimes v_{j_1}\cdots v_{j_{\frac{i}{2}s}})=\overline{\alpha v_{j_1}\cdots v_{j_{s-1}}}\otimes v_{j_s}\cdots v_{j_{\frac{i}{2}s}}+
\\
& \overline{v_{j_{\frac{i}{2}s}}\alpha v_{j_1}\cdots v_{j_{s-2}}}\otimes v_{j_{s-1}}\cdots v_{j_{\frac{i}{2}s-1}}+\cdots +\overline{v_{j_{\frac{i-2}{2}s+2}}\cdots v_{j_{\frac{i}{2}s}}\alpha}\otimes v_{j_1}\cdots v_{j_{\frac{i-2}{2}s+1}}.
\end{align*}
\end{itemize}

Therefore, $HH_*(A)$ is the homology of $(K_L,\overline{d})$.


\subsection{\texorpdfstring{Hochschild homology of $\tilde{A}_{a,b}$}{sec:hochschildhomology}}

In this section we will compute the Hochschild homology of the $(a,b)$-Koszul algebras $\tilde{A}_{a,b}$ defined in \S \ref{sec:example}, with $(a,b)=(4,5)$ or $6\le a<b$. Using the Koszul resolution of the $\tilde{A}_{a,b}$-bimodule $\tilde{A}_{a,b}$ we obtain that $HH_* (\tilde{A}_{a,b})$ is the homology of the complex
\[
0 \longrightarrow \tilde{A}_{a,b}\otimes (R_a\oplus R_b) \overset{\overline{d_2}}{\longrightarrow} \tilde{A}_{a,b}\otimes V \overset{\overline{d_1}}{\longrightarrow} \tilde{A}_{a,b} \longrightarrow 0,
\]
so $HH_n (\tilde{A}_{a,b})$ is zero for $n\ge 3$. 


\subsubsection{\texorpdfstring{Hochschild homology in degree zero}{sec:hh0}}

We know that $HH_0 (\tilde{A}_{a,b}) = \frac{\tilde{A}_{a,b}}{\im \overline{d_1}}$. The $k$-vector space $\im \overline{d_1}$ is generated by $\{ \overline{\alpha x-x\alpha}, \overline{\alpha y-y\alpha} / \overline{\alpha}\in \tilde{A}_{a,b}\}$. It is clear that $(HH_0(\tilde{A}_{a,b}))_0 =k$ and $(HH_0(\tilde{A}_{a,b}))_1=\langle x,y \rangle$.

The homology class of an element will be denoted by $[\place ]$.

\begin{rmk}\label{rmk:permutation}
As usual, the class of an element $[\overline{w_1\cdots w_n}]\in (HH_0(\tilde{A}_{a,b}))_n$ consists of the cyclic permutations of $\overline{w_1\cdots w_n}$.
\end{rmk}

Given a necklace with $n$ beads, we assign to each bead either $x$ or $y$. Two assignments are equal if there is a rotation sending one into the other. By Remark \ref{rmk:permutation} it is clear that the number of nonzero different classes of $(HH_0(\tilde{A}_{a,b}))_n$ is bounded by the number of these assignments. Note that some classes could vanish as a consequence of the relations.

We recall the following definition:

\begin{defi}\cite{A} 
The \textbf{Euler map} $\varphi :\mathbb{N}\rightarrow \mathbb{N}$ is such that $\varphi (n)$ is the number of positive integers less than or equal to $n$, coprime with $n$. In other words, $\varphi (n) = n \prod_{p|n} \left( 1-\frac{1}{p}\right)$ where $p$ is a positive prime.
\end{defi}

The number of the assignments in our case is given by $\rho (n)= \frac{1}{n}\sum_{m|n} \varphi (m) 2^{\frac{n}{m}}$.

\bigskip

We shall compute the dimensions of $(HH_0(\tilde{A}_{a,b}))_n$ for $n=2,3,4$ and we shall give an algorithm that computes the dimensions of the $n$-graded components for $n\ge 5$.

Note that $\overline{x^n}$ and $\overline{y^n}$ belong to two different nonzero unitary classes. For $0\le s\le n$, we denote  $\mathfrak{C}^n_{n-s}=[\overline{x^{n-s}y^s}]$. They are $n+1$ different classes unless they vanish in $\tilde{A}_{a,b}$.

Consider the following cases:
\begin{itemize}
\item[\textit{(i)}] For $n=2$, $\mathfrak{C}^2_2 =[\overline{x^2}]$, $\mathfrak{C}^2_1 =[\overline{xy}]$, $\mathfrak{C}^2_0 = [\overline{y^2}]$. There are $\rho (2)=3$ classes and none of them vanishes in $\tilde{A}_{a,b}$.
\item[\textit{(ii)}] For $n=3$ there are $\rho (3)=4$ nonvanishing classes.
\item[\textit{(iii)}] For $n=4$, $[\overline{xyxy}]$, $\mathfrak{C}^4_4 = [\overline{x^4}]$, $\mathfrak{C}^4_3 = [\overline{x^3y}]$, $\mathfrak{C}^4_2 = [\overline{x^2y^2}]$, $\mathfrak{C}^4_1 = [\overline{xy^3}]$, $\mathfrak{C}^4_0 = [\overline{y^4}]$ are $\rho (4)=6$ different classes and none of them vanishes except for $\tilde{A}_{4,5}$ where $\mathfrak{C}^4_2=0$.
\end{itemize}

\begin{rmk}
Let $n\ge 5$ and $w=x^{n-2}w_1w_2$. Suppose that $w_1=x$ or $w_2=x$, then $\overline{w}\in \mathfrak{C}^n_n$ if $w_1=w_2$ and $w\in \mathfrak{C}^n_{n-1}$ if $w_1\ne w_2$. If $w_1=w_2=y$, then $w\in \mathfrak{C}^n_{n-2}$.
\end{rmk}

\begin{rmk}\label{rmk:basispart}
The classes $[\overline{x^n}], [\overline{y^n}], [\overline{x^{n-1}y}], [\overline{y^{n-1}x}]\in (HH_0(\tilde{A}_{a,b}))_n$ are different and never vanish for $n\ge 3$. If $n$ is even, then $[\overline{xyx\cdots y}]\in (HH_0(\tilde{A}_{a,b}))_n$ is a new nonzero class.
\end{rmk}

The set of elements obtained in Remark \ref{rmk:basispart} may be extended to a basis of $(HH_0(\tilde{A}_{a,b}))_n$. In \S \ref{sec:appendix} we exhibit an algorithm which provides the remaining elements.

We conclude that
\[
\dim_k((HH_0(\tilde{A}_{a,b}))_n)=\begin{cases}
0 & \text{for $n<0$,}
\\
1 & \text{for $n=0$,}
\\
2 & \text{for $n=1$,}
\\
3 & \text{for $n=2$,}
\\
predim0(n,a,b)+4 & \text{for $n\ge 3$ odd,}
\\
predim0(n,a,b)+5 & \text{for $n\ge 4$ even.}
\end{cases}
\]

The integer $predim0(n,a,b)$ can be computed with the program in \S \ref{sec:appendix}.

We illustrate what we obtain from the algorithm computing bases of $(HH_0(\tilde{A}_{4,5}))_n$ for $n=5,6,7,8$:
\begin{itemize}
\item[$-$] $[\overline{x^5}]$, $[\overline{y^5}]$, $[\overline{x^4y}]$, $[\overline{y^4x}]$, $[\overline{xyxy^2}]$.
\item[$-$] $[\overline{x^6}]$, $[\overline{y^6}]$, $[\overline{x^5y}]$, $[\overline{y^5x}]$, $[\overline{xyxyxy}]$, $[\overline{x^2yx^2y}]$, $[\overline{xy^2xy^2}]$, $[\overline{xyxy^3}]$.
\item[$-$] $[\overline{x^7}]$, $[\overline{y^7}]$, $[\overline{x^6y}]$, $[\overline{y^6x}]$, $[\overline{xyxyxy^2}]$, $[\overline{x^3yx^2y}]$, $[\overline{xy^2xy^3}]$, $[\overline{xyxy^4}]$.
\item[$-$] $[\overline{x^8}]$, $[\overline{y^8}]$, $[\overline{x^7y}]$, $[\overline{y^7x}]$, $[\overline{xyxyxyxy}]$, $[\overline{xyxy^2xy^2}]$, $[\overline{x^3yx^3y}]$, $[\overline{xyxyxy^3}]$, $[\overline{xy^3xy^3}]$, $[\overline{x^4yx^2y}]$, $[\overline{xy^2xy^4}]$, $[\overline{xyxy^5}]$.
\end{itemize}


\subsubsection{\texorpdfstring{Computation of $HH_1(\tilde{A}_{a,b})$}{sec:hh1}}

Let $\overline{\alpha},\overline{\beta}\in \tilde{A}_{a,b}$. Note that $\im \overline{d_2}$ is generated by:
\begin{align*}
\overline{d_2}(\overline{\alpha}\otimes x^2w_1\cdots w_{a-4}y^2) &= \overline{\alpha x^2w_1\cdots w_{a-4}y}\otimes y+\overline{y\alpha x^2w_1\cdots w_{a-4}}\otimes y+\cdots +\overline{xw_1\cdots w_{a-4}y^2\alpha}\otimes x,
\\
\overline{d_2}(\overline{\beta}\otimes x^2y^{b-4}xy) & = \overline{\beta x^2y^{b-4}x}\otimes y+\overline{y\beta x^2y^{b-4}}\otimes x+\cdots +\overline{xy^{b-4}xy\beta}\otimes x.
\end{align*}

Moreover, $\ker \overline{d_1}$ consists of elements $\sum_j \overline{\alpha_j}\otimes x+\sum_h\overline{\beta_h}\otimes y$ where $\overline{\alpha_j},\overline{\beta_h}\in \tilde{A}_{a,b}$ are such that $\sum_j (\overline{\alpha_j x}-\overline{x\alpha_j})+\sum_h (\overline{\beta_h y}-\overline{y\beta_h})=0$. It is clear that $(\ker \overline{d_1})_0=0$ and that $(\ker \overline{d_1})_1=\langle \overline{1}\otimes x, \overline{1}\otimes y\rangle$.

Given $w_1\cdots w_n\in V^{(n)}$, we denote 
\[
\Gamma (w_1\cdots w_n)=\overline{w_1\cdots w_{n-1}}\otimes w_n+\overline{w_nw_1\cdots w_{n-2}}\otimes w_{n-1}+\cdots +\overline{w_2\cdots w_n}\otimes w_1 \in A\otimes V.
\]

For $n\ge 2$ it is easy to see that $\overline{x^{n-1}}\otimes x$, $\overline{y^{n-1}}\otimes y$, $\Gamma (x^{n-1}y)$ and $\Gamma (y^{n-1}x)$ belong to $(\ker \overline{d_1})_n$. Also, they do not vanish and clearly they do not belong to $(\im \overline{d_2})_n$ for any pair $(a,b)$. Notice that if $n$ is even, then $\overline{xyx\cdots x}\otimes y+\overline{yxy\cdots y}\otimes x\in (\ker \overline{d_1})_n\setminus (\im \overline{d_2})_n$.

In \S \ref{sec:appendix}, we give an algorithm to obtain a basis of $(HH_1(\tilde{A}_{a,b}))_n$ and the function $predim1(n)$.

For $t=2,\cdots , \big[ \frac{n}{2}\big]$ and $v_{j_l}\in \{ x,y\}$, $\Gamma (x^{n-t}yv_{j_1}\cdots v_{j_{t-2}}y)\in (\ker \overline{d_1})_n$. Note that if there exists $s>n-t$ such that $v_{j_1}\cdots v_{j_{t-2}}=v_{j_1}\cdots v_{j_h}x^s v_{j_{h+s+1}}\cdots v_{j_{t-2}}$ and $v_{j_h}=v_{j_{h+s+1}}=y$, then $\Gamma (x^{n-t}yv_{j_1}\cdots v_{j_{t-2}}y)= \Gamma (x^sv_{j_{h+s+1}}\cdots v_{j_{t-2}}yx^{n-t}yv_{j_1}\cdots v_{j_h})$. The function $predim1(n)$ gives the number of these elements which are different. We may see that
\[
\dim_k((HH_1(\tilde{A}_{a,b}))_n)=\begin{cases}
0 & \text{for $n\le 0$,}
\\
2 & \text{for $n=1$,}
\\
3 & \text{for $n=2$,}
\\
predim1(n)+4 & \text{for $n\ge 3$ odd,}
\\
predim1(n)+5 & \text{for $n\ge 4$ even.}
\end{cases}
\]

Then, the dimensions of $(HH_1(\tilde{A}_{a,b}))_n$ are for $n\ge 1$
\[
2, 3, 4, 5, 8, 14, 21, 36, 61, 107, 189, 351,\cdots 
\]
and the following are bases of $(HH_1(\tilde{A}_{a,b}))_n$ for $n\le 6$
\begin{itemize}
\item[$-$] $\{ \overline{1}\otimes x, \overline{1}\otimes y\}$ if $n=1$,
\item[$-$] $\{ \overline{x}\otimes x, \overline{y}\otimes y, \overline{x}\otimes y+\overline{y}\otimes x\}$ if $n=2$,
\item[$-$] $\{ \overline{x^2}\otimes x, \overline{y^2}\otimes y, \Gamma (x^2y), \Gamma (y^2x)\}$ if $n=3$,
\item[$-$] $\{ \overline{x^3}\otimes x, \overline{y^3}\otimes y, \Gamma (x^3y), \Gamma (y^3x), \overline{xyx}\otimes y+\overline{yxy}\otimes x\}$ if $n=4$,
\item[$-$] $\{ \overline{x^4}\otimes x, \overline{y^4}\otimes y, \Gamma (x^4y), \Gamma (y^4x), \Gamma (x^2yxy), \Gamma (xyxy^2), \Gamma (x^3y^2), \Gamma (x^2y^3)\}$ if $n=5$,
\item[$-$] $\{ \overline{x^5}\otimes x, \overline{y^5}\otimes y, \Gamma (x^5y), \Gamma (y^5x), \overline{xyxyx}\otimes y+\overline{yxyxy}\otimes x, \Gamma (x^2yx^2y), \Gamma (x^2yxy^2), \Gamma (x^2y^2xy),$ \linebreak $\Gamma (xy^2xy^2), \Gamma (x^3yxy), \Gamma (x^3y^3), \Gamma (xyxy^3), \Gamma (x^4y^2), \Gamma  (x^2y^4)\}$ if $n=6$.
\end{itemize}


\subsubsection{\texorpdfstring{Computation of $HH_2 (\tilde{A}_{a,b})$}{sec:hh2}}

Given $s\in \mathbb{N}$, $m\in \mathbb{N}_0$, $\alpha \in V^{(m)}$ and $w_1\cdots w_s\in V^{(s)}$, we denote
\[
\Gamma_s (\alpha \otimes w_1\cdots w_s)=\overline{\alpha w_1\cdots w_{s-1}}\otimes w_s+\overline{w_s\alpha w_1\cdots w_{s-2}}\otimes w_{s-1}+\cdots +\overline{w_2\cdots w_n\alpha}\otimes w_1\in A\otimes V.
\]
Note that if $m=0$, then $\Gamma_s (\alpha \otimes w_1\cdots w_s)=\alpha \Gamma (w_1\cdots w_s)$.

Since $\overline{d_3}=0$, then $HH_2 (\tilde{A}_{a,b})=\ker \overline{d_2}$.

The space $\ker \overline{d_2}$ consists of elements $\sum_j \overline{\alpha_j}\otimes x^2w_1\cdots w_{a-4}y^2+\sum_h \overline{\beta_h}\otimes x^2y^{b-4}xy$ with $\alpha_j, \beta_h \in \tilde{A}_{a,b}$, such that 
\[
\overline{d_2} (\sum_j \overline{\alpha_j} \otimes x^2w_1\cdots w_{a-4}y^2+\sum_h \overline{\beta_h} \otimes x^2y^{b-4}xy)=\sum_j \Gamma_a(\alpha_j\otimes x^2w_1\cdots w_{a-4}y^2)+\sum_h \Gamma_b(\beta_h\otimes x^2y^{b-4}xy)
\]

\noindent vanishes. It is clear that $(\ker \overline{d_2})_n=0$ for $n<a$ and that no cancellation is possible. Thus, $\overline{\alpha_j}=\overline{\beta_h}=0$ in $\tilde{A}_{a,b}$ and $HH_2 (\tilde{A}_{a,b})=0$.


\subsection{\texorpdfstring{Another example}{sec:anotherexampleofhh}}

Let $\mathcal{A}=\frac{k\langle x,y\rangle}{\langle R\rangle}$ be the $(3,4)$-Koszul algebra of \S \ref{sec:example}, where $R=\langle x^2y-2xyx+yx^2,y^2x-2yxy+xy^2,x^4,y^4\rangle$.

Its Hochschild homology is the homology of the complex
\[
\cdots \longrightarrow \mathcal{A}\otimes J_{12}^4 \overset{\overline{d_6}}{\longrightarrow} \mathcal{A}\otimes J_9^4 \overset{\overline{d_5}}{\longrightarrow} \mathcal{A}\otimes J_8^4 \overset{\overline{d_4}}{\longrightarrow} \mathcal{A}\otimes J_4^3\oplus \mathcal{A}\otimes J_5^4 \overset{\overline{d_3}}{\longrightarrow} \mathcal{A}\otimes R \overset{\overline{d_2}}{\longrightarrow} \mathcal{A}\otimes V \overset{\overline{d_1}}{\longrightarrow} \mathcal{A}\otimes k \longrightarrow 0.
\]

Given a $k$-vector space $V$, as always $\tau$ denotes the cyclic permutation in $V^{(j)}$ and $(V^{(j)})_{\tau}$ the set of coinvariants of this action. It is well-known (see for example \cite{W}) that $HH_0(\mathcal{A})= k\oplus \bigoplus_{j=1}^{\infty}(\langle x,y\rangle^{(j)})_{\tau}$.
 
Given $s\in \mathbb{N}$ and $w_1\cdots w_n\in V^{(n)}$ we denote
\[
(w_1\cdots w_n)^{+_s\tau}=\smashoperator{\sum \limits_{i=0}^{s-1}}\tau^i(w_1\cdots w_n).
\] 
 
We want to compute the Hochschild homology of $\mathcal{A}$ in higher degrees. We next list explicitely the maps $\overline{d_i}$:

\begin{itemize}
\item[$-$] $\overline{d_2}(\overline{\alpha}\otimes x^2y-2xyx+yx^2)=\overline{\alpha yx-2\alpha xy+y\alpha x+x\alpha y+xy\alpha-2yx\alpha}\otimes x+\overline{\alpha x^2-2x\alpha x+x^2\alpha}\otimes y$, 
\item[$-$] $\overline{d_2}(\overline{\beta}\otimes y^2x-2yxy+xy^2)=\overline{\beta y^2-2y\beta y+y^2\beta}\otimes x+\overline{\beta xy-2\beta yx+x\beta y+y\beta x+yx\beta-2xy\beta}\otimes y$,
\item[$-$] $\overline{d_2}(\overline{\gamma}\otimes x^4)=(\overline{\gamma x^3}+\overline{x\gamma x^2}+\overline{x^2\gamma x}+\overline{x^3\gamma})\otimes x$,
\item[$-$] $\overline{d_2}(\overline{\theta}\otimes y^4)=(\overline{\theta y^3}+\overline{y\theta y^2}+\overline{y^2\theta y}+\overline{y^3\theta})\otimes y$.
\end{itemize}

\smallskip

\begin{itemize}
\item[$-$] $\overline{d_3} (\overline{\alpha}\otimes x^2y^2-2xyxy+yx^2y+y^2x^2-2yxyx+xy^2x)=\overline{\alpha x-x\alpha}\otimes y^2x-2yxy+xy^2+\overline{\alpha y-y\alpha}\otimes x^2y-2xyx+yx^2$,
\item[$-$] $\overline{d_3} (\overline{\beta}\otimes x^5)=\overline{\beta x-x\beta}\otimes x^4$,
\item[$-$] $\overline{d_3} (\overline{\gamma}\otimes y^5)=\overline{\gamma y-y\gamma}\otimes y^4$.
\end{itemize}

For $i\ge 4$ even,
\begin{itemize}
\item[$-$] $\overline{d_i}(\overline{\alpha}\otimes x^{2i})=\overline{\alpha x^3+x\alpha x^2+x^2\alpha x+x^3\alpha}\otimes x^{2i-3}$,
\item[$-$] $\overline{d_i}(\overline{\beta}\otimes y^{2i})=\overline{\beta y^3+y\beta y^2+y^2\beta y+y^3\beta}\otimes y^{2i-3}$.
\end{itemize}

For $i\ge 5$ odd, 
\begin{itemize}
\item[$-$] $\overline{d_i}(\overline{\alpha}\otimes x^{2i-1})=\overline{\alpha x-x\alpha}\otimes x^{2i-2}$,
\item[$-$] $\overline{d_i}(\overline{\beta}\otimes y^{2i-1})=\overline{\beta y-y\beta}\otimes y^{2i-2}$.
\end{itemize}

\medskip

\begin{itemize}

\item[$\bullet$] The computation of $HH_0(\mathcal{A})$ is similar to what we have done for $HH_0(\tilde{A}_{a,b})$. Note that several generators vanish in $\mathcal{A}$, for instance $[\overline{x^n}]=[\overline{y^n}]=0$ if $n\ge 4$ and $[\overline{x^{n-1}y}]=[\overline{y^{n-1}x}]=0$ if $n\ge 5$. The function $ppredim0(n)$ in \S \ref{sec:appendix} allows us to compute $\dim_k((HH_0(\mathcal{A}))_n)$ for all $n\ge 5$.

\item[$\bullet$] The space $\ker \overline{d_1}$ is generated by $\overline{d_2}(\overline{\alpha}\otimes x^2y-2xyx+yx^2)$, $\overline{d_2}(\overline{\beta}\otimes y^2x-2yxy+xy^2)$ for $\overline{\alpha}, \overline{\beta} \in \mathcal{A}$ and by other generators obtained like in the case of $\tilde{A}_{a,b}$.

Notice that $\overline{x^{n-1}}=\overline{y^{n-1}}=0$ for $n\ge 5$, and since $\overline{v^3}\otimes v=\overline{d_2}\big( \frac{1}{4}\overline{v^3}\otimes v^4\big)$ for $v=x,y$, then $[\overline{x^3}\otimes x]=[\overline{y^3}\otimes y]=0$ in $(HH_1(\mathcal{A}))_4$.

\item[$\bullet$] It is easy to see that
\begin{itemize}
\item[$-$] $(\ker \overline{d_2})_n=0$ for $n\le 3$,
\item[$-$] $(\ker \overline{d_2})_4=\langle \overline{x}\otimes x^2y-2xyx+yx^2, \overline{y}\otimes x^2y-2xyx+yx^2, \overline{x}\otimes y^2x-2yxy+xy^2, \overline{y}\otimes y^2x-2yxy+xy^2 \rangle$,
\item[$-$] $(\ker \overline{d_2})_5=\langle \overline{xy-yx}\otimes x^2y-2xyx+yx^2, \overline{xy-yx}\otimes x^2y-2xyx+yx^2, \overline{x}\otimes x^4, \overline{y}\otimes y^4\rangle$ (note that $\overline{xy-yx}\otimes x^2y-2xyx+yx^2, \overline{xy-yx}\otimes x^2y-2xyx+yx^2 \in \im \overline{d_3}$),
\item[$-$] $(\ker \overline{d_2})_6=(\im \overline{d_3})_6+\langle \overline{x^2}\otimes x^4, \overline{y^2}\otimes y^4 \rangle$,
\item[$-$] $(\ker \overline{d_2})_7=(\im \overline{d_3})_7+\langle \overline{x^3}\otimes x^4, \overline{y^3}\otimes y^4 \rangle$,
\item[$-$] $(\ker \overline{d_2})_8=(\im \overline{d_3})_8$,
\item[$-$] $(\ker \overline{d_2})_n=(\im \overline{d_3})_n+\langle \overline{xv^1_{j_1}\cdots v^1_{j_{n-8}}x^3}\otimes x^4, \overline{x^2v^2_{j_1}\cdots v^2_{j_{n-8}}x^2}\otimes x^4, \overline{x^3v^3_{j_1}\cdots v^3_{j_{n-8}}x}\otimes x^4, \overline{yv^1_{h_1}\cdots v^1_{h_{n-8}}y^3}\otimes y^4, \overline{y^2v^2_{h_1}\cdots v^2_{h_{n-8}}y^2}\otimes y^4, \overline{y^3v^3_{h_1}\cdots v^3_{h_{n-8}}y}\otimes y^4  \rangle$ for $n\ge 9$ and $v_{j_l}^s,v_{h_l}^s\in \{x,y\}$.
\end{itemize}

\item[$\bullet$] Notice that $[\overline{xv^1_{j_1}\cdots v^1_{j_{n-8}}x^3}\otimes x^4]=[\overline{x^2v^2_{j_1}\cdots v^2_{j_{n-8}}x^2}\otimes x^4]= [\overline{x^3v^3_{j_1}\cdots v^3_{j_{n-8}}x}\otimes x^4]$ and $[\overline{yv^1_{h_1}\cdots v^1_{h_{n-8}}y^3}\otimes y^4]=[\overline{y^2v^2_{h_1}\cdots v^2_{h_{n-8}}y^2}\otimes y^4]=[\overline{y^3v^3_{h_1}\cdots v^3_{h_{n-8}}y}\otimes y^4]$ in $HH_2(\mathcal{A})$.

\item[$\bullet$] By direct computations, 
\begin{itemize}
\item[$-$] $(\ker \overline{d_3})_n=0$ for $n\le 3$,
\item[$-$] $(\ker \overline{d_3})_4=\langle \overline{1}\otimes x^2y^2-2xyxy+yx^2y+y^2x^2-2yxyx+xy^2x\rangle$,
\item[$-$] $(\ker \overline{d_3})_5=\langle \overline{yx^3+xyx^2+x^2yx+x^3y}\otimes x^5, \overline{xy^3+yxy^2+y^2xy+y^3x}\otimes y^5, \overline{1}\otimes x^5, \overline{1}\otimes y^5\rangle$ (note that $\overline{yx^3+xyx^2+x^2yx+x^3y}\otimes x^5, \overline{xy^3+yxy^2+y^2xy+y^3x}\otimes y^5\in \im \overline{d_4}$),
\item[$-$] $(\ker \overline{d_3})_6=(\im \overline{d_4})_6+\langle \overline{x}\otimes x^5, \overline{y}\otimes y^5\rangle$,
\item[$-$] $(\ker \overline{d_3})_7=(\im \overline{d_4})_7+\langle \overline{x^2}\otimes x^5, \overline{y^2}\otimes y^5\rangle$,
\item[$-$] $(\ker \overline{d_3})_8=(\im \overline{d_4})_8+\langle \overline{x^3}\otimes x^5, \overline{y^3}\otimes y^5\rangle$ (note that $\overline{x^3}\otimes x^5, \overline{y^3}\otimes y^5\in \im \overline{d_4}$),
\item[$-$] $(\ker \overline{d_3})_n=(\im \overline{d_4})_n$ for $n=9,10,11$,
\item[$-$] $(\ker \overline{d_3})_n=\langle \overline{x^3w_{j_1}\cdots w_{j_{n-11}}x^3}\otimes x^5, \overline{y^3w_{h_1}\cdots w_{h_{n-11}}y^3}\otimes y^5,
\overline{(v_{j_1}\cdots v_{j_{n-9}}yx^3)^{+_4\tau}}\otimes x^5,$ $\overline{(v_{h_1}\cdots v_{h_{n-9}}xy^3)^{+_4\tau}}\otimes y^5\rangle$ for $n\ge 12$ and $w_{j_l},w_{h_l},v_{j_l},v_{h_l}\in \{ x,y\}$ (it is clear that $\overline{(v_{j_1}\cdots v_{j_{n-9}}yx^3)^{+_4\tau}}\otimes x^5, \overline{(v_{h_1}\cdots v_{h_{n-9}}xy^3)^{+_4\tau}}\otimes y^5\in \im \overline{d_4}$).
\end{itemize}

\item[$\bullet$] Given $i\ge 4$ even, 
\begin{itemize}
\item[$-$] $(\ker \overline{d_i})_n=0$ for $n\le 2i$,
\item[$-$] $(\ker \overline{d_i})_{2i+1}=\langle \overline{x}\otimes x^{2i}, \overline{y}\otimes y^{2i}\rangle$,
\item[$-$] $(\ker \overline{d_i})_{2i+2}=\langle \overline{x^2}\otimes x^{2i}, \overline{y^2}\otimes y^{2i}, \overline{yx-xy}\otimes x^{2i}, \overline{xy-yx}\otimes y^{2i} \rangle$,
\item[$-$] $(\ker \overline{d_i})_{2i+3}=\langle \overline{x^3}\otimes x^{2i}, \overline{y^3}\otimes y^{2i}, \overline{\alpha x-x\alpha}\otimes x^{2i}, \overline{\beta y-y\beta}\otimes y^{2i} \rangle$ for $\alpha ,\beta \in V^{(2)}$,
\item[$-$] $(\ker \overline{d_i})_{2i+4}=\langle \overline{\alpha x-x\alpha}\otimes x^{2i}, \overline{\beta y-y\beta}\otimes y^{2i} \rangle$ for $\alpha ,\beta \in V^{(3)}$,
\item[$-$] $(\ker \overline{d_i})_n=\langle \overline{xv^1_{j_1}\cdots v^1_{j_{n-2i-4}}x^3}\otimes x^{2i}, \overline{x^2v^2_{j_1}\cdots v^2_{j_{n-2i-4}}x^2}\otimes x^{2i}, \overline{x^3v^3_{j_1}\cdots v^3_{j_{n-2i-4}}x}\otimes x^{2i},$ $\overline{yv^1_{h_1}\cdots v^1_{h_{n-2i-4}}y^3}\otimes y^{2i}, \overline{y^2v^2_{h_1}\cdots v^2_{h_{n-2i-4}}y^2}\otimes y^{2i}, \overline{y^3v^3_{h_1}\cdots v^3_{h_{n-2i-4}}y}\otimes y^{2i}, \overline{\alpha x-x\alpha}\otimes x^{2i},$ $\overline{\beta y-y\beta}\otimes y^{2i} \rangle$ for $n\ge 2i+5$, $v_{j_l}^s,v_{h_l}^s\in \{x,y\}$ and $\alpha ,\beta \in V^{(n-2i-1)}$.
\end{itemize}

\item[$\bullet$] It is clear that $\overline{\alpha x-x\alpha}\otimes x^{2i}, \overline{\beta y-y\beta}\otimes y^{2i}\in \ker \overline{d_i}\cap \im \overline{d_{i+1}}$ for $i$ even.

Notice that $\overline{xv^1_{j_1}\cdots v^1_{j_{n-2i-4}}x^3}\otimes x^{2i}$, $\overline{x^2v^2_{j_1}\cdots v^2_{j_{n-2i-4}}x^2}\otimes x^{2i}$, $\overline{x^3v^3_{j_1}\cdots v^3_{j_{n-2i-4}}x}\otimes x^{2i}$ belong to the same class in $\frac{(\ker \overline{d_i})_n}{(\im \overline{d_{i+1}})_n}$ for $i$ even and $n\ge 2i+5$. Taking $y$ instead of $x$ is analogous.

\item[$\bullet$] Given $i\ge 5$ odd, 
\begin{itemize}
\item[$-$] $(\ker \overline{d_i})_n=0$ for $n\le 2i-1$,
\item[$-$] $(\ker \overline{d_i})_{2i-1}=\langle \overline{1}\otimes x^{2i-1}, \overline{1}\otimes y^{2i-1}\rangle$,
\item[$-$] $(\ker \overline{d_i})_{2i}=\langle \overline{x}\otimes x^{2i-1}, \overline{y}\otimes y^{2i-1}\rangle$,
\item[$-$] $(\ker \overline{d_i})_{2i+1}=\langle \overline{x^2}\otimes x^{2i-1}, \overline{y^2}\otimes y^{2i-1}\rangle$,
\item[$-$] $(\ker \overline{d_i})_{2i+2}=\langle \overline{x^3}\otimes x^{2i-1}, \overline{y^3}\otimes y^{2i-1}\rangle$,
\item[$-$] $(\ker \overline{d_i})_{2i+3}=\langle \overline{(x^3y)^{+_4\tau}}\otimes x^{2i-1}, \overline{(y^3x)^{+_4\tau}}\otimes y^{2i-1}\rangle$,
\item[$-$] $(\ker \overline{d_i})_n=\langle \overline{(v_{j_1}\cdots v_{j_{n-2i-3}}yx^3)^{+_4\tau}}\otimes x^{2i-1}, \overline{(v_{h_1}\cdots v_{h_{n-2i-3}}xy^3)^{+_4\tau}}\otimes y^{2i-1}\rangle$ for $n=2i+4,2i+5$ and $v_{j_l},v_{h_l}\in \{ x,y\}$,
\item[$-$] $(\ker \overline{d_i})_n=\langle \overline{(v_{j_1}\cdots v_{j_{n-2i-3}}yx^3)^{+_4\tau}}\otimes x^{2i-1}, \overline{(v_{h_1}\cdots v_{h_{n-2i-3}}xy^3)^{+_4\tau}}\otimes y^{2i-1}, \overline{x^3w_{j_1}\cdots}$ $\overline{w_{j_{n-2i-5}}x^3}\otimes x^{2i-1},\overline{y^3w_{h_1}\cdots w_{h_{n-2i-5}}y^3}\otimes y^{2i-1}\rangle$ for $n\ge 2i+6$ and $w_{j_l},w_{h_l},v_{j_l},v_{h_l}\in \{ x,y\}$.
\end{itemize}

\item[$\bullet$] Notice that if $i$ is odd $(\ker \overline{d_i})_n=(\im \overline{d_{i+1}})_n$ for $2i+2\le n\le 2i+5$. Moreover, if $n\ge 2i+6$ then $\overline{(v_{j_1}\cdots v_{j_{n-2i-3}}yx^3)^{+_4\tau}}\otimes x^{2i-1}, \overline{(v_{h_1}\cdots v_{h_{n-2i-3}}xy^3)^{+_4\tau}}\otimes y^{2i-1}\in \im \overline{d_{i+1}})_n$.
\end{itemize}

We conclude that:

\begin{itemize}
\item[$\bullet$] $\dim_k((HH_0(\mathcal{A}))_n)=\begin{cases}
0 & \text{for $n<0$},
\\
1 & \text{for $n=0$},
\\
2 & \text{for $n=1$},
\\
3 & \text{for $n=2$},
\\
4 & \text{for $n=3,4$},
\\
ppredim0(n) & \text{for $n\ge 5$ odd},
\\
ppredim0(n)+1 & \text{for $n\ge 6$ even}.
\end{cases}$

\item[$\bullet$] $\dim_k((HH_1(\mathcal{A}))_n)=\begin{cases}
0 & \text{for $n\le 0$},
\\
2 & \text{for $n=1$},
\\
3 & \text{for $n=2$},
\\
predim1(n)+2 & \text{for $n\ge 3$ odd},
\\
predim1(n)+3 & \text{for $n\ge 4$ even}.
\end{cases}$

\item[$\bullet$] $\dim_k((HH_2(\mathcal{A}))_n)=\begin{cases}
0 & \text{for $n\le 3$ or $n=8$},
\\
4 & \text{for $n=4$},
\\
2 & \text{for $n=5,6,7$},
\\
2predim2(n-4) & \text{for $n\ge 9$},
\end{cases}$

\noindent where $predim2$ may be computed as in \S \ref{sec:appendix}. 

For example, for $5\le n\le 13$, $predim2(n)$ is respectively $1,2,4,6,12,22,41,74,137$.

\item[$\bullet$] $\dim_k((HH_3(\mathcal{A}))_n)=\begin{cases}
0 & \text{for $n\le 3$ or $8\le n\le 11$},
\\
1 & \text{for $n=4$},
\\
2 & \text{for $n=5,6,7$},
\\
2predim3(n-5) & \text{for $n\ge 12$}.
\end{cases}$

The function $predim3$ is computed as described in \S \ref{sec:appendix}, for $7\le n\le 15$, $predim3(n)$ is respectively $1,1,2,3,7,12,22,40,75$.

\item[$\bullet$] For $i\ge 4$ even,
\[
\dim_k((HH_i(\mathcal{A}))_n)=\begin{cases}
0 & \text{for $n\le 2i$ or $n=2i+4$},
\\
2 & \text{for $n=2i+1,2i+2,2i+3$},
\\
2predim2(n-2i) & \text{for } n\ge 2i+2.
\end{cases}
\]

\item[$\bullet$] For $i\ge 5$ odd,
\[
\dim_k((HH_i(\mathcal{A}))_n)=\begin{cases}
0 & \text{for $n\le 2i-2$ or $2i+2\le n\le 2i+5$},
\\
2 & \text{for $n=2i-1,2i,2i+1$},
\\
2predim3(n-2i+1) & \text{for $n\ge 2i+6$}.
\end{cases}
\]
\end{itemize}


\section{\texorpdfstring{Appendix}{sec:appendix}}\label{sec:appendix}

In this section we give a program using \textit{Mathematica} to compute the Hochschild homology of the algebras $\mathcal{A}$ and $\tilde{A}_{a,b}$ for $(a,b)=(4,5)$ or $6\le a<b$.

\begin{lstlisting}

prod[a_, bs_, c_] :=  Map[Join[a, #, c] &, bs];

normalform[w_] := 
  Sort[Most[NestList[RotateRight, w, Length[w]]]][[1]];

delete[s_] := 
  Union[Map[normalform, s], Map[normalform, (s /. {x -> y, y -> x})]];

relb[b_] := Join[{x, x}, Table[y, {b - 4}], {x, y}];

ker[n_, t_] := Module[{d, xt, step1, step2},
   d = n - t - 1;
   step1 = Map[
     Function[is, 
      Flatten[
       MapIndexed[
        If[EvenQ[#2[[1]]], Table[x, {#1}], Table[y, {#1}]] &, is]]
      ],
     Select[
      Flatten[Map[Permutations, IntegerPartitions[d, d, Range[1, t]]],
        1], #[[-1]] < t || EvenQ[Length[#]] &]
     ];
   xt = Table[x, { t}];
   step2 = prod[xt, step1, {y}];
   delete[step2]
   ];

generators[n_] := Flatten[Table[ker[n, t], {t, 2, n - 2}], 1];   

contain[w_, v_] := 
  0 != Length[Select[
     NestList[RotateRight, w, Length[v] - 1],
     MatchQ[#, 
       Append[Prepend[v, BlankNullSequence[]], BlankNullSequence[]]] &
     ]
    ];
relations[a_] := 
  MatchQ[a, {___, x, x, x, x, ___}]  || 
   MatchQ[a, {___, y, y, y, y, ___}];

ppredim0[n_] := 
  Module[{counter = 0, x4 = Table[x, {4}], y4 = Table[y, {4}]},
   Do[
    If[! (contain[p, x4] || contain[p, y4]), Print[p]],
    {p, generators[n]}
    ];
   Do[
    If[! (contain[p, x4] || contain[p, y4]), counter = counter + 1],
    {p, generators[n]}
    ];
   Print["ppredim0(", n, ")=", counter]
   ];
  
predim0[n_, a_, b_] := 
  Module[{aword, bword, vs, auxvs, auxrelation, counter = 0},
   auxvs = Map[
     Function[is, 
      Flatten[
       MapIndexed[
        If[OddQ[#2[[1]]], Table[x, {#1}], Table[y, {#1}]] &, is]]
      ],
     IntegerPartitions[a - 4, {a - 4}]
     ];
   vs = Flatten[auxvs];
   auxrelation = If[EvenQ[Length[vs]], RotateRight[vs], vs];
   aword = Join[{x, x}, auxrelation, {y, y}];
   bword = relb[b];
   Do[
    If[! (contain[p, aword] || contain[p, bword]), Print[p]],
    {p, generators[n]}
    ];
   Do[
    If[! (contain[p, aword] || contain[p, bword]), 
     counter = counter + 1],
    {p, generators[n]}
    ];
   Print["predim0(", n, ",", a, ",", b, ")=", counter];
   ];

predim1[n_] := Print["predim1(", n,  ")=", Length[generators[n]]];   

predim2[n_] := Module[{generators},
   generators = Map[Join[{x}, #, {x, x, x}] &, Tuples[{x, y}, n - 4]];
   Print["predim2(", n, ")=", 
    2^(n - 4) - Length[Select[generators, relations]]]
   ];
   
predim3[n_] := Module[{generators},
   generators = 
    Map[Join[{x, x, x}, #, {x, x, x}] &, Tuples[{x, y}, n - 6]];
   Print["predim3(", n, ")=", 
    2^(n - 6) - Length[Select[generators, relations]]]
   ];

\end{lstlisting}


\section*{\texorpdfstring{Acknowledgements}{sec:acknowledgements}}
\addcontentsline{toc}{section}{Acknowledgements}

The authors thank Eduardo Marcos and Mariano Su\'arez-\'Alvarez for interesting comments and suggestions.

The authors were partially supported by UBACYT X212, PIP-CONICET 112-200801-00487, PICT 2007-02182.



\bigskip

\noindent Dto. de Matem\'atica, Facultad de Ciencias Exactas y Naturales

\noindent Universidad de Buenos Aires

\noindent Ciudad Universitaria, Pabell\'on I

\noindent (1428) Buenos Aires, Argentina

\noindent e-mail: aarey@dm.uba.ar - asolotar@dm.uba.ar


\begin{thebibliography}{9999999}
\addcontentsline{toc}{section}{References}

\bibitem[A]{A} Aigner, M. \textit{Combinatorial Theory}. Reprint of the 1979 original. Classics in Mathematics. Springer-Verlag, Berlin, 1997.

\bibitem[B1]{B1} Benkart, G. \textit{Down-up algebras and Witten's deformations of the universal enveloping algebra of $\mathfrak{sl}_2$}. Contemporary Matematics, \textbf{224} (1999), 24--45, Amer. Math. Soc., Providence, RI, 1999.

\bibitem[B2]{B2} Berger, R. \textit{Koszulity for nonquadratic algebras}. J. Algebra \textbf{239} (2001), 705--734.

\bibitem[B3]{B3} Berger, R. \textit{Koszulity for nonquadratic algebras II}.

arXiv:math/0301172V1 [math.QA]

\bibitem[B4]{B4} Berger, R. \textit{La cat\'egorie des modules gradu\'es sur une alg\`ebre gradu\'ee (nouvelle version du chapitre 5 d'un cours de Master 2 \'a Lyon 1)} (2008).

http://webperso.univ-st-etienne.fr/$\sim$rberger/mes-textes.html

\bibitem[B5]{B5} Bongartz, K. \textit{Algebras and quadratic forms}. J. London Math. Soc. \textbf{28} (2) (1983), 461--469.

\bibitem[BBK]{BBK} Brenner, S.; Butler, M.; King, A. \textit{Periodic algebras which are almost Koszul}. Algebr. Represent. Theory \textbf{5} (2002), 331--367.

\bibitem[BDW]{BDW} Berger, R.; Dubois-Violette, M.; Wambst, M. \textit{Homogeneous algebras}. J. Algebra \textbf{261} (2003), 172–-185.

\bibitem[BF]{BF} Backelin, J.; Fr\"oberg, R. \textit{Koszul algebras, Veronese subrings and rings with linear resolutions}. Rev. Roumaine Math. Pures Appl. \textbf{30} (1985), 85--97.

\bibitem[BG]{BG} Berger, R.; Ginzburg, V. \textit{Higher symplectic reflection algebras and non-homogeneous {$N$}-{K}oszul property}. J. Algebra. \textbf{304} (2006), 577--601.

\bibitem[BGS1]{BGS1} Beilinson, A.; Ginzburg, V.; Schechman, V. \textit{Koszul duality}. J. Geom. Phys. \textbf{5} (1988), 317--350.

\bibitem[BGS2]{BGS2} Beilinson, A.; Ginzburg, V.; Soergel, W. \textit{Koszul duality patterns in representation theory}. J. Amer. Math. Soc. \textbf{9} (1996), no. 2, 473--527.

\bibitem[C]{C} Cartan, H. \textit{Homologie et cohomologie d'une alg\`{e}bre gradu\'ee}. S\'eminaire Cartan, Paris, 1958/59, expos\'e 15.

\bibitem[CS]{CS} Cassidy, T.; Shelton, B. \textit{Generalizing the notion of Koszul algebra}. Math. Z. \textbf{260} (2008), no. 1, 93--114.

\bibitem[F] {F} Fr\"oberg, R. \textit {Koszul Algebras}, In: {\em Advances in Commutative Ring Theory}. Proceedings of the 3rd International Conference, Fez, Lect. Notes Pure Appl. Math. \textbf{205}, Marcel Dekker, New York, (1999), 337--350.

\bibitem[FV]{FV} Fl{\o}ystad, G.; Vatne, J.E. \textit{P{BW}-deformations of $N$-Koszul algebras}, J. Algebra \textbf{302} (2006), 116--155.



\bibitem[GMMZ]{GMMZ} Green, E.; Marcos, E.; Mart\'{\i}nez-Villa, R.; Zhang, P. \textit{$D$-Koszul algebras}. J. Pure Appl. Algebra \textbf{193} (2004), no. 1-3, 141--162.

\bibitem[GMV]{GMV} Green, E.; Mart\'{\i}nez-Villa, R. \textit{Koszul and Yoneda algebras}. Representation theory of algebras (Cocoyoc, 1994), 247--297, CMS Conf. Proc., \textbf{18}, Amer. Math. Soc., Providence, RI, 1996.

\bibitem[HL]{HL} Hai, P.; Lorenz, M. \textit{Koszul algebras and the quantum MacMahon master theorem}. Bull. Lond. Math. Soc. \textbf{39} (2007), no. 4, 667--676.


\bibitem[J]{J} J\'onsson, B. \textit{Distributive sublattices of a modular lattice}. Proc. Amer. Math. Soc. \textbf{6} (1955), 682--688.

\bibitem[K]{K} Koszul, J.L. \textit{Homologie et cohomologie des alg\`{e}bres de Lie}. Bulletin de la S. M. F. \textbf{78} (1950), 65--127. 

\bibitem[M]{M} Manin, Y. \textit{Some remaks on Koszul algebras and quantum groups}. Ann. Inst. Fourier \textbf{37} (1987), 191--205.

\bibitem[O]{O} \"Ore, O. \textit{On the foundation of abstract algebra I}. Ann. of Math. (2) \textbf{36} (1935), 406--437.

\bibitem[P]{P} Priddy, S. \textit{Koszul resolutions}. Trans. Amer. Math. Soc. \textbf{152} (1970) , 39--60.

\bibitem[R]{R} Rey, A. \textit{\'Algebras $(a,b)$-Koszul y \'algebras $(a,b)$-cuasi Koszul}. Tesis doctoral, UBA, 2009.

http://cms.dm.uba.ar/academico/carreras/doctorado/desde

\bibitem[W]{W} Weibel, C. \textit{An introduction to homological algebra}. Cambridge Studies in Advanced Mathematics, \textbf{38}. Cambridge Univ. Press, Cambridge, 1994.

\bibitem[WM]{WM} Wolfram Research, Inc., Mathematica, Version 7.0, Champaign, IL (2008).

\end{thebibliography}
\end{document}